\documentclass[letterpaper,11pt,reqno]{amsart}
\usepackage[margin=1in]{geometry}
\usepackage[T1]{fontenc}
\usepackage{microtype}

\usepackage{mathrsfs,xfrac} 
\usepackage[colorlinks=true,linkcolor=blue,citecolor=blue,urlcolor=blue]{hyperref} 
\usepackage{amsmath,amssymb,amsthm,amsfonts,amsbsy,latexsym,dsfont,color,graphicx,enumitem}
\usepackage[foot]{amsaddr}
\usepackage{caption}
\usepackage{regexpatch}
\usepackage{comment}
\usepackage{todonotes}
\makeatletter
\xpatchcmd{\@todo}{\setkeys{todonotes}{#1}}{\setkeys{todonotes}{inline,#1}}{}{}
\makeatother

\newenvironment{enumeratei}{\begin{enumerate}[label=\textup{\roman*.},leftmargin=2em,itemsep=1ex]}{\end{enumerate}}
\newenvironment{enumeratea}{\begin{enumerate}[label=\textup{(\alph*)},leftmargin=2em,itemsep=1ex]}{\end{enumerate}}

\theoremstyle{remark}

\newtheorem{remark}{Remark}[section]

\renewcommand{\le}{\leqslant}  
\renewcommand{\ge}{\geqslant}

\newcommand{\ind}{\mathds{1}}
\newcommand{\eps}{\varepsilon}

\newcommand{\norm}[1]{\left\Vert#1\right\Vert}
\newcommand{\abs}[1]{\left\vert#1\right\vert}

\def\qed{ \hfill $\blacksquare$}  
\let\ga=\alpha   \let\gd=\delta 
\let\gf=\varphi    \let\gk=\kappa \let\gl=\lambda      \let\go=\omega   \let\gs=\sigma

\newcommand{\cB}{\mathcal{B}}
\newcommand{\cE}{\mathcal{E}}\newcommand{\cF}{\mathcal{F}}
\newcommand{\cG}{\mathcal{G}}\newcommand{\cH}{\mathcal{H}}
\newcommand{\cK}{\mathcal{K}}\newcommand{\cL}{\mathcal{L}}
\newcommand{\cN}{\mathcal{N}}
\newcommand{\cR}{\mathcal{R}}
\newcommand{\cT}{\mathcal{T}}\newcommand{\cU}{\mathcal{U}}
\newcommand{\cV}{\mathcal{V}}
\newcommand{\cZ}{\mathcal{Z}}  

\newcommand{\mvX}{\boldsymbol{X}}

\newcommand{\dR}{\mathds{R}}
\newcommand{\dT}{\mathds{T}}

\newcommand{\dZ}{\mathds{Z}} 

\DeclareMathOperator{\E}{\mathds{E}}
\DeclareMathOperator{\pr}{\mathds{P}}

\DeclareMathOperator{\var}{Var}

\theoremstyle{plain}
\newtheorem{theorem}{Theorem}[section]
\newtheorem{proposition}[theorem]{Proposition}
\newtheorem{lemma}[theorem]{Lemma}
\newtheorem{corollary}[theorem]{Corollary}
\theoremstyle{definition}
\newtheorem{assumption}[theorem]{Assumption}

\newcommand{\dd}{\mathrm{d}}
\newcommand{\wto}{\Rightarrow}
\newcommand{\pto}{\xrightarrow{ \pr }}
\newcommand{\eqd}{\overset{d}{=}}
\newcommand{\fa}{\mathfrak{a}}
\newcommand{\Exp}{\mathrm{Exp}}
\newcommand{\Tor}{\dT_n^d}
\newcommand{\tcol}{\tau_{\mathrm{col}}}
\newcommand{\Mit}{\mathsf{E}}
\newcommand{\dtr}{d_{\mathrm{tr}}}
\newcommand{\tW}{\widetilde W}
\newcommand{\tT}{\widetilde T}
\renewcommand{\le}{\leqslant}
\renewcommand{\ge}{\geqslant}
\DeclareMathOperator{\gen}{gen}
\DeclareMathOperator{\normal}{N}

\title[Mean-field Fluctuations and Geodesics]{Mean-field Fluctuations for First-Passage Percolation with Symmetric Kernels}
\author[Dey]{Partha S.~Dey$^\star$}
\author[Kim]{Daecheol Kim$^\dagger$}
\address{Department of Mathematics, University of Illinois Urbana-Champaign, Urbana, Illinois 61801}
\email{\{$^\star$psdey,$^\dagger$dk43\}@illinois.edu}
\date{\today}

\begin{document}
\begin{abstract}
We study first-passage percolation on finite weighted graphs whose normalized Weibull edge hazards form a symmetric Markov kernel on $N$ vertices. We prove that a uniform absolute spectral gap and an $O(N^{-1})$ bound on a fixed even-step return probability imply the stochastic mean-field passage-time limit and an independent Gaussian hopcount limit. Jointly with these quantities, the entire geodesic from an arbitrary source to an independent uniform target is approximated in total variation by an independent Markov chain stopped at the hopcount. This comparison yields laws of large numbers for fixed-length path segments and joint central limit theorems for bounded additive observables. The results extend to edge laws with a uniform relative hazard error $o(1/\log N)$, including common noise distributions with a power-law hazard and a positive-power error near zero. Applications include integrable singular torus kernels, macroscopic spread-out graphs, non-transitive regular graphs, finite groups, and discretizations of symmetric continuum kernels. 
\end{abstract}
\maketitle

\section{Introduction}\label{sec:intro}
First-passage percolation (FPP)  was  introduced by Hammersley and Welsh~\cite{HW65} to model the spread of a fluid or infection through a random medium. It assigns independent nonnegative lengths to graph edges and studies the resulting shortest-path metric. Under the usual shape-theorem assumptions, passage times on nearest-neighbor $\dZ^d$ grow linearly with distance, while their fluctuations remain open in general~\cite{Kes86,ADH17}. On the complete graph $K_N$ with exponential weights, Janson~\cite{Jan99} identified the typical-distance scale $\log N/N$ and the centered limit $G_1+G_2-G_3$ with independent standard Gumbel variables. Bhamidi and van der Hofstad~\cite{BvdH12} obtained the corresponding joint passage-time and Gaussian hopcount limits for Weibull weights of every fixed shape $\theta>0$.

Long-range models interpolate between lattice and mean-field growth, capturing spatial epidemics~\cite{Ric73,Mol77}, small-world networks~\cite{WS98,Kle00}, and spatial networks~\cite{Bar11} in which distant contacts occur more slowly.
Chatterjee and Dey~\cite{CD16} considered exponential weights $\|x-y\|^\ga\go_{x,y}$ on all vertex-pairs in $\dZ^d$ and found four growth regimes. In the instantaneous regime $\ga<d$, passage times on the infinite lattice vanish, so nontrivial asymptotics require a finite volume. On the torus $\Tor=(\dZ/n\dZ)^d$, with $N=n^d$, van der Hofstad and Lodewijks~\cite{vdHL24} proved first-order typical-distance, flooding-time and diameter asymptotics and asked about fluctuations, hopcounts and other weight laws. In a parallel direction, Dey and Kim~\cite{DK26b} proved extreme-value fluctuations for Weibull FPP on macroscopic spread-out cycles and identified a crossover at range $\ell\asymp n/\log n$.
These results motivate the following question --
    \emph{which finite spatial interaction kernels share the typical-distance and 
    hopcount limits of the stochastic mean-field model, and
    how does the underlying geometry shape the geodesic path?}

In this article, we give sufficient, geometry-independent criteria for the mean-field typical-distance and hopcount limits, and identify how spatial structure survives in the geodesic. The normalized edge intensities form a finite symmetric Markov kernel $P_N$. A uniform absolute spectral gap and an $O(N^{-1})$ fixed even-step return bound suffice. These permit non-transitive graphs, integrable singular interactions and graph sequences without a limiting graphon or metric space. The full geodesic, jointly with its passage time and hopcount, is approximated in total variation by an unconditioned $P_N$-Markov chain of that random length.
For torus costs $\norm{x-y}^\ga\go_{x,y}$ with Weibull shape $\theta$, this gives typical-distance fluctuations, a hopcount central limit theorem and a path law throughout $\ga\theta<d$ for a uniformly sampled target.

\subsection{The model}\label{sec:model}\label{sec:base-graphs}

Consider a graph $G_N$ with vertex set $\cV_N=[N]$ and edge set $\cE_N$. Let $f_N:\cV_N^2\to(0,\infty]$ be symmetric, finite and positive precisely on the allowed edges $\cE_N$, and infinite on the diagonal and on missing edges. Let $(\go_e)_{e\in\cE_N}$ be i.i.d.~nonnegative random weights. We define
\begin{align}\label{eq:weights}
W_{\{x,y\}}:=f_N(x,y)\cdot \omega_{\{x,y\}}, \quad \{x,y\}\in\cE_N.
\end{align}
We write $W_{x,y}=W_{y,x}:=W_{\{x,y\}}$ for $\{x,y\}\in\cE_N$.
For $u,v\in\cV_N$, define the first-passage time
\begin{align}\label{eq:fpt}
T_N(u,v):=\inf_{\pi:u\to v}\sum_{e\in\pi}W_e,
\end{align}
where the infimum is taken over all self-avoiding paths in $G_N$ from $u$ to $v$.
First, we assume the following condition on the edge weight distribution, which will be relaxed later.

\begin{assumption}\label{ass:weights}
The variables $(\go_e)_{e\in\cE_N}$ are independent, with
$\pr(\go_e>t)=e^{-t^\theta}$ for $t\ge0$ and a fixed $\theta>0$.
\end{assumption}

With the convention $\infty^{-\theta}=0$, define the edge intensities
and row sums
\begin{align*}
g_N(x,y):=f_N(x,y)^{-\theta},\quad
S_N(x):=\sum_{y\ne x}g_N(x,y),
\quad x,y\in\cV_N.
\end{align*}

\begin{assumption}\label{ass:SN}
There exists a constant $S_N>0$ such that $S_N(x)=S_N$ for every
$x\in\cV_N$.
\end{assumption}

Under Assumption~\ref{ass:SN}, setting
$$
P_N(x,y):={g_N(x,y)}/{S_N},\quad x,y\in\cV_N,
$$
defines a symmetric stochastic matrix satisfying
\begin{align*}
P_N(x,y)=P_N(y,x)\ge0,\quad
P_N(x,x)=0,\quad
\sum_{y\ne x}P_N(x,y)=1,\quad
\cE_N=\{\{x,y\}:P_N(x,y)>0\}.
\end{align*}
Define the Malthusian growth parameter
\begin{align}\label{eq:malthusian}
\fa_N:=\bigl(\Gamma(\theta+1)S_N\bigr)^{1/\theta}.
\end{align}
For $x\ne y$, the rescaled edge weights $\fa_NW_{x,y}$ satisfy
\begin{align}\label{eq:rescaled-law}
\pr(\fa_NW_{x,y}>s)=
\exp\biggl(-\frac{s^\theta}{\Gamma(\theta+1)}P_N(x,y)\biggr),\quad s\ge0.
\end{align}
Thus the law of $\fa_NT_N$ depends only on $P_N$ and $\theta$, and is invariant under multiplying $g_N$ and $S_N$ by a common constant. Our spatial assumptions are as follows.

\begin{assumption}\label{ass:graphs}
There exist an integer $r\ge1$ and constants $\eta\in(0,1)$ and $C<\infty$, independent of $N$, such that, for all large $N$, we have
\begin{align}\label{eq:graph-assumptions}
\|P_N-\Pi_N\|_{2\to2}\le1-\eta,\text{ and }
\max_{x\in\cV_N}P_N^{2r}(x,x)\le\frac{C}{N},
\end{align}
where $\Pi_N(x,y):=1/N$ and the operator norm is on $\ell^2(\cV_N)$.
\end{assumption}

The absolute spectral gap implies that the graph of allowed edges is connected for large $N$, and continuity of the weights yields almost sure uniqueness of geodesics. By symmetry, $P_N^{2r}(x,x)=\sum_yP_N^r(x,y)^2$, so the return bound states that in $r$ steps the walk spreads over order $N$ vertices in $L^2$. By H\"older's inequality, it implies
\begin{align}\label{eq:intro-small-masses}
\max_xP_N^2(x,x)\le(C/N)^{1/r},\text{ and }
\max_{x,y}P_N(x,y)\le(C/N)^{1/(2r)}.
\end{align}
Thus the effective minimum degree $1/\max_{x,y}P_N(x,y)$ grows at least like $N^{1/(2r)}$.

\subsection{Main Results}\label{sec:results}

Let $M_\infty$ be the rescaled population limit of the Crump--Mode--Jagers (CMJ) process with reproduction intensity $\mu_\theta(\dd s)=s^{\theta-1}\dd s/\Gamma(\theta)$ ($s>0$), defined in Section~\ref{sec:cmj-prelim}. It satisfies $M_\infty>0$ almost surely and $\E M_\infty=1/\theta$. Let $M_\infty^{(1)},M_\infty^{(2)}$ be independent copies of $M_\infty$, and let $G$ be an independent standard Gumbel variable with $\pr(G\le x)=\exp(-e^{-x})$. Define
\begin{align}\label{eq:limit-law}
Y_\theta:=-\log M_\infty^{(1)}-\log M_\infty^{(2)}-G-\log\theta.
\end{align}
For $\theta=1$, we have $M_\infty\sim\Exp(1)$ and $Y_1\eqd G_1+G_2-G_3$ for independent standard Gumbel variables $G_1,G_2,G_3$.

Fix a starting vertex $o_N\in\cV_N$. Let $U_N$ be uniform on $\cV_N\setminus\{o_N\}$, independent of the edge weights. We write
\begin{align*}
\pi_N=(x_0=o_N,x_1,\ldots,x_{H_N}=U_N),\quad
\overline T_N:=\fa_NT_N(o_N,U_N)-\log N
\end{align*}
for the geodesic and the centered passage time. The hopcount $H_N$ is defined as the number of edges in the geodesic $\pi_N$. We write $\cL$ for the law of a random element and $d_{\mathrm{TV}}(\mu,\nu):=\sup_A|\mu(A)-\nu(A)|$ for the total variation distance.

\begin{theorem}[Passage time, hopcount and geodesic]\label{thm:main}
Under Assumptions~\ref{ass:weights},~\ref{ass:SN} and~\ref{ass:graphs}, we have the following results.
\begin{enumeratea}
\item\label{main:a} As $N\to\infty$,
\begin{align}\label{eq:graph-joint-limit}
\biggl(\overline T_N,\,\frac{H_N-\theta^{-1}\log N}{\theta^{-1}\sqrt{\log N}}\biggr)
\wto(Y_\theta,Z),
\end{align}
where $Z\sim\normal(0,1)$ is independent of $Y_\theta$.
In particular, the first marginal satisfies
\begin{align}\label{eq:main-tail}
\pr(\overline T_N>x)\to
\E\exp\bigl(-\theta M_\infty^{(1)}M_\infty^{(2)}e^x\bigr),\quad x\in\dR.
\end{align}

\item\label{main:b} Let $(S_j)_{j\ge0}$ be a $P_N$-Markov chain started at $o_N$, independent of $(\overline T_N,H_N)$. Then
\begin{align}\label{eq:graph-path-tv}
\E\biggl[d_{\mathrm{TV}}\Bigl(\cL\bigl((x_j)_{j=0}^{H_N}\mid\overline T_N,H_N\bigr),\,\cL\bigl((S_j)_{j=0}^{H_N}\mid H_N\bigr)\Bigr)\biggr]\to0.
\end{align}
Here the joint law of $(\overline T_N,H_N,\pi_N)$ incorporates the uniform choice of $U_N\in\cV_N\setminus\{o_N\}$, and no terminal vertex is prescribed for $(S_j)_{j\ge0}$.

Moreover, let $(\widetilde S_j)_{j\ge0}$ be a stationary $P_N$-Markov chain, with $\widetilde S_0$ uniform on $\cV_N$. For every fixed integer $k\ge0$ and deterministic functions $\psi_N:\cV_N^{k+1}\to\dR$ with $\sup_N\|\psi_N\|_\infty<\infty$,
\begin{align}\label{eq:window-lln}
\frac1{H_N-k+1}\sum_{i=0}^{H_N-k}\psi_N(x_i,\ldots,x_{i+k})-\E\psi_N(\widetilde S_0,\ldots,\widetilde S_k)\pto0,
\end{align}
where the average is defined to be zero when $H_N<k$.
\end{enumeratea}
The limits hold uniformly in $o_N\in\cV_N$ and over kernels $P_N$ satisfying~\eqref{eq:graph-assumptions} with fixed $r,\eta,C$ (with weak convergence measured by a bounded-Lipschitz metric). For each fixed $k$ and $M<\infty$, the convergence in~\eqref{eq:window-lln} is also uniform over deterministic functions $\psi_N$ with $\|\psi_N\|_\infty\le M$.
\end{theorem}

\begin{remark}\label{rem:universality}
For each $H_N=h$, the conditional laws in~\eqref{eq:graph-path-tv} are probability measures on $\cV_N^{h+1}$. Since $(S_j)_{j\ge0}$ is independent of $(\overline T_N,H_N)$ and the two joint laws share the $(\overline T_N,H_N)$-marginal, the expectation in~\eqref{eq:graph-path-tv} equals their joint total variation distance, so~\eqref{eq:graph-path-tv} is equivalent to
\begin{align*}
d_{\mathrm{TV}}\Bigl(\cL\bigl(\overline T_N,H_N,\pi_N\bigr),\,\cL\bigl(\overline T_N,H_N,(S_j)_{j=0}^{H_N}\bigr)\Bigr)\to0
\end{align*}
on $\bigcup_{h\ge1}(\dR\times\{h\}\times\cV_N^{h+1})$. Equivalently, $(\overline T_N,H_N,\pi_N)$ can be coupled with $(S_j)_{j\ge0}$, with $(S_j)_{j\ge0}$ independent of $(\overline T_N,H_N)$, so that $\pr(\pi_N=(S_0,\ldots,S_{H_N}))\to1$. In this coupling, $(S_j)_{j\ge0}$ is not independent of $U_N$ and $\pr(S_{H_N}=U_N)\to1$; the target $U_N$ is uniform on $\cV_N\setminus\{o_N\}$ prior to conditioning on $(\overline T_N,H_N)$, though its conditional law need not be uniform. Moreover, the probability that $(S_j)_{j=0}^{H_N}$ repeats a vertex tends to zero (Remark~\ref{rem:reference-simple}), consistent with the self-avoidance of $\pi_N$.
Theorem~\ref{thm:main} implies $\fa_NT_N(o_N,U_N)/\log N\pto1$ and $\theta H_N/\log N\pto1$. Under the Weibull assumption, all positive moments also converge. A proof of uniform integrability for $\fa_NT_N(o_N,U_N)/\log N$ and $H_N/\log N$ is given in Appendix~\ref{app:moments}.  Note that any two admissible kernels with the same row-sum scale $S_N$ share the same asymptotic passage-time and hopcount fluctuations, whereas the spatial law of the geodesic, the stationary means in~\eqref{eq:window-lln}, and the CLT for additive functionals (Proposition~\ref{prop:observables}) depend on $P_N$. If $a_NT_N-\log N=O_{\pr}(1)$ and
$a'_N/a_N=1+o(1/\log N)$, then
$a'_NT_N-\log N=a_NT_N-\log N+o_{\pr}(1)$.
\end{remark}

\subsubsection{Extension to general edge-weight distributions}\label{sec:edge-laws}

The Weibull assumption can be replaced by a uniform approximation of the rescaled edge hazards up to the passage-time scale. This allows the edge laws to depend on the endpoints and on $N$, without restricting their upper tails.
For this extension, replace the Weibull weights by independent positive finite weights $(W_e)_{e\in\cE_N}$ with continuous distributions. Fix $c>0$ and define
\begin{align}\label{eq:general-edge-hazard}
\cH_{xy,N}(s):=-\log\pr\bigl(c^{1/\theta}\fa_NW_{x,y}>s\bigr),\quad \{x,y\}\in\cE_N.
\end{align}
For distinct $x,y$ with $P_N(x,y)=0$, set $\cH_{xy,N}\equiv0$, corresponding to the forbidden edge of weight $\infty$. No physical self-edges are added.

\begin{assumption}\label{ass:edge-hazards}
For all large $N$, each $\cH_{xy,N}$ is finite and absolutely continuous on $[0,2\log N]$, with $\cH_{xy,N}(0)=0$. There is a deterministic sequence $\eps_N\ge0$ such that $\eps_N \ll 1/\log N$ and, for every allowed edge

\begin{align}\label{eq:edge-hazard-approximation}
\mathop{\mathrm{ess\,sup}}_{0<s\le2\log N}\biggl|
\frac{\Gamma(\theta)\cH_{xy,N}'(s)}{P_N(x,y)s^{\theta-1}}-1
\biggr|\le\eps_N.
\end{align}
\end{assumption}

\begin{theorem}[Edge-weight universality]\label{thm:edge-universality}
Suppose that Assumptions~\ref{ass:SN},~\ref{ass:graphs} and~\ref{ass:edge-hazards} hold for the independent continuous weights above. Define
\begin{align*}
\overline T_N:=c^{1/\theta}\fa_NT_N(o_N,U_N)-\log N.
\end{align*}
Then the conclusions of Theorem~\ref{thm:main} hold for this $\overline T_N$ with the same limit $Y_\theta$. The convergence holds uniformly in $o_N\in\cV_N$ and, for a fixed sequence $\eps_N$, over edge laws and kernels satisfying the assumptions with fixed $r,\eta,C$. For each fixed $k$ and $M<\infty$, the convergence in~\eqref{eq:window-lln} is also uniform over deterministic functions $\psi_N$ with $\|\psi_N\|_\infty\le M$.
\end{theorem}

The proof is given in Section~\ref{sec:edge-proof}. It embeds the exploration into dominating CMJ trees with reproduction measure $(1+\eps_N)\mu_\theta$, where $\mu_\theta$ is defined in~\eqref{eq:mutheta}. Although the thinned proposals need not vanish in total number, their contribution to the bridge hazard tends to zero because typical ancestral lines have order $\log N$ generations. Passing from the dominating reference pair to the $(\overline T_N,H_N)$-marginal of the model yields~\eqref{eq:graph-path-tv}; together with the mixing estimate for $P_N$, this gives~\eqref{eq:window-lln} and the central limit theorem for additive functionals in Proposition~\ref{prop:observables}.
The following condition on a common noise law suffices for Assumption~\ref{ass:edge-hazards} across the spatial class of Assumption~\ref{ass:graphs}.

\begin{corollary}\label{cor:edge-noise}
Suppose that Assumptions~\ref{ass:SN} and~\ref{ass:graphs} hold. Let $(\go_e)_{e\in\cE_N}$ be independent with a common continuous distribution $F$ on $(0,\infty)$, and define the edge weights by~\eqref{eq:weights}. Suppose that $\cH_0(t):=-\log(1-F(t))$ is absolutely continuous near zero and, for some $c,\kappa>0$,
\begin{align}\label{eq:noise-hazard-expansion}
\cH_0'(t)=c\theta t^{\theta-1}\bigl(1+O(t^\kappa)\bigr)
\text{ as }t\downarrow0.
\end{align}
Then Theorem~\ref{thm:edge-universality} applies.
\end{corollary}

\begin{proof}
For each allowed edge $\{x,y\}\in\cE_N$, the cumulative hazard in~\eqref{eq:general-edge-hazard} is
\begin{align*}
\cH_{xy,N}(s)
=\cH_0\biggl(s\biggl(\frac{P_N(x,y)}{c\Gamma(\theta+1)}\biggr)^{\!1/\theta}\biggr).
\end{align*}
By~\eqref{eq:intro-small-masses}, $\max_{x,y}P_N(x,y)\le(C/N)^{1/(2r)}$. Differentiating and applying~\eqref{eq:noise-hazard-expansion} on $[0,2\log N]$ gives~\eqref{eq:edge-hazard-approximation} with
$
\eps_N=O\bigl((\log N)^\kappa N^{-\kappa/(2r\theta)}\bigr)=o(1/\log N).
$
\end{proof}

For example, a density satisfying $F'(t)=c\theta t^{\theta-1}(1+O(t^\beta))$ near zero gives~\eqref{eq:noise-hazard-expansion} with $\kappa=\min\{\beta,\theta\}$, since $\cH_0'=F'/(1-F)$. This includes $\go=U^{1/\theta}$ for $U$ uniform on $(0,1)$, and gamma noise with shape $\theta$ and rate one. A bound on the distribution function alone does not imply the derivative estimate. In cost notation $W_{x,y}=f_N(x,y)\go_{x,y}$, the intensity is $g_N(x,y)=f_N(x,y)^{-\theta}$.

The following first-order results require only lower-tail bounds on the distribution functions. Atoms are permitted; when geodesics are not unique, $H_N$ may be taken as the hopcount of any geodesic chosen by a measurable rule.

\begin{proposition}[First-order universality]
\label{prop:lower-tail-lln}
Suppose that Assumptions~\ref{ass:SN} and~\ref{ass:graphs} hold. Define the weights by~\eqref{eq:weights} using independent positive finite noises with distribution functions $F_e(t):=\pr(\go_e\le t)$.
\begin{enumeratei}
\item If, for some $c>0$,
\begin{align}\label{eq:cdf-equivalence}
\lim_{t\downarrow0}\sup_N\sup_{e\in\cE_N}
\abs{F_e(t)/t^{\theta}-c}=0,
\end{align}
then
\begin{align}\label{eq:cdf-joint-lln}
\biggl(\frac{c^{1/\theta}\fa_NT_N(o_N,U_N)}{\log N},\,
       \frac{H_N}{\log N}\biggr)
\pto \bigl(1,1/\theta\bigr).
\end{align}
\item If there exist $0<c_1\le c_2<\infty$ and $t_0>0$ such that, for all large $N$,
\begin{align}\label{eq:cdf-two-sided}
c_1t^\theta\le F_e(t)\le c_2t^\theta
\text{ for }0<t\le t_0,\ e\in\cE_N,
\end{align}
then, for every $\eps>0$,
\begin{align}\label{eq:cdf-time-bounds}
\pr\biggl(
c_2^{-1/\theta}-\eps
\le \frac{\fa_NT_N(o_N,U_N)}{\log N}
\le c_1^{-1/\theta}+\eps
\biggr)\to 1,
\end{align}
and there exists $C=C(\theta,c_1,c_2)<\infty$ such that
\begin{align}\label{eq:cdf-hop-bounds}
\pr\bigl(C^{-1}\log N\le H_N\le C\log N\bigr)\to 1.
\end{align}
\end{enumeratei}
Both conclusions hold uniformly in $o_N\in\cV_N$.
\end{proposition}

The proof is given in Section~\ref{sec:first-order-univ}. The fluctuation and path conclusions of Theorem~\ref{thm:edge-universality} rely on the hazard approximation~\eqref{eq:edge-hazard-approximation}.

\subsubsection{Applications}
Section~\ref{sec:applications} verifies Assumption~\ref{ass:graphs} and applies the results above to
\begin{enumeratea}
\item integrable singular kernels on the discrete torus and macroscopic spread-out graphs,
\item Cayley kernels on finite groups, dense regular graphs, and regular graphs of polynomially growing degree,
\item equal-measure discretizations of symmetric Markov kernels and independent point samples of bounded positive symmetric kernels with constant row integrals,
\item geodesic power-law, heat, and Bessel-potential kernels on compact Riemannian manifolds.
\end{enumeratea}

\subsection{Proof mechanism}\label{sec:heuristics}\label{sec:strategy}

We work in rescaled time. By~\eqref{eq:rescaled-law}, an edge from $x$ to $y$ is the first point of a Poisson process with intensity $P_N(x,y)\mu_\theta$. Superposition gives a CMJ process with reproduction intensity $\mu_\theta$, while offspring labels follow $P_N$ independently of reproduction times. Dijkstra's exploration~\cite{Dij59} is obtained by deleting births at occupied labels and their descendants (Lemma~\ref{lem:embedding}). The unpruned population is asymptotic to $M_\infty e^s$.

We grow two clusters simultaneously from $o_N$ and $U_N$. The passage time is twice the first bridge-collision time (Proposition~\ref{prop:half}). Given the internal two-root forest, a bridge whose endpoints were discovered at times $\gs_u,\gs_v$ has already been inspected up to $|\gs_u-\gs_v|$. Its contribution to the cumulative hazard at time $s$ is therefore
\begin{align*}
\frac{P_N(X_u,X_v)}{\Gamma(\theta+1)}
\bigl[(a_u(s)+a_v(s))^\theta-|a_u(s)-a_v(s)|^\theta\bigr].
\end{align*}
The conditional collision survival probability is $e^{-\Lambda_s}$ (Proposition~\ref{prop:cox}). Retaining the lower threshold is essential for the limiting constant.

The bridge mass is of order $e^{2s}/N$, giving the scale $s=\frac12\log N$. Lemma~\ref{lem:transfer} controls the hazard lost through deletions; a population estimate alone would not suffice for this weighted quantity. On the unpruned trees, Proposition~\ref{prop:mixing} uses the spectral gap, the return bound and genealogical correlations to replace $P_N(X_u,X_v)$ by $1/N$ in the cumulative hazard. The uniform second root makes distinct spectral modes orthogonal; conditioning it to avoid $o_N$ costs only $N/(N-1)$ in nonnegative estimates.

The stable ages in the two independent trees have law $\Exp(1)$. For independent $A_1,A_2$ with this law,
\begin{align}\label{eq:heuristic-theta}
\frac{\E((A_1+A_2)^\theta-|A_1-A_2|^\theta)}{\Gamma(\theta+1)}=\theta.
\end{align}
Consequently, at $s=\frac12(\log N+x)$ the hazard converges to $\theta M_\infty^{(1)}M_\infty^{(2)}e^x$, giving~\eqref{eq:main-tail}. Working with cumulative hazards avoids the singular instantaneous rate at age zero when $\theta<1$.

For the path law, the marked first-bridge formula is compared with its unpruned version, and its spatial survival factor is replaced by the label-independent factor $e^{-\overline\Lambda_N(s)}$. Only after this replacement do we sum the ancestral path labels. Symmetry reverses the second ancestral line and yields the exact Markov-path identity~\eqref{eq:graph-path-identity}. The remaining terminal restriction has probability at most $\max_{x,y}P_N(x,y)=o(1)$. Uniformity over bounded path marks gives total variation convergence. The age--generation profile limit (Lemma~\ref{lem:profile-clt}) supplies the independent Gaussian hopcount fluctuation, and Markov-chain mixing gives the empirical path laws. These steps adapt the two-source framework of~\cite{DK26b}; the new spatial estimates are formulated directly for finite kernels.

\subsection{Related work}\label{sec:literature}\label{rem:general-weights}

The complete-graph scalar limits are classical, see~\cite{FG85,Jan99,vdHHvM01} for exponential weights and~\cite[Theorems~1.1--1.2]{BvdH12} for Weibull weights, including the joint independence of the passage-time and hopcount limits. Diameter fluctuations and unusually short or long paths are studied in~\cite{BvdH17,EGvdHN13,EGvdHN20a,EGvdHN20b}. Related flooding and diameter results for the configuration model appear in~\cite{AL15,MS20}. Our contribution is the stability of the typical-pair limits under finite spatial modulation and the joint total variation description of the whole geodesic.

Two-source branching approximations for FPP on sparse random graphs were developed in~\cite{BvdHH10,BvdHH11,BvdHH17}, with dense random graphs treated in~\cite{Bha08}, the epidemic front in~\cite{BvdHK14}, and further weight-law results in~\cite{DSS25}; see also~\cite{AS04,vdH24}. Kolossv\'ary and Komj\'athy~\cite[Theorem~1.2]{KK15} obtained a dense-regime exponential-weight analogue for inhomogeneous random graphs with constant kernel row integrals. Thus neither the constant-row principle nor the mean-field scalar law is new. We instead impose operator conditions directly on deterministic finite matrices, allow every fixed Weibull shape, require no limiting kernel, cover integrable singular profiles, and identify the full spatial path in total variation. Generalized random shortest-path metrics on non-complete deterministic graphs are studied in~\cite{KMV21}. The branching and operator viewpoints in~\cite{BJR07,BBCR10} concern related threshold-percolation questions.

For spatial models, chemical-distance regimes in long-range percolation are studied in~\cite{Sch83,NS86,AN86,BB01,CGS02,Ber02,Bis04,Ber04,Tra10,DS13} with scale-free variants and their weighted distances appear in~\cite{DvdHH13,vdHK17,GGM22,KL20,KLLS26}. Our torus application follows the infinite-volume FPP regimes of~\cite{CD16} and finite-volume first-order results of~\cite{vdHL24}. Random vertex-weight environments with nonconstant normalized row sums require additional analysis; see~\cite{CDK26}. Small-world and gossip models~\cite{WS98,NW99,BR01,BR06,KV16,CD11,BR13} combine local spreading with distant contacts. For spread-out graphs,~\cite{DK26a} treats the microscopic line regime (see also~\cite{Ahl15,CD13}), while~\cite{DK26b} supplies the macroscopic one-dimensional fluctuation result and the two-source framework used here. Corollary~\ref{cor:spreadout} extends that macroscopic example to higher dimensions.

Our branching process inputs come from the CMJ theory of Nerman~\cite{Ner81} and Jagers and Nerman~\cite{JN84}, together with~\cite{Jag75,AN72,LP17,Big95,Shi15}. We use standard finite-group spectral analysis~\cite{Dia88}. The symmetric normalization is related to matrix scaling~\cite{SK67} and bistochastic kernels~\cite{MC19}; Proposition~\ref{prop:zero-diagonal} gives the quantitative perturbation bound needed for the centered FPP limit.

\subsection{Notation}\label{sec:notation}
Asymptotics are as $N=|\cV_N|\to\infty$. On $\Tor=(\dZ/n\dZ)^d$, we have $N=n^d$. For positive sequences $a_N,b_N$, we write $a_N\ll b_N$ if $a_N/b_N\to0$, $a_N\lesssim b_N$ if $a_N\le Cb_N$, and $a_N\asymp b_N$ if $a_N\lesssim b_N\lesssim a_N$. Constants $c,C>0$ may change from line to line and depend on the fixed model parameters $(\theta,r,\eta,d)$ and the constant in~\eqref{eq:graph-assumptions}, but never on $N$ or $o_N$. We write $\wto$ for convergence in distribution, $\pto$ for convergence in probability, $\eqd$ for equality in law, $\Exp(\gl)$ and $\mathrm{Gamma}(\ga,1)$ for exponential and gamma laws, and $\Mit_\theta(z):=\sum_{j\ge0}z^j/\Gamma(j\theta+1)$ for the Mittag-Leffler function.

\subsection{Structure of the paper}
Section~\ref{sec:applications} states verification and normalization criteria for Assumption~\ref{ass:graphs} and applies them to torus kernels, non-transitive graphs, finite groups, and continuum discretizations. Their verification is in Appendix~\ref{app:kernel-verifications}. 
Sections~\ref{sec:prelim}--\ref{sec:mixing} develop the branching, exploration and mixing estimates. Section~\ref{sec:proof-main} proves Theorem~\ref{thm:main}. Section~\ref{sec:edge-proof} treats other edge laws, proves Theorem~\ref{thm:edge-universality} and Proposition~\ref{prop:lower-tail-lln}.  Finally, Section~\ref{sec:observables} derives additive-observable and continuum path limits.  Appendix~\ref{app:moments} proves the Weibull moment statement.

\medskip
\noindent\textbf{Declaration on the use of AI.} AI tools were used as assistive tools for grammar and language editing, consistency checks of notation, and literature searches and the identification of potentially relevant references. All mathematical statements, proofs, citations, references, and final verification were written and reviewed by the authors, who take full responsibility for the content of the manuscript.

\medskip
\noindent\textbf{Acknowledgements.} The authors thank Shirshendu Chatterjee and Neeladri Maitra for helpful discussions during the early stages of this work.

\section{Kernel criteria and applications}\label{sec:applications}

In this section, we establish general criteria for verifying Assumption~\ref{ass:graphs} and apply them to convolution kernels on the torus, regular graphs without vertex transitivity, finite groups, and discretizations of continuum Markov kernels. All parameters in the hypotheses below are fixed independently of $N$. Proofs and supporting discretization estimates are given in Appendix~\ref{app:kernel-verifications}.
The spatial criteria in this section do not depend on the noise distribution. The passage-time, hopcount, and path conclusions also hold for the common noise laws in Corollary~\ref{cor:edge-noise} (using Corollary~\ref{cor:auxiliary-kernel} when a symmetric normalization is involved), with each Weibull normalizer multiplied by $c^{1/\theta}$ and the reference transition kernel unchanged. First-order limits under distribution-function bounds follow from Proposition~\ref{prop:lower-tail-lln}.

\subsection{Verification and normalization}\label{sec:kernel-criteria}

The even-step return bound in Assumption~\ref{ass:graphs} follows from a uniform one-step row $L^p$ estimate for any $p>1$.

\begin{proposition}[Row $L^p$ criterion]\label{prop:lp-smoothing}
Let $P_N$ be a symmetric stochastic matrix on a set of $N$ vertices, and write $K_N(x,y)=NP_N(x,y)$. Suppose that, for some $p>1$ and $A<\infty$,
\begin{align*}
\sup_{N,x}\biggl(\frac1N\sum_y K_N(x,y)^p\biggr)^{\!1/p}\le A.
\end{align*}
Put $\gd=1-1/p$ and $m=\lceil1/\gd\rceil$. Then
\begin{align*}
\sup_{x,y} NP_N^m(x,y)\le A^m,\quad
\max_{x,y}P_N(x,y)\le A N^{-\gd}.
\end{align*}
In particular, $P_N^{2r}(x,x)\le A^m/N$ for every $x$ and every even integer $2r\ge m$.
\end{proposition}

Proposition~\ref{prop:lp-smoothing} provides a sufficient condition for the return bound, not an equivalence: sparser expanding kernels can satisfy a fixed-step return bound even when their one-step row $L^p$ norms diverge for every $p>1$ (see Section~\ref{sec:graph-examples}).

Discretization of a continuum operator often produces an auxiliary symmetric stochastic matrix $Q_N$ with a small nonzero diagonal. Such diagonal entries do not introduce self-edges in the FPP graph. Instead, we construct a zero-diagonal reference kernel $P_N$ by symmetric diagonal scaling and compare its hazards with those of the physical off-diagonal weights.

\begin{proposition}[Removing diagonal terms]\label{prop:zero-diagonal}
Let $Q_N$ be a symmetric stochastic matrix satisfying
\begin{align*}
\|Q_N-\Pi_N\|_{2\to2}\le1-\eta,\quad
\max_xQ_N^{2r}(x,x)\le C/N
\end{align*}
with fixed $r,\eta,C$. Put $\gs_N:=\max_xQ_N(x,x)$. For all large $N$, there is a function $b_N:\cV_N\to(0,\infty)$ with $\max_x|b_N(x)-1|\le C_1\gs_N$ such that
\begin{align}\label{eq:balanced-kernel}
P_N(x,y):=
\begin{cases}
b_N(x)Q_N(x,y)b_N(y),&x\ne y,\\
0,&x=y
\end{cases}
\end{align}
is symmetric stochastic. It satisfies Assumption~\ref{ass:graphs}, with constants depending only on $r,\eta,C$, and
\begin{align*}
\|P_N-Q_N\|_{2\to2}\le C_1\gs_N,\quad
\gs_N\le(C/N)^{1/(2r)}.
\end{align*}
\end{proposition}

\begin{corollary}\label{cor:auxiliary-kernel}
Let $Q_N$ satisfy the hypotheses of Proposition~\ref{prop:zero-diagonal}, and let $P_N$ be given by~\eqref{eq:balanced-kernel}. For $S_N>0$, assign independent weights
\begin{align*}
W_{x,y}=(S_NQ_N(x,y))^{-1/\theta}\go_{x,y},\quad x\ne y,
\end{align*}
with $W_{x,y}:=\infty$ when $Q_N(x,y)=0$. Suppose the common noise law satisfies the hypotheses of Corollary~\ref{cor:edge-noise}. Then the conclusions of Theorem~\ref{thm:edge-universality} hold with reference kernel $P_N$ and normalizer $c^{1/\theta}\fa_N$, where $\fa_N=(\Gamma(\theta+1)S_N)^{1/\theta}$. In particular, Weibull noise is included with $c=1$.
\end{corollary}

This scaling modifies only the reference transition kernel, leaving the physical off-diagonal weights and their support unchanged; naive row-by-row normalization by $1-Q_N(x,x)$ would generally destroy symmetry.

\begin{proposition}[Comparison]\label{prop:comparison-products}
Suppose $Q_N$ is symmetric stochastic and satisfies
\begin{align*}
\|Q_N-\Pi_N\|_{2\to2}\le1-\eta,\quad
\sup_xQ_N^{2r}(x,x)\le \frac{C}{N}.
\end{align*}
Let $P_N$ be a symmetric stochastic matrix. If $P_N(x,y)\ge cQ_N(x,y)$ for all $x,y$ with a fixed $c>0$, then
\begin{align*}
\|P_N-\Pi_N\|_{2\to2}\le1-c\eta.
\end{align*}
If in addition $P_N(x,y)\le C_0Q_N(x,y)$ for all $x,y$ with a fixed $C_0<\infty$, then
\begin{align*}
\sup_xP_N^{2r}(x,x)\le \frac{C_0^{2r}C}{N}.
\end{align*}
Moreover, finite tensor products of admissible kernels are admissible. If $P_N=\bigotimes_{j=1}^{\ell}P_{N_j}^{(j)}$ on $\prod_{j=1}^\ell\cV_{N_j}$ with fixed $\ell$ and $N=\prod_{j=1}^\ell N_j$, and each factor satisfies the two bounds in~\eqref{eq:graph-assumptions} at its size $N_j$, with fixed constants $r_j,\eta_j,C_j$ independent of $N$, then $P_N$ satisfies Assumption~\ref{ass:graphs}.
\end{proposition}

Both matrices in Proposition~\ref{prop:comparison-products} must be stochastic. Small row-sum errors can instead be handled by Proposition~\ref{prop:zero-diagonal} and Corollary~\ref{cor:auxiliary-kernel}, after verifying their hypotheses. The tensor product updates all $\ell$ coordinates at each step, whereas a random-scan product walk updates a single coordinate per step and need not satisfy the fixed-step return bound. On $\{1,2,\ldots,n\}^{\ell}$ with fixed $\ell\ge2$, selecting one prescribed coordinate at all $2r$ updates gives $P^{2r}(x,x)\ge\ell^{-2r}/n$, which is not $O(n^{-\ell})$.

\subsection{Torus kernels and long-range models}\label{sec:torus}\label{sec:powerlaw}\label{sec:spreadout}\label{sec:complete}\label{sec:exponential}

Let $\Tor=(\dZ/n\dZ)^d$, with $N=n^d$, and write $|z|:=\min_{w\in\dZ^d}\|z+nw\|_2$ and $B_\ell:=\{z\in\Tor:|z|\le\ell\}$. On $\dT^d=(\dR/\dZ)^d$, let $\|x\|$ denote the Euclidean distance to $0$. For symmetric nonnegative weights $g_n$ on $\Tor\setminus\{0\}$, set
\begin{align}\label{eq:gSnu}
S_N=\sum_{z\ne0}g_n(z),\quad
\nu_N(z)=\frac{g_n(z)}{S_N}\ (z\ne0),\quad \nu_N(0)=0,
\quad P_N(x,y)=\nu_N(y-x).
\end{align}
The distance-dependent costs $W_{x,y}=f_n(|x-y|)\go_{x,y}$ correspond to $g_n(z)=f_n(|z|)^{-\theta}$, with $f_n=\infty$ permitted on forbidden displacements. Write $m_N:=\max_z\nu_N(z)$ and let $\overline\nu_N$ be the push-forward of $\nu_N$ to $\dT^d$ under $z\mapsto z/n$.
For $k\in\Tor$, define $\widehat\nu_N(k)=\sum_{z\in\Tor}\nu_N(z)e^{2\pi i k\cdot z/n}$. Then
\begin{align}\label{eq:graph-torus-equivalence}
\|P_N-\Pi_N\|_{2\to2}
=\max_{k\in\Tor\setminus\{0\}}|\widehat\nu_N(k)|,\quad
NP_N^{2r}(x,x)=\sum_{k\in\Tor}|\widehat\nu_N(k)|^{2r}.
\end{align}
Consequently, Assumption~\ref{ass:graphs} holds if and only if $\max_{k\ne0}|\widehat\nu_N(k)|\le1-\eta$ and $\sum_{k\in\Tor\setminus\{0\}}|\widehat\nu_N(k)|^q\le C$ for some fixed $q\in(0,\infty)$, taking any integer $r\ge q/2$. The polynomial decay of $m_N$ then follows from~\eqref{eq:intro-small-masses} without a separate hypothesis.

\begin{proposition}[Continuum profiles on the torus]\label{prop:torus-kernel}\label{prop:torus-profile}
Let $k:\dT^d\setminus\{0\}\to[0,\infty)$ be even, Riemann integrable on every closed subset of $\dT^d\setminus\{0\}$, and satisfy
\begin{align*}
0<I:=\int_{\dT^d} k(x)\,\dd x<\infty,
\quad k(x)\le C\|x\|^{-\beta},\quad 0\le\beta<d.
\end{align*}
Define $\nu_n(0)=0$ and
\begin{align*}
\nu_n(z)=\frac{k(z/n)}{\sum_{w\ne0}k(w/n)},\quad z\in\Tor\setminus\{0\}.
\end{align*}
Then the associated convolution matrices $P_N(x,y)=\nu_n(y-x)$ satisfy Assumption~\ref{ass:graphs}.
\end{proposition}

Weak convergence $\overline\nu_N\wto I^{-1}k(x)\,\dd x$ controls only fixed Fourier modes, whereas Riemann integrability away from the origin yields a discrete $L^1$ trigonometric approximation that rules out concentration at high frequencies $|k|_\infty\to\infty$. Proposition~\ref{prop:torus-kernel} covers indicators of symmetric macroscopic annuli and other symmetric Jordan measurable sets of positive volume, anisotropic profiles satisfying the displayed bound, exponential and Gaussian profiles at macroscopic bandwidth, and logarithmically modified profiles $\|x\|^{-\beta}\log^\gamma(2+\|x\|^{-1})$ for $0\le\beta<d$ and $\gamma\in\dR$. Since translation invariance preserves constant row sums, Proposition~\ref{prop:comparison-products} further shows that multiplying a torus profile by any even displacement factor uniformly bounded above and below by positive constants (possibly depending on $n$) preserves admissibility after normalization. For power-law profiles, norm equivalence also permits replacing the Euclidean norm by any fixed norm on $\dR^d$.

For convolution kernels, the Markov path comparison reduces to a comparison with i.i.d.~increments.

\begin{corollary}\label{cor:torus-path}
If the convolution matrices $P_N$ in~\eqref{eq:gSnu} satisfy Assumption~\ref{ass:graphs}, then Theorem~\ref{thm:main} applies. For the geodesic increments $\Delta_j=x_j-x_{j-1}$, let $(\xi_j)_{j\ge 1}$ be i.i.d.~increments with law $\nu_N$, independent of $(\overline T_N,H_N)$. Then
\begin{align}\label{eq:path-tv}
d_{\mathrm{TV}}\Bigl(\cL\bigl(\overline T_N,H_N,(\Delta_j)_{j=1}^{H_N}\bigr),\,
\cL\bigl(\overline T_N,H_N,(\xi_j)_{j=1}^{H_N}\bigr)\Bigr)\to 0.
\end{align}
Consequently, for every continuous $\gf:\dT^d\to\dR$,
\begin{align*}
\frac1{H_N}\sum_{j=1}^{H_N}\gf(\Delta_j/n)
-\int_{\dT^d}\gf\,\dd\overline\nu_N\pto0,
\end{align*}
and if $\overline\nu_N\wto\overline\nu_\infty$, the empirical measure $H_N^{-1}\sum_{j=1}^{H_N}\gd_{\Delta_j/n}$ converges weakly in probability to $\overline\nu_\infty$. Moreover, for every deterministic sequence $i_N\to\infty$, $d_{\mathrm{TV}}(\cL(x_{i_N\wedge H_N}),\mathrm{Unif}_N)\to0$.
\end{corollary}

\subsubsection{Long-range FPP}
For $0\le\beta<d$, define $c_\beta=\int_{[-1/2,1/2]^d}\|x\|^{-\beta}\,\dd x$. For $g_n(z)=|z|^{-\beta}$, the row sum satisfies $S_N=c_\beta n^{d-\beta}(1+o(1/\log N))$. This estimate, proved together with Corollary~\ref{cor:powerlaw} in Appendix~\ref{app:torus-verifications}, permits replacement of $\fa_N$ by its leading continuum asymptotic via Remark~\ref{rem:universality}.

\begin{corollary}\label{cor:powerlaw}
Let $f(r)=r^\ga$ with $\ga\ge0$ and $\ga\theta<d$. Then
\begin{align*}
\biggl(
(\Gamma(\theta+1)c_{\ga\theta})^{1/\theta}
n^{d/\theta-\ga}T_N(0,U_N)-d\log n,\,
\frac{H_N-\theta^{-1}d\log n}{\theta^{-1}\sqrt{d\log n}}
\biggr)\wto(Y_\theta,Z),
\end{align*}
with independent coordinates, and the empirical measure of rescaled geodesic increments converges weakly in probability to $c_{\ga\theta}^{-1}\|x\|^{-\ga\theta}\,\dd x$.
\end{corollary}

For $\theta=1$, Corollary~\ref{cor:powerlaw} establishes the passage-time fluctuations, hopcount central limit theorem, and increment law corresponding to the first-order typical distance in~\cite{vdHL24}. At the critical boundary $\ga\theta=d$, one has $S_N\asymp\log n$ and $m_N\asymp(\log n)^{-1}$, which violates the polynomial bound~\eqref{eq:intro-small-masses} required by Assumption~\ref{ass:graphs}.

\subsubsection{Macroscopic spread-out graphs}
Let $\ell_n$ be an integer satisfying $\ell_n/n\to\gl\in(0,1/2)$, and take $g_n(z)=\ind\{0<|z|\le\ell_n\}$, so that $S_N=|B_{\ell_n}|-1$.

\begin{corollary}\label{cor:spreadout}
Theorem~\ref{thm:main} and Corollary~\ref{cor:torus-path} hold for the $\ell_n$-spread-out torus with $\fa_N=(\Gamma(\theta+1)(|B_{\ell_n}|-1))^{1/\theta}$ and limiting increment law $\mathrm{Unif}(\{x\in\dT^d:\|x\|\le\gl\})$. For $d=1$, $|B_{\ell_n}|-1=2\ell_n$ for all large $n$.
\end{corollary}

In dimension $d=1$, this recovers the macroscopic spread-out cycle of~\cite{DK26b}. For $d\ge1$, we retain the cardinality $|B_{\ell_n}|-1$ in $\fa_N$, since convergence of $\ell_n/n$ without a rate does not ensure a relative error of $o(1/\log N)$ when $|B_{\ell_n}|-1$ is replaced by its continuum volume asymptotic.

\subsubsection{The complete graph}
Taking $g_n\equiv1$ gives $S_N=N-1$, $P_N=(\mathbf J_N-I)/(N-1)$, nontrivial eigenvalues $-1/(N-1)$, and $P_N^2(x,x)=1/(N-1)$. Theorem~\ref{thm:main} recovers the complete-graph passage-time and hopcount limits of~\cite{Jan99,BvdH12}.

\subsection{Non-transitive graphs and finite groups}\label{sec:graph-examples}

Here we present examples where the base graph is not complete.

\subsubsection{Dense regular graphs}
For a $k_N$-regular graph on $N$ vertices with adjacency matrix $A_N$, take $S_N=k_N$ and $P_N=A_N/k_N$. If $k_N\ge cN$ and $P_N$ has a uniform absolute spectral gap, then $P_N^2(x,x)=1/k_N\le1/(cN)$ verifies Assumption~\ref{ass:graphs} with $r=1$, without requiring vertex transitivity. For instance, replacing each vertex of a fixed connected non-bipartite $d_0$-regular graph $H$ on $k$ vertices by a block of $L$ vertices, and each edge of $H$ by $K_{L,L}$, gives
\begin{align*}
N=kL,\quad S_N=d_0L,\quad
P_N=(A_H/d_0)\otimes(\mathbf J_L/L),\quad NP_N^2(x,x)=k/d_0.
\end{align*}
The nonzero eigenvalues of $P_N$ coincide with those of $A_H/d_0$. Choosing a non-transitive base graph $H$ produces a non-transitive family whose geodesic block projection lies within $o(1)$ in total variation of a simple random walk on $H$ of the corresponding random length, whereas its centered passage-time and hopcount limits do not depend on $H$.

For a concrete weighted application, take two communities of $m=N/2$ vertices, with intensity $a>0$ within communities and $b>0$ between them. Its exact row sum is $S_N=a(m-1)+bm$; the nontrivial eigenvalues of $P_N$ are $[a(m-1)-bm]/S_N$ and $-a/S_N$. Theorem~\ref{thm:main} applies. If $C_N$ counts geodesic steps between communities and $q_N=bm/S_N$, then
\begin{align*}
\frac{C_N}{H_N}\pto\frac{b}{a+b},\qquad
\frac{C_N-q_NH_N}{\sqrt{H_N}}
\wto\frac{\sqrt{ab}}{a+b}Z',
\end{align*}
jointly with~\eqref{eq:graph-joint-limit}, where $Z'$ is independent of $(Y_\theta,Z)$. Indeed, reference-walk crossing indicators are i.i.d.~Bernoulli$(q_N)$, so the claim follows from Proposition~\ref{prop:observables}. This separates the universal passage-time fluctuations from the community structure visible along the path.

\subsubsection{Polynomially growing degree}
If $G_N$ is $k_N$-regular with $k_N\ge N^\gd$ for some fixed $\gd>0$ and every nontrivial eigenvalue of $P_N=A_N/k_N$ has modulus at most $C_0k_N^{-1/2}$, then the absolute spectral gap is $1-o(1)$ and, for any integer $r>1/\gd$,
\begin{align*}
P_N^{2r}(x,x)\le N^{-1}+C_0^{2r}k_N^{-r}\le C/N.
\end{align*}
When $k_N=o(N)$, the one-step row $L^p$ norm $(N/k_N)^{1-1/p}$ diverges for every $p>1$, showing that Assumption~\ref{ass:graphs} holds beyond the scope of Proposition~\ref{prop:lp-smoothing}.

\subsubsection{Random regular graphs}
Let $G_N$ be uniform among simple $k_N$-regular graphs on $N$ vertices, where $Nk_N$ is even and $N^\gd\le k_N\le N-1$ for some fixed $\gd\in(0,1)$. Set $S_N=k_N$ and $P_N=A_N/k_N$. By~\cite[Theorem~A]{TY19}, when $k_N\le N/2$,
\begin{align*}
\|P_N-\Pi_N\|_{2\to2}\le Ck_N^{-1/2}
\end{align*}
with probability tending to one. For $k_N>N/2$, let $\overline A_N$ denote the adjacency matrix of the complement $\overline G_N$, which is uniform $\bar k_N$-regular with $\bar k_N=N-1-k_N$. Since $A_N-k_N\Pi_N=-(I-\Pi_N)-(\overline A_N-\bar k_N\Pi_N)$, applying~\cite[Theorem~A]{TY19} to $\overline G_N$ when $\bar k_N\ge\sqrt N$ and using $\|\overline A_N-\bar k_N\Pi_N\|_{2\to2}\le\bar k_N<\sqrt N$ when $\bar k_N<\sqrt N$ extends the displayed bound to the full range $N^\gd\le k_N\le N-1$. The preceding return-probability estimate therefore verifies Assumption~\ref{ass:graphs} with probability tending to one for any fixed integer $r\ge 1/\gd$ (and with $r=1$ if $k_N\ge c_0N$ for a fixed $c_0>0$). Conditionally on $G_N$, the conclusions of Theorem~\ref{thm:main} thus hold in probability with $\fa_N=(\Gamma(\theta+1)k_N)^{1/\theta}$, and~\eqref{eq:graph-joint-limit} holds unconditionally; the common noise laws of Corollary~\ref{cor:edge-noise} are covered upon multiplying $\fa_N$ by $c^{1/\theta}$.

\begin{remark}[Independent-edge random graphs]\label{rem:indep-edge}
Let $G_N$ be a random graph on $\cV_N$ in which edges $\{x,y\}$ ($x\ne y$) are present independently with probabilities $p_N(x,y)=p_N(y,x)\in[0,1]$, and set $p_N(x,x):=0$. Assume that the expected degree $d_N:=\sum_{y\in\cV_N}p_N(x,y)$ is independent of $x\in\cV_N$ and satisfies $d_N\ge N^\gd$ for all large $N$ and some fixed $\gd\in(0,1)$, and that $\overline P_N(x,y):=p_N(x,y)/d_N$ satisfies Assumption~\ref{ass:graphs}. Conditionally on $G_N$, assign independent weights $W_e=\go_e$ to its edges under Assumption~\ref{ass:weights}. Fix $o_N\in\cV_N$ and choose $U_N$ uniformly from $\cV_N\setminus\{o_N\}$, independently of $G_N$ and the edge weights. With
\begin{align*}
\fa_N:=\bigl(\Gamma(\theta+1)d_N\bigr)^{1/\theta},
\qquad
\overline T_N:=\fa_NT_N(o_N,U_N)-\log N,
\end{align*}
the joint limit~\eqref{eq:graph-joint-limit} and the path comparison~\eqref{eq:graph-path-tv} hold conditionally on $G_N$ in probability (and hence~\eqref{eq:graph-joint-limit} also holds unconditionally), uniformly in $o_N\in\cV_N$. In~\eqref{eq:graph-path-tv}, all laws and expectations are taken conditionally on $G_N$, and the reference chain may be chosen as the simple random walk on $G_N$ started at $o_N$ and conditionally independent of $(\overline T_N,H_N)$. On the $o(1)$-probability event that $G_N$ is disconnected, these variables may be defined arbitrarily.

These conclusions follow from degree and spectral concentration, Proposition~\ref{prop:zero-diagonal}, Corollary~\ref{cor:auxiliary-kernel} and the normalization stability in Remark~\ref{rem:universality}.

Examples in this expected-degree range include Erd\H{o}s--R\'enyi graphs $G(N,p_N)$ and stochastic block models with a fixed number of equal-sized blocks and within- and between-block edge probabilities $a\rho_N$ and $b\rho_N$ for fixed $a,b>0$ with $\rho_N\max\{a,b\}\le1$. The hypotheses also hold for long-range percolation on $\Tor$ with edge probabilities $\gl|x-y|^{-\beta}$ for fixed $0<\gl\le1$ and $0\le\beta<d$ (Proposition~\ref{prop:torus-kernel}).
\end{remark}

\subsubsection{Finite groups}
Let $G_N$ be a finite group of order $N$, and let $\nu_N$ be a symmetric probability measure on $G_N$ with $\nu_N(e)=0$. The Cayley kernel $P_N(x,y)=\nu_N(x^{-1}y)$ is symmetric stochastic and has $x$-independent return probability $P_N^{2r}(x,x)=\nu_N^{*2r}(e)$. For each unitary irreducible representation $\pi\in\widehat G_N$ of dimension $d_\pi$, set $\widehat\nu_N(\pi)=\sum_{g\in G_N}\nu_N(g)\pi(g)$. Symmetry makes these matrices Hermitian, and Plancherel's formula~\cite{Dia88} gives
\begin{align*}
\|P_N-\Pi_N\|_{2\to2}
=\max_{\pi\ne\mathbf1}\|\widehat\nu_N(\pi)\|_{\mathrm{op}},\quad
N\nu_N^{*2r}(e)
=\sum_{\pi\in\widehat G_N}d_\pi\operatorname{Tr}
\bigl(\widehat\nu_N(\pi)^{2r}\bigr).
\end{align*}
Uniform upper bounds $1-\eta$ and $C$, respectively, verify Assumption~\ref{ass:graphs} without requiring commutativity. For abelian groups, all $d_\pi=1$, recovering~\eqref{eq:graph-torus-equivalence}. By Proposition~\ref{prop:lp-smoothing}, the $L^p$ bound $N^{p-1}\sum_{g\in G_N}\nu_N(g)^p\le C$ for a fixed $p>1$ is sufficient for the fixed-step return bound.

\subsection{Continuum kernels and manifolds}\label{sec:continuum-applications}\label{sec:cell-averages}\label{sec:manifolds}

Equal-measure partitions transfer a continuum spectral gap and row $L^p$ bound to the finite-graph model.

\begin{proposition}[Equal-measure discretization]\label{prop:cell-averages}
Let $(X,\mu)$ be a nonatomic probability space and let $K\ge0$ be symmetric, with
\begin{align*}
\int_XK(x,y)\,\mu(\dd y)=1
\text{ for $\mu$-almost every }x.
\end{align*}
Let $T$ be its integral operator, and suppose
\begin{align*}
\|T-\Pi\|_{L^2\to L^2}\le1-\eta,\quad
\mathop{\mathrm{ess\,sup}}_x\int_XK(x,y)^p\,\mu(\dd y)\le A^p
\end{align*}
for some $p>1$. For a measurable partition $(A_i)_{i=1}^N$ into cells of equal measure $1/N$, define
\begin{align*}
Q_N(i,j):=
N\int_{A_i}\int_{A_j}K(x,y)\,\mu(\dd y)\mu(\dd x).
\end{align*}
Then $Q_N$ satisfies the hypotheses of Proposition~\ref{prop:zero-diagonal}, with $\max_iQ_N(i,i)\le A N^{-(1-1/p)}$, and its zero-diagonal normalization $P_N$ satisfies Assumption~\ref{ass:graphs}. For any $S_N>0$, the off-diagonal intensities $g_N(i,j)=S_NQ_N(i,j)$, $i\ne j$, satisfy the conclusions of Corollary~\ref{cor:auxiliary-kernel} with reference kernel $P_N$.
\end{proposition}

With conditional expectation $E_N$ onto the partition, $(E_NTE_N)^{2r}$ need not agree with $E_NT^{2r}E_N$; instead, Jensen's inequality preserves the one-step row $L^p$ bound, from which Proposition~\ref{prop:lp-smoothing} supplies the even-step return estimate and the diagonal bound on $Q_N$. Only equal cell measures are required for the finite-graph result, while vanishing cell diameters are needed only to interpret the partition geometrically.

For bounded positive symmetric kernels with constant row integrals, independent random sampling provides a second discretization scheme.

\begin{corollary}\label{cor:iid-bounded}
Let $X_1,\ldots,X_N$ be i.i.d.~with law $\mu$ on a measurable space, and let $k(x,y)=k(y,x)$ satisfy $0<c_0\le k(x,y)\le C_0$ and $\int k(x,y)\,\mu(\dd y)=\gk_0$ for $\mu$-almost every $x$. Conditionally on $(X_1,\ldots,X_N)$, assign distinct pairs independent weights $W_{i,j}=k(X_i,X_j)^{-1/\theta}\go_{i,j}$ under Assumption~\ref{ass:weights}, and choose $U_N$ uniformly from $\{1,\ldots,N\}\setminus\{o_N\}$, independently of the edge noise. Then the passage-time and hopcount limits~\eqref{eq:graph-joint-limit} hold conditionally in probability given $(X_1,\ldots,X_N)$, and hence unconditionally, with normalizer $\fa_N=(\Gamma(\theta+1)N\gk_0)^{1/\theta}$, and the conditional path comparison~\eqref{eq:graph-path-tv} holds with the zero-diagonal reference kernel constructed in Appendix~\ref{app:iid-verification}. The same conclusions hold for the common noise laws in Corollary~\ref{cor:edge-noise} after multiplying $\fa_N$ by $c^{1/\theta}$.
\end{corollary}

For singular kernels $\rho(x,y)^{-\beta}$ on an Ahlfors $D$-regular space, continuum integrability ($\beta<D$) alone does not give uniform control of independent point samples: when $D/2<\beta<D$, the maximum of the normalized row sums $N^{-1}\sum_{j\ne i}\rho(X_i,X_j)^{-\beta}$ diverges in probability, although each fixed normalized row sum has bounded expectation and is tight. Point evaluations of singular kernels can instead be treated using equal-measure partitions of controlled diameter with separated representatives; the corresponding discretization bounds and independent-sample obstruction are proved in Appendix~\ref{sec:point-samples}.

\subsubsection{Manifold models}
Let $M$ be a compact connected Riemannian manifold without boundary, with dimension $d\ge1$, geodesic distance $\rho$, normalized volume $\mu$, and first positive eigenvalue $\gl_1$ of $-\Delta$. By~\cite[Theorem~2]{GL17}, for all large $N$ there are equal-volume cells and representatives satisfying, with $h_N=N^{-1/d}$,
\begin{align}\label{eq:manifold-mesh}
\mu(A_{N,i})=N^{-1},\quad
\operatorname{diam}(A_{N,i})\le Ch_N,\quad
\rho(z_{N,i},z_{N,j})\ge ch_N\quad(i\ne j).
\end{align}

\begin{corollary}[Manifold kernels]\label{cor:manifold-kernels}\label{cor:manifold-powerlaw}\label{cor:manifold-heat}\label{cor:manifold-Bessel}
Choose a symmetric Markov density $k$ and intensity scale $\chi_N$ from the table below, and let $T$ be its integral operator. In the power-law row, $0\le\beta=\ga\theta<d$, $M$ is homogeneous, and $s_\beta:=\int_M\rho(x,y)^{-\beta}\,\mu(\dd y)$ is independent of $x$; put $\kappa_\beta=(\operatorname{diam}M)^{-\beta}/s_\beta>0$. The heat time $t_0>0$ and all kernel parameters are fixed, while $\ell_N>0$ is arbitrary and deterministic.
\begin{center}
\renewcommand{\arraystretch}{1.2}
\begin{tabular}{lccc}
\hline
Model & $k(x,y)$ & $\chi_N$ & Bound on $\|T-\Pi\|_{2\to2}$ \\
\hline
Geodesic power law & $\rho(x,y)^{-\beta}/s_\beta$ & $\ell_N^{-\beta}s_\beta$ & $1-\kappa_\beta$ \\
Heat & $p_{t_0}(x,y)$ & $1$ & $e^{-t_0\gl_1}$ \\
Bessel potential & $K_\beta(x,y)$ & $\ell_N^{-\beta}$ & $(1+\gl_1)^{-(d-\beta)/2}$ \\
\hline
\end{tabular}
\end{center}
Here $p_t$ is the heat density with respect to $\mu$, and, for $0<\beta<d$ and $a=(d-\beta)/2$,
\begin{align}\label{eq:manifold-singular-kernel}
K_\beta(x,y)=\frac1{\Gamma(a)}\int_0^\infty t^{a-1}e^{-t}p_t(x,y)\,\dd t
\quad(x\ne y),\qquad K_\beta(x,x)=0.
\end{align}
Its operator is $(I-\Delta)^{-a}$, and
\begin{align}\label{eq:manifold-singular-bounds}
c_1\rho(x,y)^{-\beta}\le K_\beta(x,y)\le C_1\rho(x,y)^{-\beta},\qquad x\ne y.
\end{align}
For distinct vertices, take point-sampled intensities
$g_N(i,j)=\chi_Nk(z_{N,i},z_{N,j})$ under~\eqref{eq:manifold-mesh}. In the heat and Bessel cases one may instead use any equal-volume partition and
\begin{align*}
g_N(i,j)=\chi_NN^2\int_{A_{N,i}}\int_{A_{N,j}}k(x,y)\,\mu(\dd y)\mu(\dd x).
\end{align*}
Then the weights $W_{i,j}=g_N(i,j)^{-1/\theta}\go_{i,j}$ under Assumption~\ref{ass:weights} satisfy the conclusions of Theorem~\ref{thm:main}, with
\begin{align}\label{eq:manifold-centering}
\overline T_N=(\Gamma(\theta+1)N\chi_N)^{1/\theta}T_N(o_N,U_N)-\log N.
\end{align}
The reference kernel is the symmetric zero-diagonal normalization from Appendix~\ref{app:manifold-verifications}. For common noises satisfying Corollary~\ref{cor:edge-noise}, multiply the normalizer by $c^{1/\theta}$.
\end{corollary}

In the power-law case this is the cost $(\ell_N\rho(z_{N,i},z_{N,j}))^\ga\go_{i,j}$; choosing $\ell_N=N^{1/d}$ measures distance in mesh units. Homogeneity is needed only in that row: heat and Bessel densities have unit row integrals on any $M$. Under~\eqref{eq:manifold-mesh}, Corollary~\ref{cor:manifold-path} gives the continuum path and occupation limits for all three rows.

\section{Preliminaries}\label{sec:prelim}
\subsection{Rescaled hazards and deferred decisions}\label{sec:hazard}

We work with the rescaled edge weights $\tW_e=\fa_NW_e$ and the rescaled passage times $\fa_NT_N(u,v)$. For distinct vertices $x\ne y$,
\begin{align}\label{eq:Hxy}
\cH_{x,y}(s)=P_N(x,y)\frac{s^\theta}{\Gamma(\theta+1)},\quad
h_{x,y}(s)=P_N(x,y)\frac{s^{\theta-1}}{\Gamma(\theta)}.
\end{align}
Forbidden edges have zero hazard and deterministic infinite weight.

\begin{lemma}\label{lem:poisson-rep}
For a fixed label $x$, let $\Pi_x$ be a Poisson point process on $(0,\infty)\times\cV_N$ with intensity $\mu_\theta(\dd s)P_N(x,y)$. Its restrictions $\Pi_x^y=\{s:(s,y)\in\Pi_x\}$ are independent Poisson processes, and
\begin{align*}
\bigl(\min\Pi_x^y\bigr)_{y\ne x}\eqd\Bigl(\tW_{x,y}\Bigr)_{y\ne x},
\end{align*}
with $\min\emptyset=\infty$. Since $P_N(x,x)=0$, $\Pi_x^x$ is empty almost surely.
\end{lemma}

\begin{proof}
The marking theorem~\cite[Chapter~5]{LP17} gives independence and intensity $P_N(x,y)\mu_\theta$. The probability of no point in $\Pi_x^y$ up to age $s$ is $\exp(-P_N(x,y)s^\theta/\Gamma(\theta+1))$, which matches~\eqref{eq:rescaled-law}.
\end{proof}

This lemma represents the uninspected outgoing edges when a vertex is discovered; it does not assign two independent physical clocks to the two orientations of an undirected edge.

\begin{lemma}\label{lem:deferred}
Let $(W_e)_{e\in I}$ be independent random variables with continuous strictly increasing cumulative hazards $\cH_e$, and let $L_e\ge0$ be constants with $\cH_e(L_e)<\infty$. Conditionally on $W_e>L_e$ for every $e\in I$, the variables $\cH_e(W_e)-\cH_e(L_e)$ are independent unit exponentials. Hence, for $A_e\ge L_e$,
\begin{align*}
\pr(W_e>A_e\ \forall e\in I\mid W_e>L_e\ \forall e\in I)
=\exp\biggl(-\sum_{e\in I}\bigl[\cH_e(A_e)-\cH_e(L_e)\bigr]\biggr).
\end{align*}
\end{lemma}

\begin{proof}
Apply the memoryless property to the independent unit exponentials $\cH_e(W_e)$.
\end{proof}

We apply Lemma~\ref{lem:deferred} only to allowed edges; absent edges impose no condition.

\subsection{The CMJ process}\label{sec:cmj-prelim}

Let $\cZ$ be the CMJ process started with a single ancestor at time $0$, in which each individual lives forever and gives birth at ages forming an independent Poisson point process on $(0,\infty)$ with intensity
\begin{align}\label{eq:mutheta}
\mu_\theta(\dd s):=
\frac{s^{\theta-1}}{\Gamma(\theta)}\dd s,\quad s>0.
\end{align}
For $t\ge0$, let $\cZ(t)$ be the set of individuals born by time $t$ (including the ancestor). Define $Z_\theta(t):=|\cZ(t)|$. Since
\begin{align*}
\int_0^\infty e^{-s}\mu_\theta(\dd s)=1,
\quad
\int_0^\infty se^{-s}\mu_\theta(\dd s)=\theta,
\end{align*}
the process has Malthusian parameter $1$ and mean Malthusian reproduction age $\theta$.
For an individual $u$ we write $\gs_u$ for its birth time, $a_u(t):=t-\gs_u$ for its age at time $t\ge\gs_u$, and $\gen(u)$ for its generation (the root has generation $0$). We write $\dtr(u,v)$ for the graph distance between $u$ and $v$ in the family tree. Set
\begin{align*}
Q_j(t):=\mu_\theta^{*j}((0,t])\quad(j\ge1),\quad Q_0(t):=\ind\{t\ge0\},
\end{align*}
which is the expected number of individuals of generation $j$ born by time $t$. Write $\Mit_\theta(z):=\sum_{j\ge0}z^j/\Gamma(j\theta+1)$, and set $\cZ(t)=\emptyset$, $Z_\theta(t)=0$, and $Q_j(t)=0$ for $t<0$.

\begin{lemma}\label{lem:ML}
For $j\ge1$, $\mu_\theta^{*j}(\dd s)=s^{j\theta-1}\dd s/\Gamma(j\theta)$, hence $Q_j(t)=t^{j\theta}/\Gamma(j\theta+1)$. For $\zeta\in[0,1]$,
\begin{align}\label{eq:ML}
\E\biggl(\sum_{u\in\cZ(t)}\zeta^{\gen(u)}\biggr)=\sum_{j\ge0}\zeta^{j}Q_j(t)=\Mit_\theta(\zeta t^{\theta}).
\end{align}
In particular, $\E Z_\theta(t)=\Mit_\theta(t^{\theta})\sim e^{t}/\theta$ as $t\to\infty$. Moreover, for $z\ge0$,
\begin{align}\label{eq:ML-bounds}
1\le\Mit_\theta(z)\le C_\theta e^{z^{1/\theta}}
\text{ and }
\Mit_\theta(z)-1\le\frac{z}{\Gamma(\theta+1)}\Mit_\theta(z).
\end{align}
\end{lemma}

\begin{proof}
The convolution identity follows from the beta integral
\begin{align*}
\int_0^s r^{a-1}(s-r)^{b-1}\dd r
={s^{a+b-1}\Gamma(a)\Gamma(b)}/{\Gamma(a+b)}.
\end{align*}
Applying Campbell's formula generation by generation gives~\eqref{eq:ML}. The asymptotics $\Mit_\theta(z)\sim\theta^{-1}\exp(z^{1/\theta})$ as $z\to\infty$~\cite[Chapter~4]{GKMR14} yield the first bound in~\eqref{eq:ML-bounds}. For the second bound, log-convexity of $\Gamma$ gives
$
\Gamma(\theta+1)\Gamma((j-1)\theta+1)\le\Gamma(j\theta+1),$ for $j\ge1$.
Thus, $\Mit_\theta(z)-1\le z\Mit_\theta(z)/\Gamma(\theta+1)$ by comparison of coefficients.
\end{proof}

\begin{lemma}\label{lem:cmj-facts}
The process $\cZ$ satisfies the following properties.
\begin{enumerate}[label=\textup{(\roman*)},leftmargin=2em]
\item\label{cmj-i} \textup{(Growth)}
There exists a random variable $M_\infty$ such that
\begin{align}\label{eq:Minf}
e^{-t}Z_\theta(t)\to  M_\infty
\text{ a.s.\ and in }L^2,\quad
\E M_\infty=\frac{1}{\theta},\quad M_\infty>0\ \text{a.s.}
\end{align}
Moreover, $\E[Z_\theta(t)^2]\le Ce^{2t}$ for all $t\ge0$.
For $\theta=1$, $\cZ$ is the unit-rate Yule process and
$M_\infty\sim\Exp(1)$.
\item\label{cmj-ii} \textup{(Stable age distribution)} For every bounded continuous $\phi:[0,\infty)\to\dR$,
\begin{align*}
\frac1{Z_\theta(t)}\sum_{u\in\cZ(t)}\phi(a_u(t))\to \E\phi(A)\text{ a.s.},\quad A\sim\Exp(1).
\end{align*}
\item\label{cmj-iii} \textup{(Stable pedigree)} Let $\psi\ge0$ be bounded and continuous with $\E\psi(A)>0$. Let $u_t$ be sampled from $\cZ(t)$ with probability proportional to $\psi(a_u(t))$ (and chosen arbitrarily if all weights vanish). Then $\gen(u_t)/t\pto1/\theta$ and $(\gen(u_t)-t/\theta)/\sqrt{t/\theta^{2}}\wto\normal(0,1)$.
\item\label{cmj-iv} \textup{(Genealogical correlations)} For $\zeta\in[0,1)$ and $t\ge0$,
\begin{align}\label{eq:genealogy}
\begin{aligned}
\E\biggl(\sum_{\substack{u,v\in\cZ(t)\\ u\neq v}}\zeta^{\dtr(u,v)}\biggr)&\le C\zeta(1+t)^{\theta+1}\bigl(e^{t}+e^{2\zeta^{1/\theta}t}\bigr),\\
\E\biggl(\sum_{u,v\in\cZ(t)}\zeta^{\dtr(u,v)}\biggr)&\le C(1+t)^{\theta+1}\bigl(e^{t}+e^{2\zeta^{1/\theta}t}\bigr),
\end{aligned}
\end{align}
where $C$ depends only on $\theta$.
\end{enumerate}
\end{lemma}

\begin{proof}
The almost-sure and $L^1$ population limits and the stable-age law follow from the general CMJ theorems of Nerman~\cite{Ner81} and Jagers and Nerman~\cite{JN84}; see also~\cite[Section~2.1.1]{DK26b}. For the reproduction point process $\xi$ of a single individual, define $\widehat\xi:=\int_0^\infty e^{-s}\xi(\dd s)$ and $c:=2^{-\theta}$. Then $\E\widehat\xi=1$ and $\var(\widehat\xi)=c<1$, so the $x\log x$ condition holds. The exponentially tilted reproduction measure $e^{-s}\mu_\theta(\dd s)$ is the nonlattice $\mathrm{Gamma}(\theta,1)$ probability law, with mean $\theta$; the reproduction process is locally finite. These verify the renewal hypotheses used here. Since each individual has infinitely many offspring over its lifetime, extinction has probability zero; thus $M_\infty>0$ a.s.\ and $\E M_\infty=(\int_0^\infty se^{-s}\mu_\theta(\dd s))^{-1}=1/\theta$. Since individuals live forever, the stable age density is $e^{-a}$.

To prove the $L^2$ convergence in~\ref{cmj-i}, write $m(t):=\E Z_\theta(t)$ and $m_2(t):=\E[Z_\theta(t)^2]$. The Poisson branching property yields the renewal equation $m_2=m^2+\mu_\theta*m_2$. Iterating this relation using $m(t)\le Ce^t$ and $\int_0^\infty e^{-2s}\mu_\theta(\dd s)=c<1$ gives $m_2(t)\le Ce^{2t}$. Let $\cF_k$ be the $\gs$-field generated by the reproduction processes of individuals in generations less than $k$, and define the generation martingale
\begin{align}\label{eq:generation-martingale}
W_k:=\sum_{\gen(u)=k}e^{-\gs_u},\quad W_0=1.
\end{align}
By the branching property, $\E\sum_{\gen(u)=k}e^{-2\gs_u}=c^k$ and $\E[(W_{k+1}-W_k)^2]=c^{k+1}$, so $W_k$ converges in $L^2$ to a limit $W_\infty$ with $\E W_\infty^2=(1-c)^{-1}$. Decomposing $Z_\theta(t)$ at generation $k$, the generations below $k$ contribute $e^{-t}\sum_{j<k}Q_j(t)\to0$ in $L^1$, while the subtrees rooted at generation $k$ have conditional expectation $e^{-t}m(t-\gs_u)\ind\{\gs_u\le t\}\le Ce^{-\gs_u}$. Since $e^{-s}m(s)\to1/\theta$, dominated convergence gives $\E(M_\infty\mid\cF_k)=W_k/\theta$. Since $\gs(\bigcup_k\cF_k)$ contains the entire CMJ tree, taking $k\to\infty$ identifies $W_\infty=\theta M_\infty$ and hence $\E M_\infty^2=[\theta^2(1-c)]^{-1}$. On the other hand, applying dominated convergence to the renewal iteration for $m_2$ gives
\begin{align}\label{eq:cmj-second-moment-limit}
e^{-2t}m_2(t)\to \frac{1}{\theta^2}\sum_{j\ge0}c^j=\frac{1}{\theta^2(1-c)}=\E M_\infty^2.
\end{align}
Combined with $e^{-t}Z_\theta(t)\to M_\infty$ a.s., this implies $L^2$ convergence and proves~\ref{cmj-i}. Part~\ref{cmj-iii} follows from Lemma~\ref{lem:profile-clt}.

For~\ref{cmj-iv}, decompose ordered pairs $(u,v)$ with $u\neq v$ by their most recent common ancestor $w$ and generation offsets $i=\gen(u)-\gen(w)$ and $i'=\gen(v)-\gen(w)$, so that $\dtr(u,v)=i+i'$ and $(i,i')\neq(0,0)$. For $i,i'\ge1$, the two lines of descent leave $w$ through distinct children, so by the multivariate Mecke formula~\cite[Chapter~4]{LP17} the conditional expectation of the number of such pairs given $\gs_w$ is
\begin{align*}
\iint\mu_\theta(\dd r_1)\mu_\theta(\dd r_2)
Q_{i-1}(t-\gs_w-r_1)Q_{i'-1}(t-\gs_w-r_2)
=Q_i(t-\gs_w)Q_{i'}(t-\gs_w),
\end{align*}
and the right-hand side also gives the conditional expectation when $i=0$ or $i'=0$, since $Q_0(t-\gs_w)=1$ for $\gs_w\le t$. Summing against $\zeta^{i+i'}$ and applying~\eqref{eq:ML}--\eqref{eq:ML-bounds} gives
\begin{align*}
\E\biggl(\sum_{u\neq v}\zeta^{\dtr(u,v)}\biggr)
&=\E\biggl(\sum_{w\in\cZ(t)}F(t-\gs_w)\biggr),\\
F(\tau):=\Mit_\theta(\zeta\tau^\theta)^2-1
&\le C\zeta\tau^\theta e^{2\zeta^{1/\theta}\tau}.
\end{align*}
The root $w$ contributes $F(t)$. The non-root birth intensity satisfies~\cite[Chapter~4]{GKMR14}
\begin{align*}
r_\theta(\gs):=\sum_{j\ge1}\frac{\gs^{j\theta-1}}{\Gamma(j\theta)}
\le C(1+\gs^{\theta-1}\ind\{\gs\le1\})e^{\gs}.
\end{align*} Integrating against the mean birth measure yields
\begin{align*}
\E\biggl(\sum_{u\neq v}\zeta^{\dtr(u,v)}\biggr)\le C\zeta\Bigl(t^{\theta}e^{2\zeta^{1/\theta}t}+e^{t}\int_0^t\tau^{\theta}e^{(2\zeta^{1/\theta}-1)\tau}\dd\tau\Bigr)\le C\zeta(1+t)^{\theta+1}\bigl(e^{t}+e^{2\zeta^{1/\theta}t}\bigr).
\end{align*}
Adding the diagonal $\E Z_\theta(t)\le Ce^{t}$ gives the second bound in~\eqref{eq:genealogy}.
\end{proof}

\begin{lemma}\label{lem:profile-clt}
Let $\Phi$ be the standard normal distribution function. If $\phi:[0,\infty)\to\dR$ is continuous with compact support and $F:\dR\to\dR$ is bounded and continuous, then
\begin{align}\label{eq:profile-clt}
e^{-t}\sum_{u\in\cZ(t)}\phi(a_u(t))F\biggl(\frac{\theta\gen(u)-t}{\sqrt t}\biggr)
\to
M_\infty\biggl(\int_0^\infty\phi(a)e^{-a}\dd a\biggr)\biggl(\int_{\dR}F(z)\dd\Phi(z)\biggr)
\text{ in }L^1.
\end{align}
The conclusion also holds for $F=\ind_{(-\infty,z]}$, for such indicators with thresholds converging to $z$, and for bounded continuous $\phi$ without compact support.
\end{lemma}

\begin{proof}
The mean birth-time measure is $\cU_\theta(\dd t)=\gd_0(\dd t)+r_\theta(t)\dd t$, where
$
r_\theta(t)=\sum_{j\ge1}{t^{j\theta-1}}/{\Gamma(j\theta)}\sim e^t/\theta.
$
Suppose first that $\phi$ is supported in $[0,A]$. For $b\in\dR$, define $c_t:=\exp(b/\sqrt t)$. The identity
\begin{align*}
\sum_{j\ge1}\frac{e^{bj\theta/\sqrt t}(t-a)^{j\theta-1}}{\Gamma(j\theta)}
=c_tr_\theta(c_t(t-a))
\end{align*}
and $(c_t-1)t-b\sqrt t\to b^2/2$ give, for $t>A$, via the convolution formula for $\mu_\theta^{*j}$,
\begin{align*}
\E\biggl[e^{-t}\sum_{u\in\cZ(t)}\phi(a_u(t))\exp\biggl(b\frac{\theta\gen(u)-t}{\sqrt t}\biggr)\biggr]
&=e^{-t-b\sqrt t}\int_0^A\phi(a)c_tr_\theta\bigl(c_t(t-a)\bigr)\dd a\\
&\to \frac{e^{b^2/2}}{\theta}\int_0^A\phi(a)e^{-a}\dd a.
\end{align*}
By the continuity theorem for moment generating functions and linearity, for every $\phi$ and $F$ as in the statement,
\begin{align}\label{eq:profile-mean}
\E\biggl[e^{-t}\sum_{u\in\cZ(t)}\phi(a_u(t))F\biggl(\frac{\theta\gen(u)-t}{\sqrt t}\biggr)\biggr]
\to \frac{c_{\phi,F}}{\theta},
\quad
c_{\phi,F}:=\int_0^\infty\phi(a)e^{-a}\dd a\int_{\dR}F\dd\Phi.
\end{align}

To pass to the random profile, fix $k\ge1$ and let $A_t$ be the left-hand side of~\eqref{eq:profile-clt}, with $A_t^{\ge k}$ its restriction to $\gen(u)\ge k$. Conditionally on $\cF_k$, the subtrees rooted at generation $k$ are independent, so Lemma~\ref{lem:cmj-facts}\ref{cmj-i} and $\E\sum_{\gen(u)=k}e^{-2\gs_u}=c^k$ give
\begin{align*}
\E\bigl|A_t^{\ge k}-\E(A_t^{\ge k}\mid\cF_k)\bigr|^2
\le C\|\phi\|_\infty^2\norm{F}_\infty^2 c^k.
\end{align*}
For a generation-$k$ subtree born at $\gs_u$, write $v$ for an individual's local generation. Its global standardization is
\begin{align*}
\frac{\theta(k+v)-t}{\sqrt t}
=\sqrt{\frac{t-\gs_u}{t}}\,
\frac{\theta v-(t-\gs_u)}{\sqrt{t-\gs_u}}
+\frac{\theta k-\gs_u}{\sqrt t}.
\end{align*}
For fixed $k,\gs_u$, the scale factor and additive shift tend to $1$ and $0$. Tightness from~\eqref{eq:profile-mean} and bounded continuity of $F$ therefore give the same subtree limit, and the bound $Ce^{-\gs_u}$ permits dominated convergence over generation $k$. Consequently, $\E(A_t^{\ge k}\mid\cF_k)\to c_{\phi,F}W_k/\theta$ in $L^1$ as $t\to\infty$. Since $\E|A_t-A_t^{\ge k}|\le Ce^{-t}\sum_{j<k}Q_j(t)\to0$, we obtain
\begin{align*}
\limsup_{t\to\infty}\E|A_t-c_{\phi,F}M_\infty|
\le Cc^{k/2}+|c_{\phi,F}|\E|W_k/\theta-M_\infty|.
\end{align*}
Letting $k\to\infty$ proves~\eqref{eq:profile-clt}. Monotone approximation by continuous functions extends the limit to $F=\ind_{(-\infty,z_t]}$ with $z_t\to z$, since $\Phi$ is continuous, and the compact-support condition on $\phi$ is removed via the tail bound $e^{-t}\E\#\{u\in\cZ(t):a_u(t)>A\}\le e^{-t}\E Z_\theta(t-A)\le Ce^{-A}$.
\end{proof}

\begin{remark}\label{rem:exponent}
For $\zeta\in(0,1)$, the mean generation-discounted population $\E\sum_{u\in\cZ(t)}\zeta^{\gen(u)}$ coincides with the mean population of a CMJ process with reproduction intensity $\zeta\mu_\theta$, whose Malthusian parameter solves $\zeta\int_0^\infty e^{-\ga s}\mu_\theta(\dd s)=\zeta\ga^{-\theta}=1$, i.e., $\ga=\zeta^{1/\theta}$. In the common-ancestor decomposition, the contribution to $\E\sum_{u\ne v}\zeta^{d_{\mathrm{tr}}(u,v)}$ from pairs whose common ancestor is born near time $\gs$ grows like $e^{\gs}e^{2\zeta^{1/\theta}(t-\gs)}$ up to polynomial factors. When $2\zeta^{1/\theta}>1$, pairs with a common ancestor near the root dominate; when $2\zeta^{1/\theta}<1$, close relatives dominate and the off-diagonal sum is linear in $\zeta$. For $\theta=1$ (the Yule process) and each fixed $j$, the expected number of pairs at tree distance $j$ is of order $2^je^t$ as $t\to\infty$, corresponding to $\zeta=1/2$. This explains the threshold $2\zeta^{1/\theta}=1$ used in the spectral decomposition below.
\end{remark}

\subsection{Return probabilities and spectral subspaces}

Write $P=P_N$, $\rho=1-\eta$, and $\zeta=2^{-\theta}$. Choose a real orthonormal eigenbasis $(u_i)_{i=0}^{N-1}$ of $\ell^2(\cV_N)$ for $P$, with $u_0\equiv N^{-1/2}$ and eigenvalues $\gl_0=1,\gl_1,\dots,\gl_{N-1}$. Put
\begin{align*}
v_N^*:=\max_xP^2(x,x),\quad
m_N^*:=\max_{x,y}P(x,y),\quad
q_j:=\max_xP^j(x,x).
\end{align*}

\begin{lemma}\label{lem:graph-spectral}
Under Assumption~\ref{ass:graphs},
\begin{align}
v_N^*&\le(C/N)^{1/r},\quad m_N^*\le\sqrt{v_N^*},
\label{eq:graph-small-masses}\\
\sup_{x\in\cV_N}\sum_{i\ge1:|\gl_i|>\zeta}u_i(x)^2&\le\frac{C\zeta^{-2r}}{N},\label{eq:graph-projector}\\
q_j&\le\frac{C_\theta}{N}+C_\theta v_N^*\zeta^j\quad(j\ge1),\label{eq:graph-return}\\
\max_{x,y\in\cV_N}\bigl|P^j(x,y)-N^{-1}\bigr|&\le\frac{C\rho^{j-2r}}{N}\quad(j\ge2r).\label{eq:graph-mixing}
\end{align}
\end{lemma}

\begin{proof}
Since $\sum_{i=0}^{N-1} u_i(x)^2=1$, H\"older's inequality gives
\begin{align*}
P^2(x,x)=\sum_{i=0}^{N-1}\gl_i^2u_i(x)^2\le\biggl(\sum_{i=0}^{N-1}|\gl_i|^{2r}u_i(x)^2\biggr)^{\!1/r}=\bigl(P^{2r}(x,x)\bigr)^{1/r}\le\Bigl(\frac{C}{N}\Bigr)^{\!1/r},
\end{align*}
and $P(x,y)^2\le\sum_zP(x,z)^2=P^2(x,x)$ gives~\eqref{eq:graph-small-masses}. Restricting the spectral sum for $P^{2r}(x,x)$ to $|\gl_i|>\zeta$ yields~\eqref{eq:graph-projector}. Splitting $P^j(x,x)=N^{-1}+\sum_{i\ge1}\gl_i^ju_i(x)^2$ at $|\gl_i|=\zeta$ and applying~\eqref{eq:graph-projector} for $|\gl_i|>\zeta$ and $|\gl_i|^j\le\zeta^{j-2}\gl_i^2$ for $|\gl_i|\le\zeta$ gives~\eqref{eq:graph-return} for $j\ge2$, while $q_1=0$. Finally, for $j\ge2r$, Cauchy--Schwarz gives
\begin{align*}
\biggl|\sum_{i\ge1}\gl_i^ju_i(x)u_i(y)\biggr|
\le\rho^{j-2r}\biggl(\sum_{i\ge1}|\gl_i|^{2r}u_i(x)^2\biggr)^{\!1/2}\biggl(\sum_{i\ge1}|\gl_i|^{2r}u_i(y)^2\biggr)^{\!1/2}
\le\frac{C\rho^{j-2r}}{N},
\end{align*}
proving~\eqref{eq:graph-mixing}.
\end{proof}

\section{Two-source exploration}\label{sec:exploration}

We follow the two-source exploration framework of~\cite[Section~6]{DK26b}, stated here for the individual edge hazards~\eqref{eq:Hxy}; the identities apply directly to the finite-kernel model.
We work in rescaled time. For $s\ge0$ let
\begin{align*}
\cB_1(s):=\{x\in\cV_N:\fa_NT_N(o_N,x)\le s\},\quad
\cB_2(s):=\{x\in\cV_N:\fa_NT_N(x,U_N)\le s\}
\end{align*}
be the balls of radius $s$ around $o_N$ and $U_N$ in the rescaled metric $\fa_NT_N$. Dijkstra's algorithm~\cite{Dij59} reveals $\cB_1(s)$ as $s$ increases. A vertex $u\in\cB_1(s)$ has \emph{age} $a_u(s):=s-\fa_NT_N(o_N,u)$. For $u\in\cB_1(s)$ and $y\notin\cB_1(s)$, the exploration up to time $s$ reveals only that $\tW_{u,y}>a_u(s)$. For two vertices $u,y\in\cB_1(s)$ not joined by a tree edge, it reveals only that $\tW_{u,y}>|\fa_NT_N(o_N,u)-\fa_NT_N(o_N,y)|$. The same holds for $\cB_2(s)$ with ages $a_v(s):=s-\fa_NT_N(v,U_N)$. We write $N_i(s):=|\cB_i(s)|$. Each ball carries its Dijkstra tree, in which the parent of a vertex is the vertex through which it was discovered.

\begin{proposition}[Half-distance identity]\label{prop:half}
Almost surely, $\fa_NT_N(o_N,U_N)=2\tcol$, where
\begin{align*}
\tcol:=\inf\bigl\{s\ge0:\exists u\in\cB_1(s), v\in\cB_2(s), u\ne v,\text{ with }
\fa_NT_N(o_N,u)+\tW_{u,v}+\fa_NT_N(v,U_N)\le2s\bigr\}.
\end{align*}
\end{proposition}

\begin{proof}
This is~\cite[Proposition~6.5]{DK26b}. Any pair $(u,v)$ satisfying the displayed inequality yields a path from $o_N$ to $U_N$ of length at most $2s$, so $\fa_NT_N(o_N,U_N)\le2\tcol$; conversely, the unique edge of the geodesic straddling the midpoint $\frac12\fa_NT_N(o_N,U_N)$ gives $\tcol\le\frac12\fa_NT_N(o_N,U_N)$.
\end{proof}

The identity reduces the passage time to the collision of two clusters grown for the \emph{same} time. At the collision each cluster has size of order $\sqrt N$, where the branching approximation remains accurate.

\subsection{The internal process and the Cox representation}\label{sec:cox}
We now define the \emph{internal process}. It is a two-cluster exploration that coincides with $(\cB_1,\cB_2)$ up to time $\tcol$ and is defined for all times. Both clusters are grown simultaneously by Dijkstra's algorithm, each from its own root, under two rules:
\begin{enumeratea}
\item an edge $\{u,v\}$ with $u$ in one cluster and $v$ in the other (a \emph{bridge}) is never used;
\item all other edges are used as in Dijkstra's algorithm.
\end{enumeratea}
Thus each vertex is assigned permanently to the cluster that first discovers it, and the two clusters remain disjoint for all times. We retain the notation $\cB_i(s)$ and $N_i(s)$ for this internal process from now on. Let $\gs_u$ be the discovery time of $u$, and set $a_u(s):=s-\gs_u$. This is the standard multi-source Dijkstra exploration from $\{o_N,U_N\}$, so $\gs_u=\min\{\fa_NT_N(o_N,u),\fa_NT_N(U_N,u)\}$ and the assigned root realizes this minimum; before $\tcol$ the sets $\cB_1(s)$ and $\cB_2(s)$ coincide with the two metric balls, whereas for $s>\tcol$ the metric balls need not remain disjoint. Let $\cG$ be the $\gs$-field generated by the internal process at all times (the two-root forest, discovery times and spatial positions).

Consider a pair $(u,v)\in\cB_1(s)\times\cB_2(s)$ and suppose $u$ was discovered first. Before $v$ was discovered, $u$ inspected the edge $\{u,v\}$ up to age $\gs_v-\gs_u$ without reaching $v$. After $v$ is discovered, the edge is a bridge and is not inspected further by the internal process. Thus the only information that $\cG$ carries about the bridge weight is
\begin{align}\label{eq:lower-bound}
\tW_{u,v}>L_{u,v}:=|\gs_u-\gs_v|=|a_u(s)-a_v(s)|.
\end{align}
The bridge closes a path of rescaled length at most $2s$ if $\tW_{u,v}\le a_u(s)+a_v(s)$. The threshold $a_u+a_v$ equals $L_{u,v}$ when the pair is formed and grows at rate $2$ thereafter. Define the accumulated bridge hazard
\begin{align}\label{eq:Lambda-def}
\begin{aligned}
\Lambda_s&:=\sum_{u\in\cB_1(s)}\sum_{v\in\cB_2(s)}\bigl[\cH_{u,v}\bigl(a_u(s)+a_v(s)\bigr)-\cH_{u,v}\bigl(|a_u(s)-a_v(s)|\bigr)\bigr]\\
&\phantom{:}=\frac{1}{\Gamma(\theta+1)}\sum_{u\in\cB_1(s)}\sum_{v\in\cB_2(s)}P_N(u,v)\Phi_\theta\bigl(a_u(s),a_v(s)\bigr),
\end{aligned}
\end{align}
where
\begin{align*}
\Phi_\theta(a,b):=(a+b)^{\theta}-|a-b|^{\theta},\quad 0\le\Phi_\theta(a,b)\le(a+b)^{\theta}.
\end{align*}
Equivalently, $\Lambda_s=\int_0^s\sum_{(u,v)}2h_{u,v}(a_u(r)+a_v(r))\dd r$, where the sum runs over the allowed bridge pairs present at time $r$.

\begin{proposition}\label{prop:cox}
Up to time $\tcol$, the explorations $(\cB_1(s),\cB_2(s))$ coincide with the internal process. Moreover, for every $s\ge0$,
\begin{align}\label{eq:cox}
\pr(\tcol>s\mid\cG)=\exp(-\Lambda_s),\quad\text{hence}\quad\pr\bigl(\fa_NT_N(o_N,U_N)>2s\bigr)=\E\bigl(e^{-\Lambda_s}\bigr),
\end{align}
where $\Lambda_s$ is computed from the internal process.
\end{proposition}

\begin{proof}
The censoring argument of~\cite[Proposition~6.8]{DK26b} identifies the metric-ball explorations with the internal process prior to $\tcol$, and the residual calculation in~\cite[Proposition~6.11]{DK26b} extends to the edge-dependent hazards $\cH_e$ in~\eqref{eq:Hxy}. Indeed, for a fixed two-root forest $\cF$ and a given ordering of discovery times, writing $L_e=|\gs_x-\gs_y|$ for $e=\{x,y\}\in\cE_N$ and $p_e$ for the density of $\tW_e$, the joint density in discovery-time and non-tree-edge coordinates is proportional to
\begin{align*}
\prod_{e\in\cF}p_e(L_e)\prod_{e\in\cE_N\setminus\cF}p_e(w_e)\ind\{w_e>L_e\},
\end{align*}
since the parent--child change of variables has unit Jacobian and the non-tree inequalities $w_e>L_e$ characterize the Dijkstra forest $\cF$. Conditionally on $\cG$, the allowed non-tree weights $\tW_e$ are therefore independent with laws truncated below at $L_e$.

By~\cite[Proposition~6.11(ii)]{DK26b}, combining a bridge with the two forest paths gives
\begin{align*}
\fa_NT_N(o_N,U_N)=\min_{\{u,v\}\text{ a bridge}}(\gs_u+\tW_{u,v}+\gs_v).
\end{align*}
Since $\gs_u+\tW_{u,v}+\gs_v>2(\gs_u\vee\gs_v)$ by~\eqref{eq:lower-bound}, the event $\{\tcol>s\}$ coincides with $\{\tW_{u,v}>a_u(s)+a_v(s)\}$ across all bridges present at time $s$, and Lemma~\ref{lem:deferred} yields~\eqref{eq:cox}.
\end{proof}

\begin{lemma}\label{lem:first-bridge}
For an allowed bridge $e=\{u,v\}\in\cE_N$, set $K_e(s):=0$ for $s<\gs_u\vee\gs_v$ and
\begin{align*}
K_e(s):=\cH_e(2s-\gs_u-\gs_v)-\cH_e(|\gs_u-\gs_v|)\quad(s\ge\gs_u\vee\gs_v).
\end{align*}
Let $e^*$ be the bridge realizing $\tcol$. If $B_e(s)$ is bounded and measurable with respect to $(s,\cG)$, then for every Borel set $I\subseteq[0,\infty)$,
\begin{align}\label{eq:marked-first-bridge}
\E\bigl(B_{e^*}(\tcol)\ind\{\tcol\in I\}\bigm|\cG\bigr)
=\int_I e^{-\Lambda_s}\sum_e B_e(s)\dd K_e(s),
\end{align}
where, for $s>\gs_u\vee\gs_v$,
\begin{align}\label{eq:bridge-rate}
\dd K_{\{u,v\}}(s)
=\frac2{\Gamma(\theta)}P_N(u,v)\bigl(a_u(s)+a_v(s)\bigr)^{\theta-1}\dd s.
\end{align}
\end{lemma}

\begin{proof}
This is the marked form of~\cite[Proposition~6.11(iv)]{DK26b}. Conditionally on $\cG$, the candidate collision time $\frac12(\gs_u+\tW_e+\gs_v)$ of a bridge $e$ has distribution $e^{-K_e(s)}\dd K_e(s)$ on $(\gs_u\vee\gs_v,\infty)$, independently across bridges. Multiplying by the conditional survival probabilities $\prod_{e'\neq e}e^{-K_{e'}(s)}$ of the remaining bridges and summing over $e$ gives~\eqref{eq:marked-first-bridge}, and differentiation gives~\eqref{eq:bridge-rate}.
\end{proof}

\begin{remark}\label{rem:lower-bounds}
The lower bound $L_{u,v}=|\gs_u-\gs_v|=|a_u-a_v|$ in~\eqref{eq:Lambda-def} reflects the fact that the internal exploration ceases to inspect $\{u,v\}$ once its second endpoint is discovered. With independent ages $A_1,A_2\sim\Exp(1)$, the subtracted hazard at $|A_1-A_2|$ gives
\begin{align*}
\E\bigl((A_1+A_2)^\theta-|A_1-A_2|^\theta\bigr)=\Gamma(\theta+2)-\Gamma(\theta+1)=\theta\Gamma(\theta+1),
\end{align*}
which produces the constant $\theta$ in~\eqref{eq:main-tail} and Corollary~\ref{cor:hazard}. By contrast, omitting the lower bound ($L_{u,v}=0$) would reintroduce bridge weights already excluded by the discovery history and give $\Gamma(\theta+2)$, while using $a_u\vee a_v$ would treat the subsequent aging of the older endpoint as information about the bridge and give $\Gamma(\theta+1)(\theta-1+2^{-\theta})$.
\end{remark}

\section{Embedding and deleted descendants}\label{sec:cmj}

We couple the two-root internal process to two independent CMJ trees whose spatial labels evolve according to the symmetric stochastic matrix $P=P_N$. Throughout this section, write
\begin{align*}
\zeta:=2^{-\theta},\quad
q_j:=\max_xP^j(x,x),\quad
v_N^*:=\max_xP^2(x,x).
\end{align*}
By~\eqref{eq:graph-return},
\begin{align*}
q_j\le C/N+C v_N^*\zeta^j\quad(j\ge1).
\end{align*}
All constants below are uniform in $o_N\in\cV_N$.

\subsection{The branching representation}

Let $\cZ_1,\cZ_2$ be independent copies of the CMJ process with reproduction measure $\mu_\theta$, and let $\cT$ be the $\gs$-field generated by their unmarked family trees and birth times. Their roots have labels $o_N$ and $U_N$, respectively, where $U_N$ is uniform on $\cV_N\setminus\{o_N\}$ and independent of $\cT$. Conditionally on the label $x$ of a parent, assign its children independent labels with distribution $P(x,\cdot)$, independently of all reproduction times. Thus an individual at $x$ has marked offspring intensity $\mu_\theta(\dd s)P(x,y)$, with total intensity $\mu_\theta(\dd s)$. Write $X_u$ and $\gs_u$ for an individual's label and birth time, and let $\cZ_i(s)$ be the set of individuals of tree $i$ born by time $s$, with $Z_i(s):=|\cZ_i(s)|$.

We use the same two unmarked CMJ trees for every $N$ and identify $M_\infty^{(i)}=\lim_{t\to\infty}e^{-t}Z_i(t)$. Only $U_N$ and the spatial labels depend on $N$.

Retain both roots and process all subsequent births in increasing order of birth time. A birth is deleted if its parent is deleted or if its label is already occupied by a retained individual of either tree. A deleted birth with a retained parent is called a \emph{ghost}. Let $D_i(s)\subseteq\cZ_i(s)$ be the set of deleted individuals born by time $s$.

\begin{lemma}\label{lem:embedding}
The retained individuals $\cZ_i(s)\setminus D_i(s)$, together with their birth times, labels and parent relations, have the same joint law as the two-root internal process $(\cB_1(s),\cB_2(s))$ and its discovery forest.
\end{lemma}

\begin{proof}
When a vertex $x$ is discovered by the internal process, its allowed edges to vertices $y$ outside both clusters have not yet been inspected, so their weights remain independent with their original laws. The first point of a Poisson process with intensity $P(x,y)\mu_\theta$ has survival function
\begin{align*}
\exp\bigl(-P(x,y)s^\theta/\Gamma(\theta+1)\bigr),
\end{align*}
which is the law of the rescaled physical edge weight for $x\ne y$. Superposing these processes and adjoining independent processes for all currently occupied labels yields the marked offspring intensity $\mu_\theta\otimes P(x,\cdot)$.

An active edge to an unoccupied vertex retains its original clock: its first point discovers the target if that target is still unoccupied, and otherwise the proposal is suppressed. Every subsequent point for the same parent--target pair also finds the target occupied and is suppressed. Proposals to vertices already occupied when $x$ was discovered, as well as the offspring processes of deleted individuals, are not inspected by the internal exploration and may be sampled independently. Since $\E Z_\theta(t)<\infty$, only finitely many proposals occur on each bounded time interval, and induction over their birth times gives the equality in law.
\end{proof}

In the estimates below, two consequences of symmetry replace translation invariance. First, any nonnegative function of the unpruned labelled trees has expectation under $U_N\ne o_N$ at most $N/(N-1)$ times its expectation when the second root is uniform on all of $\cV_N$. Second, under this unrestricted root law, every individual of $\cZ_2$ has a uniform label, and the labelling law is invariant under rerooting: on any finite connected subtree, the joint label probability is
\begin{align*}
\frac1N\prod_{\{u,v\}\text{ a tree edge}}P(X_u,X_v),
\end{align*}
which is independent of the orientation because $P$ is symmetric. We use this property only for unpruned-tree estimates; the two physical roots remain distinct.

\subsection{Deleted populations}

\begin{proposition}\label{prop:coupling}
For every $s\ge0$ and $i=1,2$,
\begin{align}\label{eq:lost}
\E|D_i(s)|\le C\biggl(\frac{e^{2s}}N+v_N^*(1+s)^{\theta+2}e^s\biggr).
\end{align}
\end{proposition}

\begin{proof}
Every deleted individual descends from a ghost. For a non-root individual $w$, let $\cZ^{(w)}(s-\gs_w)$ denote all its descendants born by time $s$, including $w$. A ghost $w$ shares its label with an earlier retained individual, whence
\begin{align*}
|D_i(s)|\le
\sum_{\substack{w\in\cZ_i(s)\\w\text{ not a root}}}
\sum_{\substack{w'\in\cZ_1\cup\cZ_2\\\gs_{w'}<\gs_w}}
\ind\{X_w=X_{w'}\}\,|\cZ^{(w)}(s-\gs_w)|.
\end{align*}
If $w,w'$ belong to the same tree, let their distances from their most recent common ancestor be $p,q$, so $p+q=\dtr(w,w')\ge1$. Conditional on that ancestor's label $a$, their coincidence probability is
\begin{align*}
\sum_yP^p(a,y)P^q(a,y)=P^{p+q}(a,a)\le q_{p+q}.
\end{align*}
Parent--child coincidences have probability zero since $P(a,a)=0$. For individuals in different trees, stationarity of the second tree and conditioning its root to differ from $o_N$ bound the coincidence probability by $1/(N-1)$.

Since $w'$ is born before $w$, it lies outside the descendant subtree of $w$. By the branching property, the conditional expectation of $|\cZ^{(w)}(s-\gs_w)|$ is $m(s-\gs_w)\le Ce^{s-\gs_w}$, where $m(t):=\E Z_\theta(t)$. Applying the coincidence bounds above gives
\begin{align*}
\E|D_i(s)|\le C\E\sum_{w\in\cZ_i(s)}e^{s-\gs_w}
\biggl(\frac{Z_1(\gs_w)+Z_2(\gs_w)}N
+v_N^*\sum_{\substack{w'\in\cZ_i,\ \gs_{w'}<\gs_w}}
\zeta^{\dtr(w,w')}\biggr).
\end{align*}
For the population term, integration by parts against $Z_i$ and the bound $\E(Z_i(t)^2)\le Ce^{2t}$ give
\begin{align*}
\E\sum_{w\in\cZ_i(s)}e^{-\gs_w}Z_i(\gs_w)\le Ce^s.
\end{align*}
The term with $Z_{3-i}$ satisfies the same bound by independence of the unmarked trees and $m(t)\le Ce^t$, yielding a contribution of at most $Ce^{2s}/N$.

For the remaining term, define
\begin{align*}
A_i(t):=\sum_{w\in\cZ_i(t)}
\sum_{\substack{w'\in\cZ_i\\\gs_{w'}<\gs_w}}
\zeta^{\dtr(w,w')}.
\end{align*}
Since $2\zeta^{1/\theta}=1$, Lemma~\ref{lem:cmj-facts}\ref{cmj-iv} gives $\E A_i(t)\le C(1+t)^{\theta+1}e^t$. Because $A_i$ is nondecreasing with $A_i(0)=0$, integration by parts yields
\begin{align*}
\sum_{w\in\cZ_i(s)}e^{s-\gs_w}
\sum_{\substack{w'\in\cZ_i\\\gs_{w'}<\gs_w}}
\zeta^{\dtr(w,w')}
=\int_{[0,s]}e^{s-t}\dd A_i(t)
=A_i(s)+\int_0^se^{s-t}A_i(t)\dd t,
\end{align*}
whose expectation is at most $C(1+s)^{\theta+2}e^s$. This proves~\eqref{eq:lost}.
\end{proof}

\subsection{Deleted bridge mass and cumulative hazards}

Under the preceding coupling, write $a_u(s):=s-\gs_u$ for an individual present at time $s$, and define the unpruned cumulative hazard by
\begin{align*}
\Lambda_s^{\cZ}:=
\frac1{\Gamma(\theta+1)}
\sum_{u\in\cZ_1(s)}\sum_{v\in\cZ_2(s)}
P(X_u,X_v)\Phi_\theta(a_u(s),a_v(s)),
\quad
\Phi_\theta(a,b):=(a+b)^\theta-|a-b|^\theta.
\end{align*}
The corresponding sum over retained pairs equals the internal-process hazard $\Lambda_s$. Its endpoints always have distinct labels. The unpruned sum may include coincident-label pairs, which carry zero weight since $P_N(x,x)=0$. Define
\begin{align*}
\cR_N(s):=\sum_{u\in\cZ_1(s)}\sum_{v\in\cZ_2(s)}
P(X_u,X_v)\ind\{u\in D_1(s)\text{ or }v\in D_2(s)\}.
\end{align*}

\begin{lemma}\label{lem:transfer}
For all $s\ge0$,
\begin{align}
\E\cR_N(s)&\le C(1+s)^{\theta+2}
\biggl(\frac{e^{3s}}{N^2}+\frac{v_N^*e^{2s}}N\biggr),
\label{eq:lost-bridge-mass}\\
\E(\Lambda_s^{\cZ}-\Lambda_s)&\le C(1+s)^{2\theta+2}
\biggl(\frac{e^{3s}}{N^2}+\frac{v_N^*e^{2s}}N\biggr).
\label{eq:transfer-quantitative}
\end{align}
In particular, the second expectation tends to zero at $s=\frac12(\log N+x)$ for every fixed $x$. Moreover, $\dd\Lambda^{\cZ}-\dd\Lambda$ is a nonnegative measure. For a bounded pair mark $B$ and $0\le a\le b$, replacing retained pairs by all CMJ pairs in a first-collision integral over $(a,b]$, while retaining the survival factor $e^{-\Lambda_s}$, changes the integral by at most $\|B\|_\infty\E(\Lambda_b^{\cZ}-\Lambda_b)$.
\end{lemma}

\begin{proof}
Bound the indicator of a deleted endpoint by the sum of indicators that some ancestor $w$ shares its label with an earlier individual $w'$. We treat coincidences within a single tree and between the two trees separately.

For a within-tree coincidence, let $w'$ precede $w$ in the same tree, let $u\succeq w$ be a descendant present at time $s$, and let $v$ be any individual of the other tree present at time $s$. Put $j=\dtr(w,w')\ge1$. If $w\in\cZ_1$, averaging $X_v$ over the stationary second tree gives
\begin{align*}
\E\bigl(\ind\{X_w=X_{w'}\}P(X_u,X_v)\bigm|\cT\bigr)
\le q_j/(N-1).
\end{align*}
If $w\in\cZ_2$, reroot the stationary second tree at $w$. Since $w'$ lies outside the descendant subtree of $w$, the label paths from $w$ to $u$ and to $w'$ are conditionally independent given $X_w$. With $p=\dtr(w,u)$, under the unrestricted uniform-root law,
\begin{align*}
\pr(X_u=z,\ X_w=X_{w'}\mid\cT)
=\frac1N\sum_aP^p(a,z)P^j(a,a)\le q_j/N.
\end{align*}
Summing against $P(X_v,z)$ and conditioning the second root to differ from $o_N$ yields the same bound $q_j/(N-1)$. 

To sum this bound, let
\begin{align*}
L_w(s)&:=\#\{u\in\cZ_i(s):u\succeq w\},\\
B_{\alpha,i}(t)&:=\sum_{\substack{w\in\cZ_i(t)\\w\text{ not a root}}}
\sum_{\substack{w'\in\cZ_i\\\gs_{w'}<\gs_w}}
\alpha^{\dtr(w,w')},\qquad \alpha\in\{1,\zeta\}.
\end{align*}
At the birth of $w$, its reproduction process is fresh, while all the
coefficients in the corresponding increment of $B_{\alpha,i}$ are determined
by the birth history through $\gs_w$. Hence the branching property and Tonelli's
theorem give
\begin{align*}
\E\sum_{\substack{w\in\cZ_i(s)\\w\text{ not a root}}}
L_w(s)\sum_{\substack{w'\in\cZ_i\\\gs_{w'}<\gs_w}}\alpha^{\dtr(w,w')}
=\int_{(0,s]}m(s-t)\,\dd\E B_{\alpha,i}(t).
\end{align*}
Birth times are almost surely distinct, so
$B_{1,i}(t)=Z_i(t)(Z_i(t)-1)/2$ and
$\E B_{1,i}(t)\le Ce^{2t}$. The genealogical bound gives
$\E B_{\zeta,i}(t)\le C(1+t)^{\theta+1}e^t$.
Using $m(s-t)\le Ce^{s-t}$ and integration by parts therefore bounds the last
display by $Ce^{2s}$ for $\alpha=1$ and by
$C(1+s)^{\theta+2}e^s$ for $\alpha=\zeta$.
Since $q_j\le C/N+Cv_N^*\zeta^j$, independence of the other unmarked tree now
bounds the total within-tree contribution by
\begin{align*}
\frac{Cm(s)}{N-1}
\biggl(\frac{e^{2s}}N+v_N^*(1+s)^{\theta+2}e^s\biggr)
\le C\biggl(\frac{e^{3s}}{N^2}
+\frac{v_N^*(1+s)^{\theta+2}e^{2s}}N\biggr).
\end{align*}

For a between-tree coincidence, let $w$ be a birth in one tree, $v_0$ an earlier individual of the other tree, $u\succeq w$ a descendant present at time $s$, and $v$ any individual in the tree of $v_0$ present at time $s$. Put $p=\dtr(w,u)$ and $q=\dtr(v_0,v)$. We claim that
\begin{align}\label{eq:graph-cross-identity}
\E\bigl(\ind\{X_w=X_{v_0}\}P(X_u,X_v)\bigm|\cT\bigr)
\le\frac{q_{p+q+1}}{N-1}
\le\frac{C}{N-1}\bigl(N^{-1}+v_N^*\zeta^{p+q}\bigr).
\end{align}
When $w\in\cZ_1$, rerooting the stationary second tree at $v_0$ and conditioning on the common label $X_w=X_{v_0}=a$ gives
\begin{align*}
\sum_{y,z}P^p(a,y)P(y,z)P^q(a,z)=P^{p+q+1}(a,a)
\end{align*}
with stationary weight $1/N$. When $w\in\cZ_2$, rerooting $\cZ_2$ at $w$ and conditioning on the label $a$ of the most recent common ancestor of $v_0,v$ in $\cZ_1$ (with branch lengths $q_1+q_2=q$) gives
\begin{align*}
\frac1N\sum_{y,z}P^{q_1}(a,y)P^{p+1}(y,z)P^{q_2}(a,z)
=\frac1N P^{p+q+1}(a,a).
\end{align*}
Multiplying by $N/(N-1)$ for the conditioned root and applying~\eqref{eq:graph-return} proves~\eqref{eq:graph-cross-identity}; the case $p=q=0$ vanishes since $P_N(x,x)=0$.

We make the four-individual sum in~\eqref{eq:graph-cross-identity} explicit.
For $w\in\cZ_i$, write $j=3-i$, let $\cU_\theta^\circ:=\cU_\theta-\gd_0$ be the mean birth-time
measure excluding the root, and set
\begin{align*}
L_{\zeta,w}(s)&:=\sum_{\substack{u\in\cZ_i(s)\\u\succeq w}}
\zeta^{\dtr(w,u)},\\
J_j(t,s)&:=\sum_{\substack{v_0\in\cZ_j\\\gs_{v_0}<t}}
\sum_{v\in\cZ_j(s)}\zeta^{\dtr(v_0,v)}.
\end{align*}
For a birth $w$ at time $t$, the expectations of $L_w(s)$ and
$L_{\zeta,w}(s)$, conditional on the birth history of $\cZ_i$ through $t$, are
$m(s-t)$ and
$m_\zeta(s-t):=\E\sum_{u\in\cZ(s-t)}\zeta^{\gen(u)}$, respectively.
The birth of $w$ is determined by its ancestors and does not condition its own
reproduction. The other unmarked tree is independent of this history and of
the subtree rooted at $w$. Summing over the birth times of $w$, or equivalently
summing generation by generation and using
$\cU_\theta^\circ=\sum_{k\ge1}\mu_\theta^{*k}$, gives
\begin{align*}
\E\sum_{\substack{w\in\cZ_i(s)\\w\text{ not a root}}}
L_w(s)Z_j(\gs_w-)Z_j(s)
&=\int_{(0,s]}m(s-t)
\E\bigl[Z_\theta(t-)Z_\theta(s)\bigr]\cU_\theta^\circ(\dd t),\\
\E\sum_{\substack{w\in\cZ_i(s)\\w\text{ not a root}}}
L_{\zeta,w}(s)J_j(\gs_w,s)
&=\int_{(0,s]}m_\zeta(s-t)\E J_j(t,s)\cU_\theta^\circ(\dd t).
\end{align*}
These identities use only the birth history of $w$ and the independent other
tree, rather than conditioning on the entire future family tree.
For the first identity, Cauchy--Schwarz and
Lemma~\ref{lem:cmj-facts}\ref{cmj-i} give
$\E(Z_\theta(t-)Z_\theta(s))\le Ce^{t+s}$.
Together with $m(s-t)\le Ce^{s-t}$ and
$\cU_\theta^\circ((0,s])=m(s)-1\le Ce^s$, this bounds the $N^{-2}$ term of
\eqref{eq:graph-cross-identity} by $Ce^{3s}/N^2$.

For the second identity, dropping the restriction $\gs_{v_0}<t$ gives
\begin{align*}
\E J_j(t,s)\le G_\zeta(s)
:=\E\sum_{v_0,v\in\cZ(s)}\zeta^{\dtr(v_0,v)}
\le C(1+s)^{\theta+1}e^s.
\end{align*}
Moreover, the convolution formula from Lemma~\ref{lem:ML} yields
\begin{align*}
\int_{[0,s]}m_\zeta(s-t)\cU_\theta(\dd t)
&=\sum_{n\ge0}Q_n(s)\sum_{k=0}^n\zeta^k
\le\frac{m(s)}{1-\zeta}\le Ce^s.
\end{align*}
Consequently the second term of~\eqref{eq:graph-cross-identity}, summed over
both orientations of the coincidence, is at most
\begin{align*}
\frac{Cv_N^*}{N}G_\zeta(s)
\int_{[0,s]}m_\zeta(s-t)\cU_\theta(\dd t)
\le\frac{Cv_N^*}{N}(1+s)^{\theta+1}e^{2s},
\end{align*}
which proves~\eqref{eq:lost-bridge-mass}.

Multiplying~\eqref{eq:lost-bridge-mass} by $\Phi_\theta(a_u(s),a_v(s))\le(2s)^\theta$ gives~\eqref{eq:transfer-quantitative}, which vanishes at $s=\frac12(\log N+x)$ since $v_N^*=O(N^{-1/r})$. Finally, since retention is permanent at birth, $\Lambda_s^{\cZ}-\Lambda_s$ is a sum of nondecreasing pair cumulative hazards, so $\dd\Lambda^{\cZ}-\dd\Lambda\ge0$ as a Stieltjes measure, and the marked integral bound follows from $e^{-\Lambda_s}\le1$.
\end{proof}

\section{Spatial mixing and cumulative hazards}\label{sec:mixing}

We first compare a bilinear form on the unpruned labelled trees with its uniform average. The second tree is stationary under an unrestricted uniform root; conditioning its root to differ from $o_N$ increases expectations of non-negative quantities by at most a factor $N/(N-1)$.

\begin{proposition}[Bilinear spectral mixing]\label{prop:mixing}
Suppose Assumptions~\ref{ass:weights} and~\ref{ass:graphs} hold. Let $w_u^i\in[-1,1]$ be $\cT$-measurable weights. Define
\begin{align*}
R_N^P(s):=\sum_{u\in\cZ_1(s)}\sum_{v\in\cZ_2(s)}P(X_u,X_v)w_u^1w_v^2-\frac1N\biggl(\sum_{u\in\cZ_1(s)}w_u^1\biggr)\biggl(\sum_{v\in\cZ_2(s)}w_v^2\biggr).
\end{align*}
With $\eta'=\min\{1/2,1-(1-\eta)^{1/\theta}\}$,
\begin{align}\label{eq:graph-bilinear}
\E\bigl(|R_N^P(s)|^2\bigr)\le C(1+s)^{2\theta+2}\biggl(\frac{v_N^*e^{2s}}N+\frac{e^{4s}}{N^2}e^{-4\eta's}\biggr).
\end{align}
All constants are independent of $o_N\in\cV_N$.
\end{proposition}

\begin{proof}
Set
\begin{align*}
A_{i,k}&:=\sum_{u\in\cZ_i(s)}w_u^iu_k(X_u),\qquad
G_i(a):=\sum_{u,u'\in\cZ_i(s)}a^{\dtr(u,u')},\\
R_N^P(s)&=\sum_{k=1}^{N-1}\gl_k A_{1,k}A_{2,k}.
\end{align*} Since $|R_N^P(s)|^2\ge0$, its conditional expectation given $\cT$ under $U_N\sim\mathrm{Unif}(\cV_N\setminus\{o_N\})$ is at most $N/(N-1)$ times that under an unrestricted uniform root on $\cV_N$ for $\cZ_2$. Under the unrestricted root law, conditioning on the label of the most recent common ancestor $v\wedge v'$ and using $P^du_k=\gl_k^du_k$ and the orthonormality of $(u_k)$ gives
\begin{align}\label{eq:graph-stationary-orthogonality}
\E\bigl(u_k(X_v)u_l(X_{v'})\bigm|\cT\bigr)=N^{-1}\ind\{k=l\}\gl_k^{\dtr(v,v')}.
\end{align}
Thus the cross terms $k\neq l$ vanish when averaging over $\cZ_2$, and $\E(A_{2,k}^2\mid\cT)\le N^{-1}G_2(|\gl_k|)$. In the first tree, for $u,u'\in\cZ_1(s)$ with common ancestor $w=u\wedge u'$ at generation $a=\gen(w)$, conditional independence of the two branches from $w$ gives
\begin{align*}
\E\bigl(u_k(X_u)u_k(X_{u'})\bigm|\cT\bigr)=\gl_k^{\dtr(u,u')}(P^a u_k^2)(o_N).
\end{align*}
Since $P^a$ is a Markov operator, the spectral expansion and~\eqref{eq:graph-projector} give
\begin{align*}
\sum_{|\gl_k|\le\zeta}\gl_k^2(P^au_k^2)(o_N)&\le v_N^*,\\
\sum_{|\gl_k|>\zeta}(P^au_k^2)(o_N)&\le C\zeta^{-2r}/N.
\end{align*} Splitting the sum over $1\le k\le N-1$ at $|\gl_k|=\zeta$ therefore yields
\begin{align*}
\E\bigl(|R_N^P(s)|^2\bigm|\cT\bigr)
\le\frac{Cv_N^*}{N}G_1(\zeta)G_2(\zeta)+\frac{C}{N^2}G_1(\rho)G_2(\rho).
\end{align*}
Taking expectations over $\cT$ and applying Lemma~\ref{lem:cmj-facts}\ref{cmj-iv} proves~\eqref{eq:graph-bilinear}.
\end{proof}

For later use, define the cumulative hazard without spatial weights by
\begin{align}\label{eq:meanfield-hazard}
\overline\Lambda_N(s):=\frac1{N\Gamma(\theta+1)}\sum_{u\in\cZ_1(s)}\sum_{v\in\cZ_2(s)}\Phi_\theta\bigl(a_u(s),a_v(s)\bigr),
\end{align}
which is $\cT$-measurable. Write
\begin{align*}
s_N(x):=\tfrac12(\log N+x),\quad
L(x):=\theta M_\infty^{(1)}M_\infty^{(2)}e^x.
\end{align*}

\begin{corollary}\label{cor:hazard}
Suppose Assumptions~\ref{ass:weights} and~\ref{ass:graphs} hold. For every fixed $x\in\dR$,
\begin{align}\label{eq:hazard-limit}
\Lambda_{s_N(x)},\quad\Lambda^{\cZ}_{s_N(x)},\quad\overline\Lambda_N\bigl(s_N(x)\bigr)\to  L(x)\text{ in }L^1.
\end{align}
For every fixed $K<\infty$, the three convergences also hold uniformly in probability on $[-K,K]$. In particular, $\{\Lambda^{\cZ}_{s_N(K)}:N\text{ large}\}$ is uniformly integrable.
\end{corollary}

\begin{proof}
By Lemma~\ref{lem:transfer}, it suffices to consider the unpruned hazards $\Lambda^{\cZ}_{s_N(x)}$ and $\overline\Lambda_N(s_N(x))$. For any age function $F\ge0$, stationarity of the second tree and conditioning its root to differ from $o_N$ give
\begin{align}\label{eq:first-moment}
\E\biggl(\sum_{u\in\cZ_1(s)}\sum_{v\in\cZ_2(s)}P_N(X_u,X_v)F(a_u,a_v)\biggm|\cT\biggr)
\le\frac1{N-1}\sum_{u\in\cZ_1(s)}\sum_{v\in\cZ_2(s)}F(a_u,a_v).
\end{align}
For $A>0$, let $\chi_A$ be continuous, equal to $1$ on $[0,A]$, to $0$ on $[A+1,\infty)$, and linear on $[A,A+1]$, and set $\Phi_\theta^A(a,b):=\Phi_\theta(a,b)\chi_A(a)\chi_A(b)$. Since $\Phi_\theta-\Phi_\theta^A\le2^\theta(a^\theta\ind\{a>A\}+b^\theta\ind\{b>A\})$ and
\begin{align*}
\E\sum_{u\in\cZ(s)}a_u(s)^\theta\ind\{a_u(s)>A\}
\le\sum_{j\ge\lfloor A\rfloor}(j+1)^\theta\E Z_\theta(s-j)
\le Ce^s\sum_{j\ge\lfloor A\rfloor}(j+1)^\theta e^{-j},
\end{align*}
it follows from~\eqref{eq:first-moment} that truncating $\Phi_\theta$ to $\Phi_\theta^A$ in either unpruned hazard at $s=s_N(x)$ produces an $L^1$ error of at most $Ce^x\sum_{j\ge\lfloor A\rfloor}(j+1)^\theta e^{-j}$, which vanishes as $A\to\infty$ uniformly in $N$.

For fixed $A$, approximate $\Phi_\theta^A$ uniformly on $[0,\infty)^2$ within error $\eps>0$ by a finite sum $\sum_j\phi_j(a)\psi_j(b)$ of continuous functions supported in $[0,A+1]$. By~\eqref{eq:first-moment} with $F\equiv1$, this uniform approximation error contributes at most $C\eps e^x$ in $L^1$. For each product term, Proposition~\ref{prop:mixing} replaces $\sum_{u,v}P_N(X_u,X_v)\phi_j(a_u)\psi_j(a_v)$ in $L^2$ by
\begin{align*}
\frac1N\biggl(\sum_{u\in\cZ_1(s)}\phi_j\bigl(a_u(s)\bigr)\biggr)\biggl(\sum_{v\in\cZ_2(s)}\psi_j\bigl(a_v(s)\bigr)\biggr).
\end{align*}
By Lemma~\ref{lem:profile-clt} with $F\equiv1$ and the independence of $\cZ_1$ and $\cZ_2$, this product converges in $L^1$ at $s=s_N(x)$ to $e^xM_\infty^{(1)}M_\infty^{(2)}\E\phi_j(A_1)\E\psi_j(A_2)$, where $A_1,A_2\sim\Exp(1)$ are independent. Letting $\eps\to0$ and then $A\to\infty$, and using $\E\Phi_\theta(A_1,A_2)=\Gamma(\theta+2)-\Gamma(\theta+1)=\theta\Gamma(\theta+1)$, we obtain~\eqref{eq:hazard-limit}.

Finally, since the three hazards are nondecreasing in $x$ and the limit $L(x)$ is continuous, pointwise convergence on a finite partition of $[-K,K]$ yields uniform convergence in probability on $[-K,K]$, and $L^1$ convergence at $x=K$ yields uniform integrability of $\{\Lambda^{\cZ}_{s_N(K)}\}$.
\end{proof}

\section{Proof of Theorem~\ref{thm:main}}\label{sec:proof-main}

The cumulative-hazard limit in Corollary~\ref{cor:hazard} identifies the passage-time distribution. To obtain the joint hopcount limit and the law of the geodesic, we compare the first bridge in the internal process with a reference collision on the unmarked CMJ trees, separating the selection of the collision time and generations from the spatial labels.

We first define the reference collision law and state the path-comparison lemma. Sections~\ref{sec:proof-main-i} and~\ref{sec:path-comparison} then prove parts~\ref{main:a} and~\ref{main:b} of Theorem~\ref{thm:main}, respectively. The proof of the comparison lemma is deferred to Section~\ref{sec:proof-path-comparison}. Its central spatial step is the path identity~\eqref{eq:graph-path-identity}, where symmetry of $P_N$ reverses the second ancestral line and averaging over the uniform second root supplies the transition across the bridge.

Conditionally on $\cT$, assign to each pair $u\in\cZ_1$, $v\in\cZ_2$ an independent exponential clock in cumulative-hazard time, with cumulative hazard zero before $\gs_u\vee\gs_v$ and
\begin{align*}
\overline K_{u,v}(s):=
\frac{\Phi_\theta\bigl(a_u(s),a_v(s)\bigr)}{N\Gamma(\theta+1)}
\quad(s\ge\gs_u\vee\gs_v).
\end{align*}
These cumulative hazards sum to $\overline\Lambda_N(s)$ in~\eqref{eq:meanfield-hazard}. Since there are finitely many individuals at each finite time and $\overline K_{o_1,o_2}(s)\to\infty$ as $s\to\infty$, the reference collision time $\widetilde\tau_N$ and pair $(\widetilde u_N,\widetilde v_N)$ realizing the minimum are almost surely finite and unique. Put
\begin{align*}
\tT_N:=2\widetilde\tau_N-\log N,
\quad
\widetilde H_N:=\gen(\widetilde u_N)+\gen(\widetilde v_N)+1.
\end{align*}
Independently of the unmarked trees and the reference clocks, let $(\widetilde S_j)_{j\ge0}$ be a $P_N$-Markov chain started at $o_N$. The reference spatial path is $(\widetilde S_j)_{j=0}^{\widetilde H_N}$.

\begin{lemma}[Path comparison]\label{lem:path-comparison}
Under Assumptions~\ref{ass:weights} and~\ref{ass:graphs}, on the disjoint union $\bigcup_{h\ge1}(\dR\times\{h\}\times\cV_N^{h+1})$,
\begin{align}\label{eq:full-path-comparison}
d_{\mathrm{TV}}\Bigl(
\cL\bigl(\overline T_N,H_N,(x_j)_{j=0}^{H_N}\bigr),\,
\cL\bigl(\tT_N,\widetilde H_N,
(\widetilde S_j)_{j=0}^{\widetilde H_N}\bigr)
\Bigr)\to 0,
\end{align}
uniformly in $o_N\in\cV_N$ and over kernels satisfying~\eqref{eq:graph-assumptions} with fixed $r,\eta,C$.
\end{lemma}

\subsection{Proof of Theorem~\ref{thm:main}\ref{main:a}}
\label{sec:proof-main-i}

\begin{proof}
For the passage-time marginal, fix $x\in\dR$ and define $s=\frac12(\log N+x)$. By Proposition~\ref{prop:half} and the Cox representation~\eqref{eq:cox},
\begin{align*}
\pr\bigl(\fa_NT_N(o_N,U_N)-\log N>x\bigr)=\pr(\tcol>s)=\E\bigl(e^{-\Lambda_s}\bigr).
\end{align*}
By Corollary~\ref{cor:hazard}, $\Lambda_s\pto\theta M^{(1)}_\infty M^{(2)}_\infty e^{x}$. Since $z\mapsto e^{-z}$ is bounded and continuous on $[0,\infty)$, the right-hand side converges to $\E\exp(-\theta M^{(1)}_\infty M^{(2)}_\infty e^{x})$, establishing~\eqref{eq:main-tail}. For $G$ a standard Gumbel variable independent of $(M^{(1)}_\infty,M^{(2)}_\infty)$, we have
\begin{align*}
\pr(Y_\theta>x)=\pr\bigl(G<-x-\log(\theta M^{(1)}_\infty M^{(2)}_\infty)\bigr)=\E\exp\bigl(-\theta M^{(1)}_\infty M^{(2)}_\infty e^{x}\bigr).
\end{align*}
Since $M^{(i)}_\infty>0$ a.s., $Y_\theta$ is a.s.\ finite and the passage-time convergence follows.

\noindent\textbf{hopcount fluctuations.}
Let $(u^*,v^*)$ be the bridge pair in the internal process realizing $\tcol$. By Proposition~\ref{prop:half}, the geodesic from $o_N$ to $U_N$ is the concatenation of the ancestral line of $u^*$ in the first tree, the bridge $\{X_{u^*},X_{v^*}\}$, and the ancestral line of $v^*$ in the second tree, traversed backwards. Thus
\begin{align}\label{eq:H-decomp}
H_N=\gen(u^*)+\gen(v^*)+1,
\quad \overline T_N=2\tcol-\log N.
\end{align}
By Lemma~\ref{lem:path-comparison}, it suffices to establish the joint limit for the reference pair $(\tT_N,\widetilde H_N)$. Put $r_N:=\frac12\log N$ and, for $z\in\dR$, set
\begin{align*}
j_N(z):=\biggl\lfloor\frac{r_N}{\theta}
+\frac{\sqrt{r_N}}{\theta}z\biggr\rfloor.
\end{align*}
For fixed $z_1,z_2\in\dR$, define the generation-restricted cumulative hazard
\begin{align*}
\overline\Lambda_N^{z_1,z_2}(s)
:=\frac1{N\Gamma(\theta+1)}
\sum_{u\in\cZ_1(s)}\sum_{v\in\cZ_2(s)}
\Phi_\theta\bigl(a_u(s),a_v(s)\bigr)
\ind\{\gen(u)\le j_N(z_1),\ \gen(v)\le j_N(z_2)\}.
\end{align*}
Since the thresholds $j_N(z_i)$ do not depend on $s$, this function is nondecreasing in $s$, and $0\le\dd\overline\Lambda_N^{z_1,z_2}\le\dd\overline\Lambda_N$ as Stieltjes measures.

For $s=s_N(x)=r_N+x/2$, the shift $s-r_N=x/2$ and the integer part in $j_N(z_i)$ alter the standardized generation threshold by $o(1)$. Thus Lemma~\ref{lem:profile-clt} gives, for every continuous compactly supported $\phi$,
\begin{align*}
e^{-s}\sum_{u\in\cZ_i(s)}\phi\bigl(a_u(s)\bigr)
\ind\{\gen(u)\le j_N(z_i)\}
\to
M_\infty^{(i)}\Phi(z_i)\int_0^\infty\phi(a)e^{-a}\dd a
\text{ in }L^1.
\end{align*}
Applying the same age truncation and finite-product approximation as in Corollary~\ref{cor:hazard}, using the independence of $\cZ_1$ and $\cZ_2$, yields
\begin{align}\label{eq:marked-cumulative-limit}
\overline\Lambda_N^{z_1,z_2}\bigl(s_N(x)\bigr)
\to \Phi(z_1)\Phi(z_2)L(x)
\text{ in }L^1,
\quad
L(x)=\theta M_\infty^{(1)}M_\infty^{(2)}e^x.
\end{align}
Monotonicity and continuity of $L$ therefore yield uniform convergence in probability on compact $x$-intervals.

Set
\begin{align*}
A_{1,N}:=\frac{\theta\gen(\widetilde u_N)-r_N}{\sqrt{r_N}},
\quad
A_{2,N}:=\frac{\theta\gen(\widetilde v_N)-r_N}{\sqrt{r_N}}.
\end{align*}
For $K>\vert{}x\vert{}$, conditioning on $\cT$ gives
\begin{align*}
\pr\bigl(-K<\tT_N\le x,\ A_{1,N}\le z_1,\ A_{2,N}\le z_2\bigr)
=\E\int_{s_N(-K)}^{s_N(x)}
e^{-\overline\Lambda_N(s)}
\dd\overline\Lambda_N^{z_1,z_2}(s).
\end{align*}
Since $\overline\Lambda_N(s_N(\cdot))$ and $\overline\Lambda_N^{z_1,z_2}(s_N(\cdot))$ converge uniformly in probability on $[-K,K]$ to continuous limits, the Stieltjes integral converges in probability to
\begin{align*}
\Phi(z_1)\Phi(z_2)\int_{-K}^{x}e^{-L(y)}\dd L(y)
=\Phi(z_1)\Phi(z_2)\bigl(e^{-L(-K)}-e^{-L(x)}\bigr).
\end{align*}
Each integral lies in $[0,1]$ because $\dd\overline\Lambda_N^{z_1,z_2}\le\dd\overline\Lambda_N$, so bounded convergence gives convergence of the expectations. Letting $K\to\infty$ and using tightness of $\tT_N$ yields
\begin{align*}
(\tT_N,A_{1,N},A_{2,N})
\wto(Y_\theta,Z_1,Z_2),
\end{align*}
where $Z_1,Z_2\sim\normal(0,1)$ are independent of each other and of $Y_\theta$. Since
\begin{align*}
\frac{\widetilde H_N-\theta^{-1}\log N}
{\theta^{-1}\sqrt{\log N}}
=\frac{A_{1,N}+A_{2,N}}{\sqrt2}
+\frac{\theta}{\sqrt{\log N}},
\end{align*}
Lemma~\ref{lem:path-comparison} establishes~\eqref{eq:graph-joint-limit}, uniformly in $o_N\in\cV_N$ and over the kernel class in the theorem.
\end{proof}

\subsection{Proof of Theorem~\ref{thm:main}\ref{main:b}}\label{sec:path-comparison}

Lemma~\ref{lem:path-comparison} uses the reference collision time and hopcount $(\tT_N,\widetilde H_N)$. We first replace this pair by $(\overline T_N,H_N)$ via contraction of total variation under a Markov kernel, and then apply the mixing bound~\eqref{eq:graph-mixing} to obtain the empirical averages in~\eqref{eq:window-lln}.

\begin{proof}
Let $\mu_N=\cL(\overline T_N,H_N)$ and $\widetilde\mu_N=\cL(\tT_N,\widetilde H_N)$. Define a Markov kernel $\cK_N$ mapping $(t,h)$ to the law of $(t,h,(S_j)_{j=0}^{h})$, where $(S_j)_{j\ge0}$ is a $P_N$-Markov chain started at $o_N$. The reference law in Lemma~\ref{lem:path-comparison} is $\widetilde\mu_N\cK_N$, whereas $\mu_N\cK_N$ is the law of $(\overline T_N,H_N,(S_j)_{j=0}^{H_N})$ with $(S_j)_{j\ge0}$ independent of $(\overline T_N,H_N)$. Write $e_N$ for the total variation distance in~\eqref{eq:full-path-comparison}. Projecting onto the first two coordinates gives $d_{\mathrm{TV}}(\mu_N,\widetilde\mu_N)\le e_N$, so the triangle inequality and contraction under $\cK_N$ yield
\begin{align*}
d_{\mathrm{TV}}\Bigl(\cL(\overline T_N,H_N,\pi_N),\,\mu_N\cK_N\Bigr)\le e_N+d_{\mathrm{TV}}(\widetilde\mu_N\cK_N,\mu_N\cK_N)\le2e_N\to0.
\end{align*}
Disintegrating over the common $(\overline T_N,H_N)$-marginal identifies this distance with the expectation in~\eqref{eq:graph-path-tv}.

For~\eqref{eq:window-lln}, fix $k\ge0$ and set $M:=\sup_N\|\psi_N\|_\infty$ and $\overline\psi_N:=\E\psi_N(\widetilde S_0,\ldots,\widetilde S_k)$, where $(\widetilde S_j)_{j\ge0}$ is the stationary chain in the theorem. For the chain $(S_j)_{j\ge0}$ started at $o_N$, define
\begin{align*}
\varphi_N(z):=\E_z^{P_N}\psi_N(S_0,\ldots,S_k)-\overline\psi_N,\qquad Y_i:=\psi_N(S_i,\ldots,S_{i+k})-\overline\psi_N.
\end{align*}
Then $N^{-1}\sum_{z\in\cV_N}\varphi_N(z)=0$ and $\|\varphi_N\|_\infty\le2M$. For $j-i\ge k$, conditioning on $(S_0,\ldots,S_{i+k})$ and applying the Markov property gives
\begin{align*}
\E_{o_N}^{P_N}(Y_iY_j)=\E_{o_N}^{P_N}\bigl[Y_i(P_N^{j-i-k}\varphi_N)(S_{i+k})\bigr].
\end{align*}
By~\eqref{eq:graph-mixing}, $|\E_{o_N}^{P_N}(Y_iY_j)|\le4CM^2(1-\eta)^{j-i-k-2r}$ when $j-i\ge k+2r$, while $|\E_{o_N}^{P_N}(Y_iY_j)|\le4M^2$ for all $(i,j)$. Summing over $0\le i,j\le m-1$ gives, for every $m\ge1$,
\begin{align}\label{eq:window-second-moment}
\sup_{o_N\in\cV_N}\E_{o_N}^{P_N}\biggl[\biggl(\frac1m\sum_{i=0}^{m-1}Y_i\biggr)^2\biggr]\le\frac{C_kM^2}{m},
\end{align}
where $C_k$ depends only on $k,r,\eta,C$. Take $(S_j)_{j\ge0}$ independent of $(\overline T_N,H_N)$. Since $H_N\pto\infty$, conditioning on $H_N$ and applying~\eqref{eq:window-second-moment} with $m=H_N-k+1$ on $\{H_N\ge k\}$ establishes~\eqref{eq:window-lln} for $(S_j)_{j=0}^{H_N}$, and~\eqref{eq:graph-path-tv} transfers this convergence to $\pi_N$. All bounds hold uniformly as stated in the theorem.
\end{proof}

\begin{remark}
\label{rem:reference-simple}
The $P_N$-chain has no holding steps, but it may revisit an earlier label. For deterministic $h\ge1$, the Markov property and~\eqref{eq:graph-return} give, uniformly in $o_N\in\cV_N$,
\begin{align*}
\pr_{o_N}^{P_N}\bigl(S_i=S_j\text{ for some }0\le i<j\le h\bigr)
&\le\sum_{k=2}^{h}(h+1-k)\max_xP_N^k(x,x)
\le Ch^2/N+C v_N^*h.
\end{align*}
Fix $L>1/\theta$ and define $h_N=\lceil L\log N\rceil$. Since $(S_j)_{j\ge0}$ is independent of $H_N$ and $\pr(H_N>h_N)\to0$, it follows from~\eqref{eq:graph-small-masses} that
\begin{align*}
\pr\bigl(S_i=S_j\text{ for some }0\le i<j\le H_N\bigr)
\le\pr(H_N>h_N)+Ch_N^2/N+C v_N^*h_N\to 0.
\end{align*}
\end{remark}

\subsection{Proof of Lemma~\ref{lem:path-comparison}}
\label{sec:proof-path-comparison}
Fix $K<\infty$ and write $a=s_N(-K)$, $b=s_N(K)$. Let $F$ be a measurable function on the state space in Lemma~\ref{lem:path-comparison}, with $\vert{}F\vert{}\le1$. For each pair $(u,v)\in\cZ_1\times\cZ_2$, define $h(u,v):=\gen(u)+\gen(v)+1$. Let $\mvX(u,v)\in\cV_N^{h(u,v)+1}$ be the vertex path consisting of the ancestral line from $o_N$ to $X_u$, the bridge $(X_u,X_v)$, and the second ancestral line traversed from $X_v$ back to $U_N$. This path is defined on the unpruned labelled trees, where repeated labels are permitted. Set $q_{u,v}(s):=0$ for $s\le\gs_u\vee\gs_v$ and
\begin{align}\label{eq:selection}
q_{u,v}(s):=\frac{2}{\Gamma(\theta)}P_N(X_u,X_v)
\bigl(a_u(s)+a_v(s)\bigr)^{\theta-1}
\quad(s>\gs_u\vee\gs_v).
\end{align}

For a retained bridge, every ancestor of either endpoint is retained. Its ancestral paths and generations are therefore determined by $\cG$, so the mark below is $(s,\cG)$-measurable. Applying Lemma~\ref{lem:first-bridge} first to this forest mark, and then representing the retained forest under the CMJ coupling, the expectation of $F(\overline T_N,H_N,(x_j)_{j=0}^{H_N})\ind\{\tcol\in(a,b]\}$ equals
\begin{align*}
\E\int_a^b e^{-\Lambda_s}
\sum_{\substack{u\in\cZ_1(s)\setminus D_1(s)\\
v\in\cZ_2(s)\setminus D_2(s)}}
q_{u,v}(s)
F\bigl(2s-\log N,h(u,v),\mvX(u,v)\bigr)\dd s.
\end{align*}

The following extensions are comparisons of integrals on the auxiliary coupling, after the physical first-bridge identity has been applied. By Lemma~\ref{lem:transfer}, extending the sum to all pairs in $\cZ_1(s)\times\cZ_2(s)$ changes this integral by at most $\E(\Lambda_b^{\cZ}-\Lambda_b)=o(1)$ (pairs with $X_u=X_v$ carry zero weight since $P_N(x,x)=0$). Next, replacing $e^{-\Lambda_s}$ by $e^{-\overline\Lambda_N(s)}$ incurs an error bounded by
\begin{align*}
\E\biggl(
\sup_{a\le s\le b}
\bigl|e^{-\Lambda_s}-e^{-\overline\Lambda_N(s)}\bigr|
\Lambda_b^{\cZ}\biggr)=o(1),
\end{align*}
since the supremum is bounded by $1$ and tends to $0$ in probability while $\{\Lambda_b^{\cZ}\}$ is uniformly integrable by Corollary~\ref{cor:hazard}. Both errors are uniform in $\vert{}F\vert{}\le1$ and in $o_N$.

\noindent\textbf{The Markov path identity.}
Since the survival factor $e^{-\overline\Lambda_N(s)}$ is $\cT$-measurable, we may average over the spatial labels conditionally on $\cT$ before integrating over the collision time. Fix $\cT$ and a pair $(u,v)$, and write $k=\gen(u)$, $\ell=\gen(v)$ and $h=k+\ell+1$. For any sequence $(z_0,\ldots,z_h)\in\cV_N^{h+1}$ with $z_0=o_N$ and $z_h\ne o_N$, its conditional weight in $\E(P_N(X_u,X_v)f(\mvX(u,v))\mid\cT)$ is
\begin{align*}
\frac1{N-1}
\Biggl(\prod_{j=1}^{k}P_N(z_{j-1},z_j)\Biggr)
P_N(z_k,z_{k+1})
\Biggl(\prod_{j=k+2}^{h}P_N(z_j,z_{j-1})\Biggr)
=\frac1{N-1}\prod_{j=1}^{h}P_N(z_{j-1},z_j),
\end{align*}
where the three factors come from the first ancestral line, the bridge, and the second ancestral line, and the equality uses the symmetry of $P_N$ (with empty products equal to $1$, covering root endpoints). Labels off the two ancestral lines sum to $1$. Consequently, for every bounded function $f$ on $\cV_N^{h+1}$,
\begin{align}\label{eq:graph-path-identity}
\E\bigl(P_N(X_u,X_v)f(\mvX(u,v))
\bigm|\cT\bigr)
=\frac1{N-1}\E_{o_N}^{P_N}
\bigl(f(S_0,\ldots,S_h)\ind\{S_h\ne o_N\}\bigr).
\end{align}

For every $h\ge1$, stochasticity gives $P_N^h(o_N,o_N)\le m_N^*$. Applying~\eqref{eq:graph-path-identity} to $f(z)=F(2s-\log N,h,z)$, we may replace its right-hand side by $N^{-1}\E_{o_N}^{P_N}F(2s-\log N,h,(S_j)_{j=0}^h)$ with conditional error at most $(m_N^*+N^{-1})/(N-1)$. Multiplying by the age factors, summing over pairs, and integrating over $(a,b]$ bounds the total error by
\begin{align*}
C(m_N^*+N^{-1})
\E\int_a^b e^{-\overline\Lambda_N(s)}
\dd\overline\Lambda_N(s)
\le C(m_N^*+N^{-1})=o(1)
\end{align*}
by~\eqref{eq:graph-small-masses}. By the competing-clock formula for the reference process, the resulting integral equals
\begin{align*}
\E\bigl(
F\bigl(\tT_N,\widetilde H_N,
(\widetilde S_j)_{j=0}^{\widetilde H_N}\bigr)
\ind\{\widetilde\tau_N\in(a,b]\}\bigr).
\end{align*}

Finally, since $\pr(\overline T_N>x)=\E e^{-\Lambda_{s_N(x)}}$ and $\pr(\tT_N>x)=\E e^{-\overline\Lambda_N(s_N(x))}$, Corollary~\ref{cor:hazard} implies that both $\overline T_N$ and $\tT_N$ converge in distribution to $Y_\theta$ and are tight uniformly in $o_N$ and over the kernel class in the lemma. Taking the supremum over $\vert{}F\vert{}\le1$, and letting $N\to\infty$ and then $K\to\infty$, proves~\eqref{eq:full-path-comparison}.\qed

\section{Proof of Theorem~\ref{thm:edge-universality}}\label{sec:edge-proof}

We prove Theorem~\ref{thm:edge-universality} without modifying the spatial mixing estimates of Sections~\ref{sec:cmj} and~\ref{sec:mixing}, comparing the first collision in the general-weight model with a deterministic time change of the reference collision law from Section~\ref{sec:proof-main}. Two errors enter the comparison: rejected proposals remove their descendants, and retained bridge pairs have hazards that differ from the dominating Weibull hazards. We construct the coupling and state the hazard comparison lemma below. Section~\ref{sec:edge-main-argument} proves Theorem~\ref{thm:edge-universality}, Section~\ref{sec:first-order-univ} proves Proposition~\ref{prop:lower-tail-lln}, and Section~\ref{sec:edge-thinning} proves the comparison lemma.

For Theorem~\ref{thm:edge-universality} and Lemma~\ref{lem:edge-law-transfer}, we work in the normalization $c^{1/\theta}\fa_N$ and define
\begin{align*}
c_N:=1+\eps_N,\quad
\gl_N:=c_N^{1/\theta},\quad
\gd_N:=\frac{2\eps_N}{1+\eps_N}.
\end{align*}
We may assume $\eps_N<1$. First extend each cumulative hazard beyond $\ell_N:=2\log N$ by setting
\begin{align*}
\widehat{\cH}_{xy,N}(s)
:=\cH_{xy,N}(\ell_N)
+\frac{c_NP_N(x,y)}{\Gamma(\theta+1)}
  (s^\theta-\ell_N^\theta),\quad s>\ell_N,
\end{align*}
and leaving it unchanged on $[0,\ell_N]$. Coupling the original and extended weights through the same independent unit exponential variables makes their truncations at $\ell_N$ coincide: on the event that the extended model has rescaled passage time at most $\ell_N$, its passage time and geodesic agree with those of the original model. Until the end of the proof of Theorem~\ref{thm:edge-universality} we work with the extended hazards and omit the hat, so that~\eqref{eq:edge-hazard-approximation} holds for almost every $s>0$.

Let $\cZ_1^+,\cZ_2^+$ be independent CMJ trees with reproduction measure $c_N\mu_\theta$, obtained from the standard trees $\cZ_1,\cZ_2$ by dividing every birth time by $\gl_N$. Assign $P_N$-labels rooted at $o_N$ and $U_N$ as in Section~\ref{sec:cmj}. At a proposal from $x$ to $y$ at parental age $s$, accept the proposal with probability
\begin{align}\label{eq:edge-thinning}
p_{xy,N}(s):=
\frac{\Gamma(\theta)\cH_{xy,N}'(s)}
{c_NP_N(x,y)s^{\theta-1}},
\quad 1-\gd_N\le p_{xy,N}(s)\le 1,
\end{align}
using an independent $\mathrm{Unif}(0,1)$ mark for each proposal. Retain both roots, and retain a subsequent proposal if and only if its parent is retained, the proposal is accepted, and its target label is unoccupied.

By Poisson thinning, accepted proposals from $x$ to $y$ form a Poisson process with intensity $\cH_{xy,N}'(s)\dd s$, so their first points have the rescaled edge laws specified by $\cH_{xy,N}$. Active edges to unoccupied vertices retain their original clocks, while subsequent points for the same parent--target pair find the target occupied and are suppressed. As in Lemma~\ref{lem:embedding}, the retained individuals have the joint law of the two-root internal process, whereas unused proposals are auxiliary and residual bridge weights are sampled conditionally on the internal forest as in Proposition~\ref{prop:cox}.

Write $\gs_u^+$ for birth times and $a_u^+(s):=s-\gs_u^+$ for ages in the dominating trees, and let $\cB_i^F(s)$ denote the set of retained individuals in tree $i$ at time $s$. The internal-process cumulative hazard is
\begin{align}\label{eq:general-internal-hazard}
\Lambda_s^F:=\sum_{u\in\cB_1^F(s)}\sum_{v\in\cB_2^F(s)}
\bigl[\cH_{X_uX_v,N}(a_u^+(s)+a_v^+(s))
-\cH_{X_uX_v,N}(|a_u^+(s)-a_v^+(s)|)\bigr].
\end{align}
The deferred-decision and marked first-bridge formulas of Section~\ref{sec:exploration} apply to these hazards; each retained pair $(u,v)$ has pair hazard $K_{u,v}^F$ with Stieltjes measure
\begin{align*}
\dd K_{u,v}^F(s)
=2\cH_{X_uX_v,N}'\bigl(a_u^+(s)+a_v^+(s)\bigr)\dd s
\quad(s>\gs_u^+\vee\gs_v^+).
\end{align*}

Define the unpruned dominating hazard and its spatial average by
\begin{align}\label{eq:dominating-hazards}
\Lambda_s^{\cZ,+}
&:=\frac{c_N}{\Gamma(\theta+1)}
\sum_{u\in\cZ_1^+(s)}\sum_{v\in\cZ_2^+(s)}
P_N(X_u,X_v)\Phi_\theta\bigl(a_u^+(s),a_v^+(s)\bigr),\notag\\
\overline\Lambda_N^+(s)
&:=\frac{c_N}{N\Gamma(\theta+1)}
\sum_{u\in\cZ_1^+(s)}\sum_{v\in\cZ_2^+(s)}
\Phi_\theta\bigl(a_u^+(s),a_v^+(s)\bigr).
\end{align}
Since $\gl_N^\theta=c_N$ and $a_u^+(s)=\gl_N^{-1}a_u(\gl_Ns)$, the homogeneity of $\Phi_\theta$ yields
\begin{align}\label{eq:hazard-time-change}
\Lambda_s^{\cZ,+}=\Lambda_{\gl_Ns}^{\cZ},\quad
\overline\Lambda_N^+(s)=\overline\Lambda_N(\gl_Ns).
\end{align}

\begin{lemma}[Hazard error under thinning]\label{lem:edge-law-transfer}
Let $v_N^*$ be as in~\eqref{eq:graph-small-masses}. With $t=\gl_Ns$,
\begin{align}\label{eq:edge-transfer}
\E(\Lambda_s^{\cZ,+}-\Lambda_s^F)
&\le C(1+t)^{2\theta+2}
\biggl(\frac{e^{3t}}{N^2}
+\frac{v_N^*e^{2t}}{N}\biggr)
+C\eps_N(1+t)\frac{e^{2t}}{N}.
\end{align}
Moreover, $\dd\Lambda^{\cZ,+}-\dd\Lambda^F$ is a nonnegative measure.
For a bounded pair mark $B$ and $0\le a\le b$, replacing the retained pair measures by the dominating pair measures in a first-bridge integral over $(a,b]$, with survival factor $e^{-\Lambda_s^F}$ unchanged, changes the integral by at most $\|B\|_\infty\E(\Lambda_b^{\cZ,+}-\Lambda_b^F)$.
\end{lemma}

\subsection{Proof of Theorem~\ref{thm:edge-universality}}
\label{sec:edge-main-argument}
For the passage-time limit, define $s_N(x)=\frac12(\log N+x)$. Since $\eps_N\log N\to0$, we have
\begin{align*}
2\gl_Ns_N(x)-\log N
=\gl_Nx+(\gl_N-1)\log N\to  x.
\end{align*}
By Corollary~\ref{cor:hazard}, monotonicity, and~\eqref{eq:hazard-time-change},
\begin{align*}
\Lambda_{s_N(x)}^{\cZ,+},\quad
\overline\Lambda_N^+(s_N(x))
\to  L(x):=\theta M_\infty^{(1)}M_\infty^{(2)}e^x
\text{ in }L^1.
\end{align*}
At the collision scale, Lemma~\ref{lem:edge-law-transfer} gives
\begin{align*}
\E\bigl(\Lambda_{s_N(x)}^{\cZ,+}-\Lambda_{s_N(x)}^F\bigr)
\le C_x(\log N)^{2\theta+2}\bigl(N^{-1/2}+v_N^*\bigr)
+C_x\eps_N\log N=o(1),
\end{align*}
by~\eqref{eq:graph-small-masses} and $\eps_N\log N\to0$. Thus $\Lambda_{s_N(x)}^F$ has the same $L^1$ limit. Monotonicity and continuity of $L$ then yield uniform convergence in probability on compact $x$-intervals. The Cox representation yields~\eqref{eq:main-tail}, and $\{\Lambda_{s_N(K)}^{\cZ,+}\}$ is uniformly integrable for each fixed $K$.

\medskip\noindent\textbf{Transfer of the path law.}
To compare the marked first-bridge integrals, fix $K<\infty$ and restrict $s$ to $(s_N(-K),s_N(K)]$. For a retained bridge, its two ancestral lines and generations are determined by the internal forest. As in Section~\ref{sec:proof-path-comparison}, apply the first-bridge formula to this forest-measurable mark before extending the resulting integral to dominating pairs. By Lemma~\ref{lem:edge-law-transfer}, replacing the retained pair measures by the unpruned dominating measures changes the integral by $o(1)$ uniformly over marks bounded by $1$. Replacing the survival factor $e^{-\Lambda_s^F}$ by $e^{-\overline\Lambda_N^+(s)}$ incurs an error of at most
\begin{align*}
\E\biggl(
\sup_{|x|\le K}
\bigl|e^{-\Lambda_{s_N(x)}^F}
-e^{-\overline\Lambda_N^+(s_N(x))}\bigr|
\Lambda_{s_N(K)}^{\cZ,+}
\biggr)=o(1).
\end{align*}
After these replacements, the integrand is free of thinning marks and the survival factor is $\cT$-measurable, so the path identity~\eqref{eq:graph-path-identity} applies and removing the endpoint restriction costs $O(m_N^*+N^{-1})$ as in Lemma~\ref{lem:path-comparison}.

Let $(\tT_N,\widetilde H_N)$ be the unmarked reference collision pair of Section~\ref{sec:proof-main}. By~\eqref{eq:hazard-time-change}, the dominating reference collision time is $\widetilde\tau_N/\gl_N$, with the same selected pair and hopcount $\widetilde H_N$. Setting
\begin{align*}
\tT_N^+:=
\frac{\tT_N+\log N}{\gl_N}-\log N
\end{align*}
and taking an independent $P_N$-Markov chain $(\widetilde S_j)_{j\ge0}$ started at $o_N$, the integral comparison above followed by $N\to\infty$ and $K\to\infty$ yields
\begin{align}\label{eq:edge-reference-tv}
d_{\mathrm{TV}}\Bigl(
\cL\bigl(\overline T_N,H_N,(x_j)_{j=0}^{H_N}\bigr),\,
\cL\bigl(\tT_N^+,\widetilde H_N,
(\widetilde S_j)_{j=0}^{\widetilde H_N}\bigr)
\Bigr)\to 0.
\end{align}
Since $\tT_N^+-\tT_N=(\gl_N^{-1}-1)(\tT_N+\log N)\pto0$, the joint limit for $(\tT_N,\widetilde H_N)$ established in Section~\ref{sec:proof-main-i} together with~\eqref{eq:edge-reference-tv} gives~\eqref{eq:graph-joint-limit}. Projecting onto the first two coordinates, applying the Markov-path kernel $\cK_N$ as in Section~\ref{sec:path-comparison} to replace $(\tT_N^+,\widetilde H_N)$ in~\eqref{eq:edge-reference-tv} by $(\overline T_N,H_N)$, and disintegrating over the common $(\overline T_N,H_N)$-marginal establish~\eqref{eq:graph-path-tv}. Combining~\eqref{eq:graph-path-tv} with~\eqref{eq:window-second-moment} and $H_N\pto\infty$ then yields~\eqref{eq:window-lln}.

Finally, in the extended model $c^{1/\theta}\fa_NT_N=\log N+O_{\pr}(1)$, so $\pr(c^{1/\theta}\fa_NT_N>2\log N)\to0$, and the truncation coupling at $\ell_N=2\log N$ transfers the conclusions to the original edge weights. All bounds hold uniformly in $o_N\in\cV_N$ and over the classes specified in the theorem.\qed

\subsection{Proof of Proposition~\ref{prop:lower-tail-lln}}
\label{sec:first-order-univ}

The first-order conclusions follow from distribution-function bounds. A local quantile coupling compares passage times with their Weibull counterparts, while a first-moment estimate excludes paths with too few or too many edges.

\begin{proof}
Put $L=\log N$ and $m_N^*=\max_{x,y}P_N(x,y)$. Take independent $\mathrm{Unif}(0,1)$ variables $V_e$ and couple $\xi_e=(-\log(1-V_e))^{1/\theta}$ with $\go_e=F_e^{-1}(V_e)$, where $F_e^{-1}(u)=\inf\{t:F_e(t)\ge u\}$. Let $T_N^{\mathrm W}$ be the passage time from $o_N$ to $U_N$ under the Weibull weights $g_N(e)^{-1/\theta}\xi_e$.

First assume~\eqref{eq:cdf-two-sided}. For $0<\gd<c_1$, comparing $F_e(t)$ with $1-e^{-(c_1-\gd)t^\theta}$ and $1-e^{-(c_2+\gd)t^\theta}$ gives
\begin{align*}
(c_2+\gd)^{-1/\theta}\xi_e
\le \go_e\le(c_1-\gd)^{-1/\theta}\xi_e
\end{align*}
whenever $\xi_e$ is sufficiently small, uniformly in $N$ and $e\in\cE_N$; similarly, $\xi_e\le(c_2+\gd)^{1/\theta}\go_e$ whenever $\go_e$ is sufficiently small. By Theorem~\ref{thm:main}, $\fa_NT_N^{\mathrm W}\le2L$ with probability $1-o(1)$. On this event, every edge noise along the Weibull geodesic is at most
\begin{align*}
2\biggl(\frac{m_N^*}{\Gamma(\theta+1)}\biggr)^{1/\theta}L=o(1)
\end{align*}
by~\eqref{eq:intro-small-masses}, which yields $T_N(o_N,U_N)\le(c_1-\gd)^{-1/\theta}T_N^{\mathrm W}$. This upper bound implies that every edge noise along any geodesic in the $\go$-environment is $O((m_N^*)^{1/\theta}L)=o(1)$, so the reverse inequality gives
\begin{align}\label{eq:local-quantile-sandwich}
(c_2+\gd)^{-1/\theta}T_N^{\mathrm W}
\le T_N(o_N,U_N)
\le(c_1-\gd)^{-1/\theta}T_N^{\mathrm W}
\end{align}
with probability $1-o(1)$. Combining~\eqref{eq:local-quantile-sandwich} with $\fa_NT_N^{\mathrm W}/L\pto1$ and letting $\gd\downarrow0$ proves~\eqref{eq:cdf-time-bounds}. Under~\eqref{eq:cdf-equivalence}, applying~\eqref{eq:local-quantile-sandwich} with $c_1=c(1-\gd)$ and $c_2=c(1+\gd)$ and letting $\gd\downarrow0$ gives the passage-time convergence in~\eqref{eq:cdf-joint-lln}.

\medskip\noindent\textbf{Controlling the hopcount.}
For the hopcount in part~\textup{(i)}, define $X_e=c^{1/\theta}\fa_NW_e$. For fixed $K,\gd>0$, uniformly for $0\le s\le KL$ and all large $N$,~\eqref{eq:cdf-equivalence} gives
\begin{align*}
\pr(X_e\le s)\le
\frac{(1+\gd)P_N(e)s^\theta}{\Gamma(\theta+1)}.
\end{align*}
For $\gl>0$ and $t\le KL$, Stieltjes integration by parts yields
\begin{align*}
\E\bigl(e^{-\gl X_e}\ind\{X_e\le t\}\bigr)
\le(1+\gd)P_N(e)\gl^{-\theta}.
\end{align*}
Let $\cN_{N,k}(t)$ denote the number of self-avoiding $k$-edge paths from $o_N$ to $U_N$ of total $X$-weight at most $t$. Independence of edge weights along each path, averaging over $U_N$, and stochasticity of $P_N$ give
\begin{align}\label{eq:short-path-first-moment}
\E\cN_{N,k}(t)
\le\frac{e^{\gl t}}{N-1}
   \bigl((1+\gd)\gl^{-\theta}\bigr)^k.
\end{align}
Fix $0<\eps<1/\theta$, take $t=(1+\gd)L$, and set $\gl_\pm=1\pm\theta\eps$. Summing~\eqref{eq:short-path-first-moment} over $k\ge(\theta^{-1}+\eps)L$ with $\gl=\gl_+$ and over $1\le k\le(\theta^{-1}-\eps)L$ with $\gl=\gl_-$ gives geometric series whose leading exponents of $N$ at $\gd=0$ are
\begin{align*}
\gl_\pm-1-\gl_\pm\log\gl_\pm<0.
\end{align*}
Both sums therefore vanish as $N\to\infty$ for $\gd>0$ sufficiently small, which together with $$c^{1/\theta}\fa_NT_N(o_N,U_N)/L\pto1$$ proves the hopcount convergence in~\eqref{eq:cdf-joint-lln}.

For the hopcount bounds in part~\textup{(ii)}, set $X_e=\fa_NW_e$ and replace $1+\gd$ in~\eqref{eq:short-path-first-moment} by $c_2$. Fix $K>c_1^{-1/\theta}$, so that $\fa_NT_N(o_N,U_N)\le KL$ with probability $1-o(1)$, and choose $\gl_+>c_2^{1/\theta}$ and $0<\gl_-<\min\{K^{-1},c_2^{1/\theta}\}$. For $C<\infty$ sufficiently large, the geometric sums of~\eqref{eq:short-path-first-moment} at $t=KL$ over $k\ge CL$ (with $\gl_+$) and over $1\le k\le L/C$ (with $\gl_-$) both tend to zero, proving~\eqref{eq:cdf-hop-bounds}. All estimates hold uniformly in $o_N\in\cV_N$ and apply to every geodesic.
\end{proof}

\subsection{Proof of Lemma~\ref{lem:edge-law-transfer}}
\label{sec:edge-thinning}
Call a proposal marked if its uniform mark exceeds $1-\gd_N$, which occurs independently of the unmarked trees and spatial labels. By~\eqref{eq:edge-thinning}, every rejected proposal is marked. A vertex of generation $k$ has a marked proposal along its ancestral line with probability at most $\gd_Nk$.

For a standard CMJ tree $\cZ$ with $m(t)=\E Z_\theta(t)$, the convolution formula $\mu_\theta^{*j}(\dd r)=r^{j\theta-1}\dd r/\Gamma(j\theta)$ and the beta integral give $\E\sum_{\gen(u)=j}a_u(t)^\theta=\Gamma(\theta+1)Q_{j+1}(t)$ for every $j\ge0$. Together with $Q_j(t)=t^{j\theta}/\Gamma(j\theta+1)$ and the birth-intensity bound from Lemma~\ref{lem:cmj-facts}, this yields
\begin{align*}
\E\sum_{u\in\cZ(t)}\gen(u)
&=\sum_{j\ge1}jQ_j(t)
=\frac{t}{\theta}r_\theta(t)
\le C(1+t)e^t,\\
\E\sum_{u\in\cZ(t)}a_u(t)^\theta
&=\Gamma(\theta+1)\bigl(m(t)-1\bigr)\le Ce^t,\\
\E\sum_{u\in\cZ(t)}\gen(u)a_u(t)^\theta
&=\Gamma(\theta+1)\sum_{j\ge2}(j-1)Q_j(t)
\le C(1+t)e^t.
\end{align*}
For each unpruned pair $(u,v)$, conditional spatial averaging gives $\E(P_N(X_u,X_v)\mid\cT)\le 1/(N-1)$ as in~\eqref{eq:first-moment}. By the marking bound, independence of $\cZ_1$ and $\cZ_2$, $\Phi_\theta(a,b)\le C(a^\theta+b^\theta)$, and~\eqref{eq:hazard-time-change}, the expected contribution to $\Lambda_s^{\cZ,+}$ from pairs with a marked vertex on either ancestral line is at most
\begin{align*}
\frac{C\gd_N}{N-1}
\E\sum_{u\in\cZ_1(t)}\sum_{v\in\cZ_2(t)}
\bigl(\gen(u)+\gen(v)\bigr)\bigl(a_u(t)^\theta+a_v(t)^\theta\bigr)
\le C\eps_N(1+t)\frac{e^{2t}}{N}.
\end{align*}

If a deleted vertex has no marked vertex on its ancestral line, then some ancestor shares its label with an earlier retained vertex, and hence with an earlier unpruned vertex. By~\eqref{eq:hazard-time-change}, the bound of Lemma~\ref{lem:transfer} for ancestor label coincidences on the unpruned trees applies at time $t=\gl_Ns$, yielding the first term on the right-hand side of~\eqref{eq:edge-transfer}.

For each pair $(u,v)$, let $K_{u,v}^+$ denote its summand in $\Lambda^{\cZ,+}$. On retained pairs, since the lower threshold $\vert{}\gs_u^+-\gs_v^+\vert{}$ is constant after $\gs_u^+\vee\gs_v^+$,~\eqref{eq:edge-thinning} gives
\begin{align*}
0\le\dd K_{u,v}^+-\dd K_{u,v}^F
\le\gd_N\dd K_{u,v}^+.
\end{align*}
Summing over retained pairs and using $\E\Lambda_s^{\cZ,+}\le Ce^{2t}/N$ bounds the expected retained-pair discrepancy by $C\eps_Ne^{2t}/N$, proving~\eqref{eq:edge-transfer}.

Finally, deleted pairs contribute $\dd K_{u,v}^+\ge0$ and retained pairs contribute $\dd K_{u,v}^+-\dd K_{u,v}^F\ge0$. Since retention is permanent at birth, $\dd\Lambda^{\cZ,+}-\dd\Lambda^F$ is a nonnegative measure, and the marked integral bound follows from $e^{-\Lambda_s^F}\le1$.\qed

\section{Additive functionals of the geodesic}\label{sec:observables}

The path comparisons in Theorems~\ref{thm:main} and~\ref{thm:edge-universality} transfer bounded additive observables of the geodesic to an independent $P_N$-Markov chain stopped at the sampled hopcount. We allow the observable to depend on both endpoints of each directed step, covering vertex occupation counts, transition counts between subsets, and, on the torus, functions of the increments. The law of large numbers follows immediately from the case $k=1$ of~\eqref{eq:window-lln}; below we establish a joint central limit theorem and then derive the resulting occupation-measure limits and continuum path comparisons. We retain the finite-$N$ stationary mean in the central-limit centering, as convergence of the mean without a rate does not suffice on the $\sqrt{\log N}$ scale.

Work in the setting of Theorem~\ref{thm:main} or Theorem~\ref{thm:edge-universality}, with $\overline T_N$ centered using the normalizer of the corresponding theorem. Let $\psi_N:\cV_N^2\to\dR$ be deterministic with $\sup_N\|\psi_N\|_\infty<\infty$, and define
\begin{align*}
A_N(\psi):=\sum_{j=0}^{H_N-1}\psi_N(x_j,x_{j+1}),\qquad\overline\psi_N:=\frac1N\sum_{x,y\in\cV_N}P_N(x,y)\psi_N(x,y).
\end{align*}
Taking $k=1$ in~\eqref{eq:window-lln} gives
\begin{align}\label{eq:observable-lln}
\frac{A_N(\psi)}{H_N}-\overline\psi_N\pto0.
\end{align}
For the fluctuations, set
\begin{align*}
b_{\psi,N}(x):=\sum_{y\in\cV_N}P_N(x,y)\psi_N(x,y)-\overline\psi_N.
\end{align*}
Let $u_N$ be the unique mean-zero solution of the Poisson equation
\begin{align}\label{eq:observable-poisson}
(I-P_N)u_N=b_{\psi,N},\qquad\frac1N\sum_{x\in\cV_N}u_N(x)=0,
\end{align}
and set
\begin{align}\label{eq:observable-variance}
\gs_{\psi,N}^2:=\frac1N\sum_{x,y\in\cV_N}P_N(x,y)\bigl(\psi_N(x,y)-\overline\psi_N+u_N(y)-u_N(x)\bigr)^2.
\end{align}

\begin{proposition}[Additive observables]\label{prop:observables}
Under the hypotheses of Theorem~\ref{thm:main} or Theorem~\ref{thm:edge-universality}, if $\gs_{\psi,N}^2\to\gs_\psi^2\ge0$, then
\begin{align}\label{eq:observable-joint-clt}
\biggl(\overline T_N,\,\frac{H_N-\theta^{-1}\log N}{\theta^{-1}\sqrt{\log N}},\,\frac{A_N(\psi)-\overline\psi_NH_N}{\sqrt{H_N}}\biggr)\wto(Y_\theta,Z,\gs_\psi Z'),
\end{align}
where $Z,Z'\sim\normal(0,1)$ are independent of each other and of $Y_\theta$.
\end{proposition}

\begin{proof}
Since $N^{-1}\sum_{x\in\cV_N}b_{\psi,N}(x)=0$, the bound~\eqref{eq:graph-mixing} gives $\|P_N^jb_{\psi,N}\|_\infty\le C(1-\eta)^{j-2r}\|b_{\psi,N}\|_\infty$ for all $j\ge2r$. Hence the series $u_N:=\sum_{j\ge0}P_N^jb_{\psi,N}$ converges and solves~\eqref{eq:observable-poisson}, with $\sup_N\|u_N\|_\infty<\infty$. For a $P_N$-Markov chain $(S_j)_{j\ge0}$ started at any vertex, set
\begin{align*}
D_{j+1,N}^{\psi}:=\psi_N(S_j,S_{j+1})-\overline\psi_N+u_N(S_{j+1})-u_N(S_j).
\end{align*}
By~\eqref{eq:observable-poisson}, $(D_{j+1,N}^{\psi})_{j\ge0}$ is a uniformly bounded martingale difference sequence with respect to the natural filtration of $(S_j)$, and
\begin{align}\label{eq:observable-decomposition}
\sum_{j=0}^{h-1}\bigl(\psi_N(S_j,S_{j+1})-\overline\psi_N\bigr)=\sum_{j=0}^{h-1}D_{j+1,N}^{\psi}+u_N(S_0)-u_N(S_h).
\end{align}
Its conditional variance is $v_{\psi,N}(S_j)$, where
\begin{align*}
v_{\psi,N}(x):=\sum_{y\in\cV_N}P_N(x,y)\bigl(\psi_N(x,y)-\overline\psi_N+u_N(y)-u_N(x)\bigr)^2
\end{align*}
is uniformly bounded with stationary mean $\gs_{\psi,N}^2$. Applying~\eqref{eq:window-second-moment} with $k=0$ to $v_{\psi,N}$ gives, uniformly in the initial vertex $S_0$,
\begin{align*}
\E\biggl[\biggl(\frac1h\sum_{j=0}^{h-1}v_{\psi,N}(S_j)-\gs_{\psi,N}^2\biggr)^2\biggr]\le\frac{C}{h}.
\end{align*}
Since $\gs_{\psi,N}^2\to\gs_\psi^2$, it follows that for every deterministic sequence $h_N\to\infty$, the average conditional variance $h_N^{-1}\sum_{j=0}^{h_N-1}v_{\psi,N}(S_j)$ converges in $L^2$ to $\gs_\psi^2$, uniformly in $S_0$. Uniform boundedness of $D_{j+1,N}^{\psi}$ yields the Lindeberg condition, so the martingale central limit theorem~\cite[Chapter~3]{HH80} and~\eqref{eq:observable-decomposition} give
\begin{align*}
\frac1{\sqrt{h_N}}\sum_{j=0}^{h_N-1}\bigl(\psi_N(S_j,S_{j+1})-\overline\psi_N\bigr)\wto\gs_\psi Z',
\end{align*}
uniformly in the initial vertex.

The deterministic-length convergence is uniform in bounded-Lipschitz distance over $h\ge\lfloor\sqrt{\log N}\rfloor$; otherwise a contradicting sequence $h_N\to\infty$ could be selected. Now let $(S_j)_{j\ge0}$ be a $P_N$-Markov chain started at $o_N$, independent of $(\overline T_N,H_N)$. Since $\pr(H_N\ge\lfloor\sqrt{\log N}\rfloor)\to1$, conditioning on $(\overline T_N,H_N)$ and using~\eqref{eq:graph-joint-limit} gives~\eqref{eq:observable-joint-clt} for the reference path $(S_j)_{j=0}^{H_N}$, and~\eqref{eq:graph-path-tv} transfers the joint limit to $\pi_N$.
\end{proof}

\begin{remark}\label{rem:observable-centering}
If $\overline\psi_N\to m$, then centering $A_N(\psi)$ by $\theta^{-1}\overline\psi_N\log N$ and scaling by $\sqrt{\log N}$ yields the joint limit
\begin{align*}
\biggl(\overline T_N,\,
\frac{H_N-\theta^{-1}\log N}{\theta^{-1}\sqrt{\log N}},\,
\frac{A_N(\psi)-\theta^{-1}\overline\psi_N\log N}{\sqrt{\log N}}\biggr)
\wto\bigl(Y_\theta,Z,\theta^{-1}mZ+\theta^{-1/2}\gs_\psi Z'\bigr).
\end{align*}
Replacing $\overline\psi_N$ by $m$ in this centering requires $\sqrt{\log N}(\overline\psi_N-m)\to0$. On the torus, if $\psi_N(x,y)$ depends only on the increment $y-x$, then $b_{\psi,N}\equiv0$, $u_N\equiv0$, and $\gs_{\psi,N}^2$ equals the variance under $\nu_N$. For position-dependent observables, the potential $u_N$ in~\eqref{eq:observable-variance} accounts for temporal correlations along the chain; in particular, for a coboundary $\psi_N(x,y)=\chi_N(x)-\chi_N(y)$ with $\chi_N$ bounded, the sum telescopes and $\gs_{\psi,N}^2=0$.
\end{remark}

\begin{corollary}\label{cor:position-increment}
In the torus setting of Corollary~\ref{cor:torus-path}, with edge weights satisfying either Assumption~\ref{ass:weights} or the hypotheses of Theorem~\ref{thm:edge-universality}, suppose $\overline\nu_N\wto\overline\nu_\infty$. Then for every continuous $\gf:\dT^d\times\dT^d\to\dR$,
\begin{align*}
\frac1{H_N}\sum_{j=0}^{H_N-1}\gf(x_j/n,\Delta_{j+1}/n)
\pto\int_{\dT^d}\int_{\dT^d}\gf(x,z)\overline\nu_\infty(\dd z)\dd x.
\end{align*}
In particular, the empirical position measure $H_N^{-1}\sum_{j=0}^{H_N-1}\gd_{x_j/n}$ converges weakly in probability to Haar probability measure on $\dT^d$.
\end{corollary}

\begin{proof}
Apply~\eqref{eq:observable-lln} to $\psi_N(x,y)=\gf(x/n,(y-x)/n)$. By uniform continuity of $\gf$, Riemann-sum approximation in the first variable and weak convergence $\overline\nu_N\wto\overline\nu_\infty$ in the second give $\overline\psi_N\to\int_{\dT^d}\int_{\dT^d}\gf(x,z)\overline\nu_\infty(\dd z)\dd x$. Testing against a countable convergence-determining class of continuous functions on $\dT^d$ yields weak convergence of the empirical position measures.
\end{proof}

For each row of Corollary~\ref{cor:manifold-kernels}, let $k$ denote the Markov density in the table.

\begin{corollary}
\label{cor:manifold-path}
Assume~\eqref{eq:manifold-mesh} for any row of Corollary~\ref{cor:manifold-kernels}, with $\overline T_N$ centered using the corresponding normalizer, and let $\iota_N:M\to\{1,\ldots,N\}$ map each point to the index of the cell containing it, with $\iota_N(z_{N,i})=i$. Let $(\mathsf X_j)_{j\ge0}$ be a Markov chain on $M$ with transition kernel $k(x,y)\mu(\dd y)$ started at $\mathsf X_0=z_{N,o_N}$, independent of $(\overline T_N,H_N)$. Then
\begin{align*}
d_{\mathrm{TV}}\Bigl(
\cL\bigl(\overline T_N,H_N,(x_j)_{j=0}^{H_N}\bigr),\,
\cL\bigl(\overline T_N,H_N,(\iota_N(\mathsf X_j))_{j=0}^{H_N}\bigr)
\Bigr)\to 0.
\end{align*}
Consequently, for every continuous $\psi:M\times M\to\dR$,
\begin{align*}
\frac1{H_N}\sum_{j=1}^{H_N}
\psi(z_{N,x_{j-1}},z_{N,x_j})
\pto\int_M\int_M\psi(x,y)k(x,y)\,\mu(\dd y)\mu(\dd x).
\end{align*}
In particular, the empirical position measure $H_N^{-1}\sum_{j=0}^{H_N-1}\gd_{z_{N,x_j}}$ converges weakly in probability to $\mu$, and $z_{N,x_{i_N\wedge H_N}}\wto\mu$ for every deterministic sequence $i_N\to\infty$.
\end{corollary}

\begin{proof}
By Lemma~\ref{lem:point-discretization} (and the same shell bound for $\overline Q_N$ in the cell-averaged case, together with Proposition~\ref{prop:zero-diagonal}), we have
\begin{align*}
\sup_i\sup_{x\in A_{N,i}}\sum_{j=1}^N
\biggl|P_N(i,j)-\int_{A_{N,j}}k(x,y)\,\mu(\dd y)\biggr|
\le C\eps_N,
\end{align*}
where $\eps_N=h_N^{\min\{1,d-\beta\}}\log(2/h_N)$ for the two singular kernels and $\eps_N=h_N$ for the heat kernel. Thus, conditionally on $\mathsf X_m=x\in A_{N,i}$, the law of $\iota_N(\mathsf X_{m+1})$ differs in total variation from $P_N(i,\cdot)$ by at most $C\eps_N$. Maximal step-by-step coupling of a $P_N$-Markov chain $(S_j)_{j\ge0}$ started at $o_N$ with $(\mathsf X_j)_{j\ge0}$ therefore yields $\pr(S_j\ne\iota_N(\mathsf X_j)\text{ for some }0\le j\le m)\le Cm\eps_N$. Since $H_N=O_{\pr}(\log N)$ and $\eps_N\log N\to0$, combining this coupling with~\eqref{eq:graph-path-tv} gives the total variation path comparison.

Applying~\eqref{eq:window-lln} with $k=1$ to $\psi_N(i,j)=\psi(z_{N,i},z_{N,j})$ and writing $\omega_\psi$ for the modulus of continuity of $\psi$ on $M\times M$, the row bound above yields
\begin{align*}
\biggl|\frac1N\sum_{i,j}P_N(i,j)\psi_N(i,j)
-\int_M\int_M\psi(x,y)k(x,y)\,\mu(\dd y)\mu(\dd x)\biggr|
\le C\|\psi\|_\infty\eps_N+\omega_\psi(2Ch_N)\to0,
\end{align*}
which proves the law of large numbers for transitions and, via a countable convergence-determining class in $C(M)$, weak convergence of the empirical position measures. Finally,~\eqref{eq:graph-path-tv} and~\eqref{eq:graph-mixing} give $d_{\mathrm{TV}}(\cL(x_{i_N\wedge H_N}),\mathrm{Unif}_N)\to0$, and $N^{-1}\sum_{i=1}^N\gd_{z_{N,i}}\wto\mu$ by~\eqref{eq:manifold-mesh}, whence $z_{N,x_{i_N\wedge H_N}}\wto\mu$.
\end{proof}

The total variation comparison in Corollary~\ref{cor:manifold-path} is stated at the level of cell indices, since the discrete laws on $\{z_{N,1},\ldots,z_{N,N}\}$ are mutually singular with respect to $\mu$.

\bibliographystyle{abbrv}
\bibliography{lrfpp_torus}

\appendix
\section{Verification of the kernel criteria}\label{app:kernel-verifications}

We prove the criteria and applications stated in Section~\ref{sec:applications}, and record the discretization estimates used there and in Section~\ref{sec:observables}.

\subsection{Row bounds, comparison and symmetric normalization}\label{app:comparison-products}

\begin{proof}[Proof of Proposition~\ref{prop:lp-smoothing}]
Equip the vertex set with the uniform probability measure. By H\"older's inequality, symmetry of $K_N$, and Minkowski's inequality,
\begin{align*}
 \norm{P_N f}_\infty\le A\norm{f}_{p'},
 \quad
 \norm{P_N f}_p\le A\norm{f}_1,\quad p'=\frac{p}{p-1}=\frac{1}{\gd}.
\end{align*}
By Riesz--Thorin interpolation, $\norm{P_N f}_s\le A\norm{f}_t$ whenever $1\le t\le p'$ and $1/s=1/t-\gd$. Starting from $t_0=1$, $m-1$ iterations reach an exponent $t_{m-1}$ with $1/t_{m-1}=1-(m-1)\gd\le\gd=1/p'$, so $t_{m-1}\ge p'$. Applying the $L^{p'}\to L^\infty$ bound together with the inclusion $L^{t_{m-1}}\subseteq L^{p'}$ at the $m$-th step gives $\norm{P_N^m}_{1\to\infty}\le A^m$. Since $P_N$ has nonnegative entries, this operator norm equals $\max_{x,y}NP_N^m(x,y)$, and subsequent Markov steps preserve the upper bound. Bounding the sum in the hypothesis by a single term gives $K_N(x,y)\le A N^{1/p}$, which yields the entrywise bound.
\end{proof}

\begin{proof}[Proof of Proposition~\ref{prop:zero-diagonal}]
Write $D_N=\operatorname{diag}(Q_N(x,x))$ and $R_N=Q_N-D_N$. The spectral decomposition and Cauchy--Schwarz give $\|Q_N^j-\Pi_N\|_{\infty\to\infty}\le C(1-\eta)^{j-2r}$ for $j\ge2r$. Consequently, the series
\begin{align*}
(I+Q_N)^{-1}
=\tfrac12\Pi_N+\sum_{j\ge0}(-1)^j(Q_N^j-\Pi_N)
\end{align*}
converges in operator norm, with $\|(I+Q_N)^{-1}\|_{\infty\to\infty}\le M$ for a constant $M$ independent of $N$.

Writing $b_N=\mathbf1+z$, the unit row-sum condition for~\eqref{eq:balanced-kernel} takes the form
\begin{align*}
(I+Q_N)z
=D_N\mathbf1+2D_Nz-\operatorname{diag}(z)R_Nz.
\end{align*}
On the ball $\|z\|_\infty\le R:=2M\gs_N$, applying $(I+Q_N)^{-1}$ to the right-hand side defines a map with $\ell^\infty$ norm at most $M(\gs_N+2\gs_NR+R^2)$ and Lipschitz constant at most $M(2\gs_N+2R)$. Since $\gs_N\le\max_x Q_N^{2r}(x,x)^{1/(2r)}\to0$, this map is a contraction on the ball for all large $N$, and its unique fixed point yields $b_N>0$ satisfying~\eqref{eq:balanced-kernel}.

Finally, the symmetric Schur bound gives $\|P_N-Q_N\|_{2\to2}\le C_1\gs_N$, so the absolute spectral gap of $P_N$ is at least $\eta/2$ for large $N$. Since $P_N\le(1+C_1\gs_N)Q_N$ entrywise, we also obtain $P_N^{2r}(x,x)\le(1+C_1\gs_N)^{2r}Q_N^{2r}(x,x)\le C'/N$.
\end{proof}

\begin{proof}[Proof of Corollary~\ref{cor:auxiliary-kernel}]
For $x\ne y$ with $Q_N(x,y)>0$, the cumulative hazard $$\cH_{xy,N}(s)=-\log\pr(c^{1/\theta}\fa_N W_{x,y}>s)$$ satisfies
\begin{align*}
\cH_{xy,N}(s)=\cH_0\biggl(s\biggl(\frac{Q_N(x,y)}{c\Gamma(\theta+1)}\biggr)^{\!1/\theta}\biggr).
\end{align*}
By~\eqref{eq:noise-hazard-expansion} and~\eqref{eq:balanced-kernel}, for almost every $0<s\le2\log N$,
\begin{align*}
\frac{\Gamma(\theta)\cH_{xy,N}'(s)}{P_N(x,y)s^{\theta-1}}
=\frac{1+O\bigl(s^\kappa Q_N(x,y)^{\kappa/\theta}\bigr)}{b_N(x)b_N(y)}.
\end{align*}
Since $\gs_N\le\max_{x,y}Q_N(x,y)\le(C/N)^{1/(2r)}$, the relative error on $[0,2\log N]$ is bounded by
\begin{align*}
C\gs_N+C(\log N)^\kappa N^{-\kappa/(2r\theta)}=o(1/\log N),
\end{align*}
and Theorem~\ref{thm:edge-universality} applies.
\end{proof}

\begin{proof}[Proof of Proposition~\ref{prop:comparison-products}]
Summing over $y$ gives $c\le1$. For $c<1$, write $P_N=cQ_N+(1-c)R_N$, where $R_N=(1-c)^{-1}(P_N-cQ_N)$ is symmetric stochastic and thus has operator norm at most $1$ on the mean-zero subspace. The triangle inequality gives $\|P_N-\Pi_N\|_{2\to2}\le c(1-\eta)+(1-c)=1-c\eta$, and the case $c=1$ is immediate. Under the upper bound $P_N\le C_0Q_N$, expanding the $2r$-step transition probability over paths gives $P_N^{2r}(x,x)\le C_0^{2r}Q_N^{2r}(x,x)$.

For the tensor product, every nontrivial eigenvalue of $P_N$ is a product $\prod_{j=1}^\ell \gl_j$, where each $\gl_j$ is an eigenvalue of $P_{N_j}^{(j)}$ and at least one factor is nontrivial; since $|\gl_j|\le 1$ for all $j$, the product has modulus at most $\max_j \|P_{N_j}^{(j)}-\Pi_{N_j}\|_{2\to 2}$. Finally, for any symmetric stochastic matrix, $r\mapsto P^{2r}(x,x)=\sum_i \gl_i^{2r}u_i(x)^2$ is nonincreasing in $r\ge 1$. Taking a common even exponent $2r$ at least as large as those of the $\ell$ factors yields $P_N^{2r}(x,x)\le \prod_{j=1}^\ell (C_j/N_j)=C/N$.
\end{proof}

\subsection{Torus kernels and normalizers}\label{app:torus-verifications}

\begin{proof}[Proof of Proposition~\ref{prop:torus-kernel}]
Summing $k(z/n)\le C(|z|/n)^{-\beta}$ over shells shows that the region $0<|z|\le \gd n$ contributes $O(\gd^{d-\beta})$ to $N^{-1}\sum_{z\ne 0}k(z/n)$ uniformly in $n$, while on $|z|>\gd n$ the sum is a Riemann sum for $k$. Letting $n\to\infty$ and then $\gd\downarrow 0$ gives $N^{-1}\sum_{z\ne0}k(z/n)\to I$. Choosing $p>1$ with $\beta p<d$, the same shell bound shows that $N^{-1}\sum_{z\ne0} k(z/n)^p$ is uniformly bounded, so Proposition~\ref{prop:lp-smoothing} yields the fixed-step return bound.

To verify the uniform spectral gap, let $\eps>0$. By truncating $k$ near $0$, approximating the bounded Riemann-integrable remainder in $L^1(\dT^d)$ by a continuous function, and then approximating uniformly by a trigonometric polynomial $F(x)=\sum_{|m|_\infty\le M}c_m e^{-2\pi i m\cdot x}$, we obtain
\begin{align*}
 \limsup_{n\to\infty}\frac1N\sum_{z\ne0}|k(z/n)-F(z/n)|<\eps.
\end{align*}
Represent frequencies $\xi\in\Tor$ in $(-n/2,n/2]^d$. For $n>2M$ and any $\xi\in\Tor$ with $|\xi|_\infty>M$, discrete orthogonality on $\Tor$ gives $\sum_{z\in\Tor}F(z/n)e^{2\pi i \xi\cdot z/n}=0$, so $\bigl|N^{-1}\sum_{z\ne 0}F(z/n)e^{2\pi i \xi\cdot z/n}\bigr|=|F(0)|/N=O(N^{-1})$. Choosing $\eps<I/4$ therefore ensures $\sup_{|\xi|_\infty>M}|\widehat\nu_n(\xi)|\le 1/2$ for all large $n$. For each of the finitely many nonzero frequencies with $|\xi|_\infty\le M$, $\widehat\nu_n(\xi)$ converges to $I^{-1}\int_{\dT^d}k(x)\cos(2\pi \xi\cdot x)\,\dd x\in(-1,1)$, which proves the uniform absolute spectral gap.
\end{proof}

\begin{proof}[Proof of Corollary~\ref{cor:torus-path}]
Since a random walk on $\Tor$ with transition matrix $P_N$ has i.i.d.~increments with law $\nu_N$, mapping the vertex path in~\eqref{eq:graph-path-tv} to its successive differences gives~\eqref{eq:path-tv}. Conditioning on $H_N$ bounds the conditional variance of $H_N^{-1}\sum_{j=1}^{H_N}\gf(\xi_j/n)$ by $\|\gf\|_\infty^2/H_N$. Since $H_N\pto\infty$, Chebyshev's inequality and~\eqref{eq:path-tv} yield the law of large numbers for each continuous $\gf$, and testing against a countable dense subset of $C(\dT^d)$ gives weak convergence of the empirical measures. Finally, by Lemma~\ref{lem:graph-spectral} (see~\eqref{eq:graph-mixing}), the $j$-step walk distribution lies within $C(1-\eta)^{j-2r}$ of $\mathrm{Unif}_N$ in total variation for all $j\ge 2r$. For any fixed $J\ge 2r$ and all large $N$ with $i_N\ge J$, combining this bound with~\eqref{eq:path-tv} yields $d_{\mathrm{TV}}(\cL(x_{i_N\wedge H_N}),\mathrm{Unif}_N)\le\pr(H_N<J)+C(1-\eta)^{J-2r}+o(1)$; letting $N\to\infty$ and then $J\to\infty$ completes the proof.
\end{proof}

\begin{proof}[Proof of Corollary~\ref{cor:powerlaw}]
Put $\beta=\ga\theta\in[0,d)$ and $g_n(z)=|z|^{-\beta}$ on $\Tor\setminus\{0\}$. Applying Proposition~\ref{prop:torus-kernel} to $k(x)=\|x\|^{-\beta}$ establishes Assumption~\ref{ass:graphs}, since the factor $n^{-\beta}$ cancels in $\nu_N(z)=g_n(z)/S_{N}$. In the Riemann-sum comparison for $n^{-(d-\beta)}S_N$, the cells meeting the origin contribute $O(n^{-(d-\beta)})$, and on the remaining cells the gradient bound $|\nabla(\|x\|^{-\beta})|\le C\|x\|^{-\beta-1}$ bounds the discretization error by
\begin{align*}
Cn^{-1}\int_{1/n}^{C_d}t^{d-\beta-2}\dd t
\le Cn^{-\min\{1,d-\beta\}}\log(2+n),
\end{align*}
so $S_N=c_\beta n^{d-\beta}(1+O(n^{-\min\{1,d-\beta\}}\log(2+n)))$. Truncating near the origin similarly yields $\overline\nu_N\wto c_\beta^{-1}\|x\|^{-\beta}\,\dd x$. Combining Theorem~\ref{thm:main} and Corollary~\ref{cor:torus-path}, and using Remark~\ref{rem:universality} to replace $\fa_N$ by $(\Gamma(\theta+1)c_{\ga\theta})^{1/\theta}n^{d/\theta-\ga}$ since the relative error in $S_N$ is $o(1/\log N)$, completes the proof.
\end{proof}

\begin{proof}[Proof of Corollary~\ref{cor:spreadout}]
For $g_n(z)=\ind\{0<|z|\le\ell_n\}$, Riemann sums give
\begin{align*}
S_N=|B_{\ell_n}|-1,\quad
\overline\nu_N\wto\mathrm{Unif}(\{x\in\dT^d:\|x\|\le\gl\}).
\end{align*}
The densities $N\nu_N$ are uniformly bounded, so Proposition~\ref{prop:lp-smoothing} gives the return bound with $r=1$. Moreover, $\nu_N$ differs in total variation by $o(1)$ from the kernel sampled from the indicator of the fixed ball $\{\|x\|\le\gl\}$. Proposition~\ref{prop:torus-kernel} gives a uniform spectral gap for the latter, which is preserved under $o(1)$ total variation perturbations. Theorem~\ref{thm:main} and Corollary~\ref{cor:torus-path} now apply. For $d=1$, $\ell_n<n/2$ for all large $n$, so $|B_{\ell_n}|-1=2\ell_n$.
\end{proof}

\subsection{Cell averages and independent samples}\label{app:iid-verification}

\begin{proof}[Proof of Proposition~\ref{prop:cell-averages}]
The matrix $Q_N$ is symmetric stochastic. On cell-constant functions it represents $E_NTE_N$, where $E_N$ is conditional expectation onto the partition. Since $E_N$ preserves constants and is an $L^2$ contraction, the absolute spectral gap is inherited. By Jensen's inequality,
\begin{align*}
\frac1N\sum_j(NQ_N(i,j))^p
\le N\int_{A_i}\int_XK(x,y)^p\,\mu(\dd y)\mu(\dd x)\le A^p.
\end{align*}
Proposition~\ref{prop:lp-smoothing} gives the return bound and the diagonal estimate; Proposition~\ref{prop:zero-diagonal} and Corollary~\ref{cor:auxiliary-kernel} then apply.
\end{proof}

\begin{lemma}[Near-constant row sums]\label{lem:near-regular-rows}
Let $A_N$ be a symmetric nonnegative matrix with zero diagonal, and write $r_N(x)=\sum_yA_N(x,y)$ and $s_N=\max_xr_N(x)$. Suppose that, for fixed $\gk>0$ and $p_0>1$,
\begin{align*}
\max_x|r_N(x)-1|\le C N^{-\gk},
\quad
\sup_x\frac1N\sum_y(NA_N(x,y))^{p_0}\le C.
\end{align*}
Define
\begin{align}\label{eq:auxiliary-point-kernel}
Q_N(x,y)=\frac{A_N(x,y)}{s_N}\quad(x\ne y),
\quad
Q_N(x,x)=1-\frac{r_N(x)}{s_N}.
\end{align}
If $Q_N$ has a uniform absolute spectral gap, it satisfies the hypotheses of Proposition~\ref{prop:zero-diagonal}, and the resulting reference kernel obeys $P_N(x,y)=(1+O(N^{-\gk}))A_N(x,y)$ for $x\ne y$. For physical intensities $g_N(x,y)=S_NA_N(x,y)$ and the common noise laws of Corollary~\ref{cor:edge-noise}, the conclusions of Corollary~\ref{cor:auxiliary-kernel} hold with normalizer $(c\Gamma(\theta+1)S_N)^{1/\theta}$.
\end{lemma}

\begin{proof}
The matrix $Q_N$ is symmetric stochastic, with $s_N=1+O(N^{-\gk})$ and $\gs_N=\max_xQ_N(x,x)=O(N^{-\gk})$. Choose $1<p\le p_0$ with $1-1/p\le\gk$. H\"older's inequality bounds the off-diagonal row $L^p$ norm by the $L^{p_0}$ hypothesis, and the diagonal contributes at most $CN^{p-1-\gk p}\le C$. Proposition~\ref{prop:lp-smoothing} gives the even-step return bound, and Proposition~\ref{prop:zero-diagonal} gives $\max_x|b_N(x)-1|\le C_1\gs_N=O(N^{-\gk})$, whence $P_N(x,y)=(1+O(N^{-\gk}))A_N(x,y)$ for $x\ne y$. Applying Corollary~\ref{cor:auxiliary-kernel} with scale $S_Ns_N$ and removing the factor $s_N^{1/\theta}=1+O(N^{-\gk})$ shifts the centered passage time by $O(N^{-\gk}\log N)=o(1)$, while applying the same affine map to both time coordinates preserves the total variation path comparison.
\end{proof}

\begin{proof}[Proof of Corollary~\ref{cor:iid-bounded}]
Conditionally on $X_i$, the variables $(k(X_i,X_j))_{j\ne i}$ are independent and bounded, with mean $\gk_0$. Hoeffding's inequality and a union bound give
\begin{align*}
\max_i\biggl|\frac1N\sum_{j\ne i}k(X_i,X_j)-\gk_0\biggr|
=O_{\pr}\bigl(\sqrt{\log N/N}\bigr).
\end{align*}
Put $A_N(i,j)=(N\gk_0)^{-1}k(X_i,X_j)$ for $i\ne j$ and $A_N(i,i)=0$, and let $Q_N$ be given by~\eqref{eq:auxiliary-point-kernel}. On an event $\Omega_N$ of probability $1-o(1)$, $\max_i|r_N(i)-1|\le C\sqrt{\log N/N}\le N^{-1/3}$ and $c_0/(2\gk_0N)\le Q_N(i,j)\le 2C_0/(\gk_0N)$ for $i\ne j$. Comparison with $(\mathbf J_N-I)/(N-1)$ via Proposition~\ref{prop:comparison-products} gives a uniform absolute spectral gap on $\Omega_N$. Define $P_N$ on $\Omega_N$ by Lemma~\ref{lem:near-regular-rows} (with $S_N=N\gk_0$), and on $\Omega_N^c$ take the simple random-walk kernel on $K_N$. Lemma~\ref{lem:near-regular-rows} then yields the conditional fluctuation and path limits on $\Omega_N$ with normalizer $(c\Gamma(\theta+1)N\gk_0)^{1/\theta}$.
\end{proof}

\subsection{Point discretizations}\label{sec:metric-kernels}\label{sec:point-samples}

\begin{lemma}[Point discretization]\label{lem:point-discretization}
Let $(X,\rho,\mu)$ be a compact metric probability space of diameter $L>0$ satisfying $\mu(B(x,t))\le Ct^D$ for $0<t\le L$ and some $D>0$. Let $K$ be symmetric and nonnegative, with $\int_XK(x,y)\,\mu(\dd y)=1$ for every $x$, and suppose its integral operator satisfies $\|T-\Pi\|_{L^2\to L^2}\le1-\eta$ for some fixed $\eta>0$. Assume that $K(x,y)\le C\rho(x,y)^{-\beta}$ for $x\ne y$ and some $0\le\beta<D$, and that
\begin{align}\label{eq:point-kernel-modulus}
|K(x,y)-K(x',y')|
\le C\bigl(\rho(x,x')+\rho(y,y')\bigr)
\rho(x,y)^{-\beta-1}
\end{align}
whenever $x\ne y$ and $\rho(x,x')+\rho(y,y')\le\rho(x,y)/2$. Let $(A_i)_{i=1}^N$ be a partition with $\mu(A_i)=1/N$, $\operatorname{diam}(A_i)\le Ch$, and representatives $x_i\in A_i$ satisfying $\rho(x_i,x_j)\ge ch$ for $i\ne j$, where $h=N^{-1/D}$. Put $A_N(i,j)=N^{-1}K(x_i,x_j)$ for $i\ne j$ and $A_N(i,i)=0$. Then
\begin{align}\label{eq:point-cell-comparison}
\max_i\sup_{x\in A_i}\sum_{j=1}^N
\biggl|A_N(i,j)-\int_{A_j}K(x,y)\,\mu(\dd y)\biggr|
\le C\eps_N,\quad
\eps_N=h^{\min\{1,D-\beta\}}\log(2/h).
\end{align}
The auxiliary matrix~\eqref{eq:auxiliary-point-kernel} satisfies the hypotheses of Proposition~\ref{prop:zero-diagonal}, and~\eqref{eq:point-cell-comparison} also holds with $A_N$ replaced by its balanced reference kernel $P_N$. If $K$ is bounded and Lipschitz on $X\times X$, the error $\eps_N$ improves to $h+N^{-1}$.
\end{lemma}

\begin{proof}
Choose $p_0>1$ with $\beta p_0<D$. Decomposing into dyadic annuli and using the volume bound gives $\sup_x\int_XK(x,y)^{p_0}\,\mu(\dd y)<\infty$, so the cell-averaged matrix $\overline Q_N$ satisfies Proposition~\ref{prop:cell-averages}. Equal cell measures and the diameter bound imply $\#\{j:\rho(x_i,x_j)\le t\}\le CNt^D$ for $t\ge h$. Fix $C_1>0$ large enough that $\operatorname{diam}(A_i)+\operatorname{diam}(A_j)\le\frac12\rho(x_i,x_j)$ whenever $\rho(x_i,x_j)>C_1h$. The cells with $\rho(x_i,x_j)\le C_1h$ contribute $O(h^{D-\beta})$ to~\eqref{eq:point-cell-comparison} uniformly over $x\in A_i$, using separation of representatives for $A_N(i,j)$ and the volume bound for the integrals. On the remaining cells,~\eqref{eq:point-kernel-modulus} bounds the summed discrepancy by
\begin{align*}
Ch\int_h^L t^{D-\beta-2}\,\dd t\le C\eps_N.
\end{align*}
When $K$ is bounded and Lipschitz, the off-diagonal discrepancy is $O(h)$ and omission of the diagonal cell $j=i$ costs $O(N^{-1})$.

Averaging~\eqref{eq:point-cell-comparison} over $x\in A_i$ gives $\max_i\sum_j|A_N(i,j)-\overline Q_N(i,j)|\le C\eps_N$. Thus the row sums of $A_N$ are $1+O(\eps_N)$, and the symmetric Schur bound gives $\|Q_N-\overline Q_N\|_{2\to2}=O(\eps_N)$, so the uniform spectral gap persists. The counting bound above also gives a uniform row $L^{p_0}$ bound for $NA_N$. Since $\eps_N=O(N^{-\gk})$ for some $\gk>0$, Lemma~\ref{lem:near-regular-rows} applies. Moreover, $s_N=1+O(\eps_N)$ and $\gs_N=O(\eps_N)$; the bound $\|b_N-1\|_\infty\le C\gs_N$ in Proposition~\ref{prop:zero-diagonal} therefore gives $P_N(i,j)=(1+O(\eps_N))A_N(i,j)$, whence $\max_i\sum_j|P_N(i,j)-A_N(i,j)|\le C\eps_N$ and~\eqref{eq:point-cell-comparison} holds for $P_N$.
\end{proof}

\paragraph{\bfseries Singular kernels at independent sample points.}
For i.i.d.~sample points $X_1,\ldots,X_N$ with law $\mu$ on an Ahlfors $D$-regular probability space, let $0<\beta<D$. Conditionally on $X_i=x$,
\begin{align*}
\E\biggl[\frac1N\sum_{j\ne i}\rho(X_i,X_j)^{-\beta}\biggm|X_i=x\biggr]
=\frac{N-1}{N}\int_X\rho(x,y)^{-\beta}\,\mu(\dd y)\le C,
\end{align*}
uniformly in $N$ and $x$, by Ahlfors regularity and $\beta<D$. Thus each fixed normalized row sum is uniformly tight. For $\beta>D/2$, define $r_N=N^{-(2-\eps)/D}$ with $\eps>0$. The expected number of pairs within distance $r_N$ is of order $N^\eps$; indicators of disjoint pairs are independent, and overlapping pairs contribute $O(N^3r_N^{2D})$ to the variance. By the second-moment method, a pair at distance at most $r_N$ exists with probability $1-o(1)$. For $\eps>0$ small enough that $(2-\eps)\beta>D$, such a pair contributes
\begin{align*}
 N^{-1}\rho(X_i,X_j)^{-\beta}
 \ge N^{(2-\eps)\beta/D-1}\to \infty
\end{align*}
to the corresponding normalized row sum, so $\max_i N^{-1}\sum_{j\ne i}\rho(X_i,X_j)^{-\beta}\pto\infty$ and Lemma~\ref{lem:near-regular-rows} does not apply.

\subsection{Manifold models}\label{app:manifold-verifications}

\begin{proof}[Proof of Corollary~\ref{cor:manifold-kernels}]
For the power-law row, isometry invariance and $\mu(B(x,t))\le Ct^d$ make $s_\beta$ finite and constant. The normalized density is bounded below by $\kappa_\beta$, so its operator has the table's absolute-gap bound. The triangle inequality and the mean-value bound for $t^{-\beta}$ give~\eqref{eq:point-kernel-modulus}.

The heat density $p_{t_0}$ is smooth, symmetric, strictly positive, and has unit row integrals. Its operator has nonconstant norm $e^{-t_0\gl_1}$~\cite{Dav89}, and its bounded Lipschitz regularity gives the improved point-discretization error.

For Bessel densities, Tonelli's theorem gives unit row integrals, and integration against Laplace eigenfunctions identifies $T_\beta=(I-\Delta)^{-a}$ and the stated nonconstant norm. Gaussian heat-kernel and derivative bounds yield~\eqref{eq:manifold-singular-bounds} and $|\nabla_xK_\beta(x,y)|\le C\rho(x,y)^{-\beta-1}$~\cite[Lemmas~2.5--2.6]{BCCGST14}, hence~\eqref{eq:point-kernel-modulus}.

The volume bound gives uniform row $L^p$ bounds for $1<p<d/\beta$ when $\beta>0$, and every finite $p>1$ for heat densities or $\beta=0$. Proposition~\ref{prop:cell-averages} treats cell averages. Lemmas~\ref{lem:point-discretization} and~\ref{lem:near-regular-rows} treat point samples, with scale $S_N=N\chi_N$ in either construction. Corollary~\ref{cor:auxiliary-kernel} then gives~\eqref{eq:manifold-centering} and the common-noise extension.
\end{proof}

\section{Moment convergence}\label{app:moments}

We prove the moment statement in Remark~\ref{rem:universality} under the Weibull assumptions. By Theorem~\ref{thm:main}, there is $b>1$ such that, for all large $N$ and every root $o$, a uniform independent target $V$ on $\cV_N$ satisfies
\begin{align*}
\pr(\fa_NT_N(o,V)>b\log N)\le1/4.
\end{align*}
Fix distinct vertices $o,x$. Take $2k$ independent copies of the edge environment and $k$ independent uniform vertices $V_j$. Route $j$ goes from $o$ to $V_j$ in copy $2j-1$ and from $V_j$ to $x$ in copy $2j$. Symmetry and a union bound show that its rescaled cost exceeds $2b\log N$ with probability at most $1/2$. These events are independent over $j$.

The edgewise minimum of the $2k$ environments, multiplied by $(2k)^{1/\theta}$, has the original Weibull law. Its passage time is at most the minimum route cost, after loop removal. Therefore
\begin{align*}
\sup_{o\ne x}\pr\bigl(\fa_NT_N(o,x)>2b(2k)^{1/\theta}\log N\bigr)\le2^{-k}.
\end{align*}
Choosing $k$ proportional to $z^\theta$ and adjusting constants for bounded $z$ yields
\begin{align}\label{eq:uniform-time-tail}
\sup_{o\ne x}\pr(\fa_NT_N(o,x)>z\log N)\le Ce^{-cz^\theta},\qquad z\ge1.
\end{align}

For a self-avoiding path $\pi=(x_0,\ldots,x_j)$, independence and the density bound from~\eqref{eq:rescaled-law} give
\begin{align*}
\pr\biggl(\sum_{e\in\pi}\fa_NW_e\le t\biggr)
\le\frac{t^{j\theta}}{\Gamma(j\theta+1)}\prod_{i=1}^jP_N(x_{i-1},x_i).
\end{align*}
Summing the product over self-avoiding $j$-step paths from a fixed root is at most one. For sufficiently small $c_0>0$, Stirling's formula consequently gives
\begin{align*}
\pr(H_N\ge h,\ \fa_NT_N(o_N,U_N)\le c_0h)
\le\sum_{j\ge h}\frac{(c_0h)^{j\theta}}{\Gamma(j\theta+1)}\le Ce^{-c_1h}.
\end{align*}
Together with~\eqref{eq:uniform-time-tail}, this proves uniform integrability of every positive power of $\fa_NT_N/\log N$ and $H_N/\log N$, uniformly in the root. Their laws of large numbers give the claimed moment limits.

These moment statements do not extend under lower-tail assumptions alone. For example, a continuous noise law that agrees with a Weibull law near zero but has survival asymptotic to $1/\log t$ at infinity satisfies Corollary~\ref{cor:edge-noise}. At each fixed $N$, the minimum of the finitely many incident weights has no finite positive moment, and neither does the passage time to a distinct target.

\end{document}